\documentclass[preprint,12pt]{elsarticle}
\usepackage{xcolor}
\usepackage{lineno,hyperref}
\usepackage{mathtools}
\usepackage{symbols}
\usepackage{amsfonts}
\usepackage{xspace}
\usepackage{empheq}
\usepackage{cancel}  
\usepackage{rotating} 
\usepackage{placeins}
\usepackage{tikz}
\biboptions{sort&compress}
\usepackage{pgfplots}

\modulolinenumbers[5]

\usepackage{ifpdf}

\ifpdf
\usepackage{pdfsync}
\fi
\usepackage{todonotes}

\usepackage[refpage]{nomencl}
\makenomenclature

\begin{document}

\begin{frontmatter}

\title{A geometry aware arbitrary order collocation Boundary Element Method solver for the potential flow past three dimensional lifting surfaces}

\author{Luca Cattarossi\fnref{myfootnote1}}
\author{Filippo Sacco\fnref{myfootnote4}}
\author{Nicola Giuliani\fnref{myfootnote2}}
\author{Nicola Parolini\fnref{myfootnote3}}
\author{Andrea Mola\fnref{myfootnote1}\corref{mycorrespondingauthor}}
\fntext[myfootnote1]{MUSAM Lab, Scuola IMT Alti Studi Lucca, Piazza S. Ponziano, 6 - 55100 Lucca, Italy, andrea.mola@imtlucca.it}
\fntext[myfootnote2]{previously Scuola Internazionale Sueriore di Studi Avanzati (SISSA), now Applied Materials Inc., Santa Clara, CA 95054-3299 USA, nicola\_giuliani@amat.com}
\fntext[myfootnote3]{MOX, Dipartimento di Matematica, Politecnico di Milano,  Piazza Leonardo da Vinci, 32 – 20133 Milano, Italy, nicola.parolini@polimi.it}
\fntext[myfootnote4]{Collegio San Carlo, Corso Magenta, 71, 20123 Milano, filippogd.sacco@gmail.com}


\cortext[mycorrespondingauthor]{Corresponding author}


\begin{abstract}

This work presents a numerical model for the simulation of potential flow past three dimensional lifting surfaces. The solver is
based on the collocation Boundary Element Method, combined with Galerkin variational formulation of the nonlinear Kutta condition imposed
at the trailing edge. A similar Galerkin variational formulation is also used for the computation of the fluid velocity at the wake collocation
points, required by the relaxation algorithm which aligns the wake with the local flow. The use of such a technique, typically associated
with the Finite Element Method, allows in fact for the evaluation of the solution derivatives in a way that is independent
of the local grid topology. As a result of this choice, combined with the direct interface with CAD surfaces, the solver is able to use
arbitrary order Lagrangian elements on automatically refined grids. Numerical results on a rectangular wing with NACA 0012 airfoil sections
are presented to compare the accuracy improvements obtained by grid spatial refinement or by discretization degree increase. Finally,   
numerical results on rectangular and swept wings with NACA 0012 airfoil section confirm that the model is able to reproduce
experimental data with good accuracy.

\end{abstract}


\end{frontmatter}



\section{Introduction and literature review}\label{sec:intro}

Accurate force prediction on three dimensional lifting bodies is paramount for the design of wings, tailplanes, rudders, aircraft and naval propellers,
hydrofoils and keels. The study of the forces and loads generated by such appendages through their interaction with the surrounding fluid is in
fact not only necessary for reliable equilibrium and stability estimates for cruising aircraft or vessels, but also for the the proper sizing of
their structure.  For most of the 20th century, potential flow theory has been the backbone of fluid dynamic computations in aerospace and naval
engineering. This is a result of the fact that --- when used on the correct problems --- potential flow solvers are capable of very accurate pressure force
estimates while requiring very limited computational resources. Despite the fact that in the last decades the constant increase in computational power made
affordable more accurate models, such as Navier--Stokes equations, potential flow theory is still widespread among aircraft and vessel designers. Simpler
models are in fact extremely valuable in the earlier stages of design, when both faster feedback is needed and not enough data is available to feed more
accurate solvers.

There are of course several different mathematical and numerical models based on potential flow theory that were developed to compute forces generated by
lifting bodies. For the most part, these models differ in the way the boundary conditions of the potential flow governing equation --- the Laplace equation ---
are handled. The lifting line \cite{PhillipsEtAl2001} and lifting surface methods \cite{LiuEtAl1996}, for instance, exploit slenderness of the lifting body
in one or two directions to linearize the
Laplace problem boundary conditions. Due to their extremely favorable balance between accuracy and computational cost these methods still enjoy significant
success in the engineering community. However, the linearized boundary conditions makes them unable to account for significantly thick surfaces, or for the
presence of arbitrarily shaped streamlined bodies surrounding the lifting bodies. The latter cases --- which include bulbs or fuselages --- can be instead
treated by models based on the exact formulation of the --- Neumann --- non penetration boundary condition for the velocity potential, and are in most cases 
discretized via the Boundary Element Method (BEM). Since the early works of Hess et al. \cite{HESS19671}, BEM discretization of the Laplace equation
for the flow velocity potential with exact non penetration boundary conditions represents a very interesting compromise, as it is able to both model
complex streamlined geometries and obtain fast predictions. 

This work describes the development of a numerical tool based on the BEM library $\pi$-BEM \cite{pi-BEM} for the computation of the flow past three
dimensional lifting bodies. As opposed to most BEM discretizations typically used for this application \cite{gaggeroBrizzolara2007}, which make use
of virtual singularity distribution formulations, $\pi$-BEM library implements an arbitrary order isoparametric discretization based on Lagrangian finite elements.
To allow for the automatic generation of the grids required by high order formulations, the library is also geometry aware, as described
in \cite{heltaiEtAlACMTOMS2021}. In the present context, this means that the simulations are directly interfaced with CAD data structures --- in the form
of IGES files --- so that the solver is able to generate on demand new grid points on the user-prescribed surfaces. The BEM solver has been shown
in \cite{pi-BEM} to have the expected convergence properties
on academic test cases characterized by --- regular enough --- analytic solutions. The potential flow on the lifting bodies application tackled in this work
poses instead several challenges, associated with the complex geometry features and lower solution regularity when compared to the academic test cases one. 
In addition, further complications had to be addressed to adapt the methods to the lifting body problem. Based on quasi-potential flow
model \cite{morinoGennaretti1992,GENNARETTI1996281}, the simulations have to account for the presence of the wake vortex sheet detaching from the wing
trailing edge. In this regard, among the many different forms of the Kutta condition \cite{XU1998415} used to determine the velocity potential jump on the wake
vortex sheet, we elected to employ the nonlinear form of the condition, imposing continuous pressure across the wing trailing edge. At the numerical level,
this condition involves the potential gradient evaluation at the trailing edge collocation points. In all BEM formulations based on standard Lagrangian finite
elements, such potential gradients are discontinuous on all the edges of the computational grid cells. For such a reason, BEM formulations typically used in this
application are based on piecewise constant discretizations, in which the solution gradient is computed making use of finite
differences \cite{gaggeroBrizzolara2007}. Clearly, such methodology results in a formulation heavily dependent on the specific topology and orientation of
the computational grid used. Moreover, of course, it limits the degree of the finite elements that can be employed. A possible interesting alternative
presented in \cite{CHOULIARAS2021113556}
is represented by the use of isogeometric formulation, in which smooth NURBS shape functions result in continuous gradients at the collocation points. However,
isogeometric formulation typically requires water tight CAD surfaces in most cases composed by NURBS patches all having similar size, which is are quite difficult
to obtain in
industrial applications. So, to make the present isoparametric formulation independent on the kind of mesh used, and suitable for arbitrary degree Lagrangian finite
elements, the potential gradients in the Kutta condition are here computed via Galerkin $L_2$ projection \cite{waveBem}. To the best of the authors' knowledge,
such a procedure --- borrowed from FEM --- has never been coupled to the BEM solution of quasi-potential flow simulations. In addition, the geometry
aware implementation allows for fully automated refinement of a coarse initial mesh on CAD surfaces that are not necessarily water tight. We refer the interested
reader to \cite{molaIsope2016,DASSI2014253,heltaiEtAlACMTOMS2021} for more detail on such a methodology. Again, this represents an aspect of novelty, as
no arbitrary order and CAD interfaced solver has been developed for this aeronautical and naval engineering application.    
   
A further extension to adapt $\pi$-BEM functions to the simulation of the flow past lifting surfaces is related to the wake
geometry relaxation. In the quasi-potential flow model \cite{morinoGennaretti1992,GENNARETTI1996281}, the wake vortex sheet surface must be
a flow surface everywhere parallel to the local velocity field. The wake relaxation algorithm used in this work to align the wake geometry to the local fluid
velocity requires the computation of the the potential gradient on all the wake collocation points. Since on the wake vortex sheet the only variable available
is the jump of the potential across such surface of discontinuity, the Galerkin $L_2$ projection method used to obtain the potential gradient on the body
cannot be used in this context.
As in other works \cite{bernasconiPhD,gaggeroBrizzolara2007} we then resorted to the integration of the hypersingular Boundary Integral Equation (BIE) on the wake
collocation points. Quadrature rules by Guiggiani et al. \cite{guiggianiEtAl1992} were used to evaluate the hypersingular integrals resulting from the gradient
BIE. However, such quadrature rules are based on the assumption that the potential gradient is H\"older continuous at each collocation point considered. This
of course poses a problem for finite elements which feature collocation points located on cell edges. In particular, this problem affects all Lagrangian
finite elements
of order higher than zero. Also in the present case, it is possible to abandon Lagrangian shape functions and resort to
NURBS shape functions as in \cite{bernasconiPhD} to build a single patch for the wake in which the potential gradient has the desired regularity. In the present
work, we elected instead to write a weak form of the hypersingular BIE and solve it based on a Galerkin method. As will be shown, at the numerical level this
will only require
the evaluation of the hypersingular BIE on the standard Gauss quadrature nodes. These are located inside the grid cells,
where the potential derivatives are H\"older continuous. Once again, the Galerkin approach results in a formulation that is independent of the grid or finite element
used, allowing for user prescribed discretization order selection at the start of each simulation.

Given the CAD data structures integration and use of Galerkin variational formulation for both the Kutta condition and wake velocity computation,
the BEM solver developed is able to make use of Lagrangian shape functions of arbitrary order on structured or unstructured meshes with possible
non conformal cells. As such, the present work shows that all the necessary ingredients have been put in place towards the implementation of
$hp$ refinement for the simulation of quasi-potential flow past lifting bodies. The derivation of adequate error estimators and the test of such a refinement
strategy will be object of future work.

The content of this paper is organized as follows. Section \ref{sec:math_model} describes the physical problem
and discusses the derivation of the mathematical model equations. Section \ref{sec:numerical_model} presents the main features of the BEM solver used,
reporting details of the Galerkin formulation for the Kutta condition and wake velocity computation. Section \ref{sec:results} presents and discusses
numerical results and comparison with experimental data. Finally, Section \ref{sec:conclusions} is devoted to some concluding remarks and future
perspectives. 

\section{Mathematical model equations}\label{sec:math_model}

The equations of motion that describe the velocity and pressure
fields $\vb$ and $p$ in a fluid region surrounding a moving body are the
incompressible Navier--Stokes equations. In the present context, such equations
are written in the moving domain $\Omega(t) \subset\mathbb{R}^3$
which represents the simply connected region of air and water surrounding
the lifting body --- or wing. Making use of the right handed orthogonal coordinate
system $\widehat{xyz}$, where $\hat{\ib},\hat{\jb},\hat{\kb}$ are the unit
vectors aligned with axes $x$, $y$ and $z$, respectively, the equations read:


\begin{subequations}
  \label{eq:incompressible-euler}
  \begin{alignat}{3}
    \label{eq:conservation-momentum}
    & \rho \left( \frac{\partial \vb}{\partial t} + \vb \cdot \nablab
    \vb \right)= \nablab\cdot\sigmab + \bb \qquad & \text{ in } \Omega(t)\\
    \label{eq:incompressibility}
    & \nablab \cdot \vb = 0 & \text{ in } \Omega(t)\\
    \label{eq:non-penetration-condition}
    & \vb = \vb_g & \text{ on } \Gamma(t)
  \end{alignat}
\end{subequations}
in which $\rho$ is the --- constant --- density of the fluid and $\bb$ are external volume
forces which in the present framework typically are related to gravity and possible inertial
forces due to possible acceleration of the reference frame. The term $\sigmab = -p\Ib +
\mu(\nablab\vb+\nablab\vb^T)$ is the stress tensor for an
incompressible Newtonian fluid, $\Gamma(t) := \partial\Omega(t)$ is
the boundary of the region of interest, and $\nb$ is the outer normal
to the boundary $\Gamma(t)$. On the boundaries of the domain, the prescribed
velocity $\vb_g$ is either equal to the body velocity, or to a given velocity field
(for surfaces far away from body).


Equation~(\ref{eq:conservation-momentum}) is usually referred to as
the momentum balance equation, while (\ref{eq:incompressibility}) is
referred to as the incompressibility constraint, or continuity
equation.

\subsection{Potential flow model}\label{sec:potential_flow}

In the flow field past a slender body --- such as a wing, as depicted in Figure~\ref{fig:wholeDomain} --- aligned with the
free stream velocity, vorticity is confined to the boundary layer region and to a thin
wake following the body. In such conditions, the assumption of irrotational flow
and non viscous fluid are fairly accurate in the vast majority of the computational domain.

If simply connected domains are considered, under the assumptions of irrotational flow and inviscid
fluid, at each point $\xb=x\hat{\ib}+y\hat{\jb}+z\hat{\kb}$ 
the velocity field $\vb(\xb,t)$ admits a representation through a scalar potential function
$\Phi(\xb,t)$, namely
\begin{equation}
\label{eq:potential-definition}
\vb(\xb,t) = \nablab\Phi(\xb,t) \qquad\qquad    \text{in} \ \ \Omega(t).
\end{equation}
%
In such a case, the equations of motion simplify to the unsteady
Bernoulli equation and to the Laplace equation for the flow potential:

\begin{subequations}
  \label{eq:incompressible-euler-potential}
  \begin{alignat}{3}
    \label{eq:bernoulli}
    & \frac{\Pa\Phi}{\Pa t}+\frac{1}{2}\left|\nablab\Phi\right|^2
      +\frac{p}{\rho}+\bbeta
      \cdot \xb = C(t) \qquad & \text{ in } \Omega(t)\\
    \label{eq:incompressibility-potential}
    & \Delta \Phi = 0 & \text{ in } \Omega(t)
  \end{alignat}
\end{subequations}
%
where $C(t)$ is an arbitrary function of time. In
Equation~\eqref{eq:bernoulli} we have assumed
that all body forces can be expressed as $\bb = -\nablab (\bbeta
\cdot \xb)$, i.e., they are all of potential type. This is true for
gravitational body forces and for inertial body forces due to uniform
accelerations along fixed directions of the frame of reference.

We remark that the unknowns of System~\eqref{eq:incompressible-euler-potential}
$\Phi$ and $p$ are uncoupled. Once the solution of Laplace
Problem~\eqref{eq:incompressibility-potential} is obtained, the pressure can
be in fact recovered by means of a postprocessing step through  Bernoulli's
Equation~(\ref{eq:bernoulli}).


\subsection{Perturbation potential}
\label{sec:pert_pot_prob}

\begin{figure}[htb!]
\centerline{
  \ifpdf
  \resizebox{0.7\textwidth}{!}{
    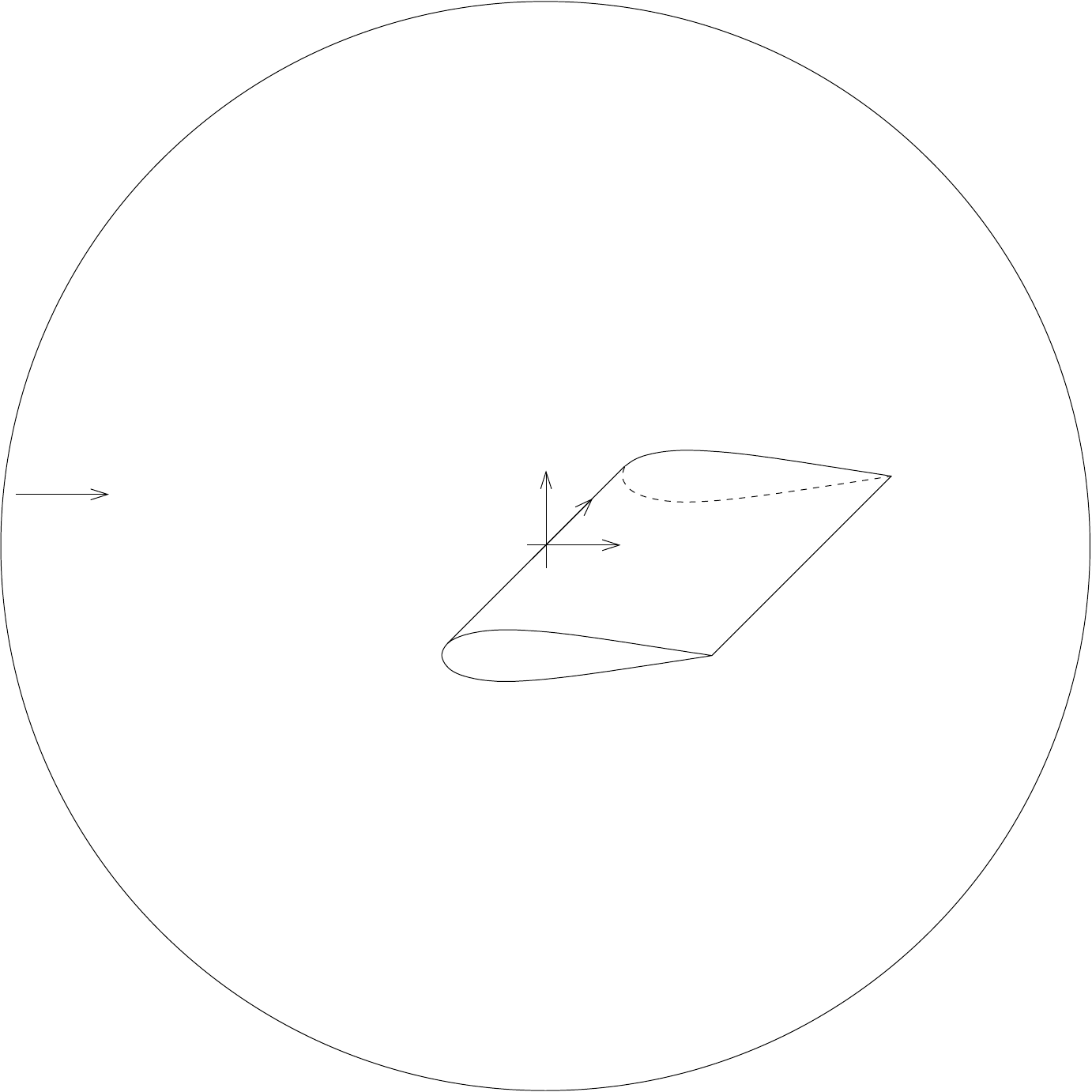
  }
  \else
  \resizebox{0.7\textwidth}{!}{
    \input{./figures/naca_domain.pstex_t}
  }
  \fi
}       	
\caption{The computational domain \label{fig:wholeDomain}}
\end{figure}

As mentioned, we solve Problem~\eqref{eq:incompressible-euler-potential} in the
region $\Omega(t)$ surrounding the wing. In this work, we will assume that the wing is moving in a fluid at rest
at velocity $\Vb_{\infty}=V_{\infty}\hat{\ib}$. 
Exploiting Galilean relativity, we then solve the governing equations in a 
relative frame of reference moving with the wing average speed. In such an inertial frame, the $x$ axis of the
right-handed coordinate system used will be aligned with the asymptotic fluid velocity, the $z$ is aligned with
the vertical --- gravity acceleration --- direction, and the $y$ axis identifies a lateral
direction normal to both $x$ and $z$. Figure~\ref{fig:wholeDomain} depicts
the domain $\Omega(t)$, along with the asymptotic velocity vector and coordinate system.
We remark that possible accelerations of the wing in the inertial frame chosen can be accounted for by
specifying, at each time instant, the updated position and velocity of the wing surface points.
Such a strategy would allow for the simulation of possible pitching or plunging motions of the wing under
the action of fluid dynamic forces, or rotations around a shaft experienced by propeller blades.
However, in all the numerical tests and simulations presented in this work the
wing is not experiencing time dependent motions in the relative reference frame described.
We also point out that, due to the absence of inertia in the potential flow model, whenever the boundary
conditions are independent of time, a steady solution is expected. For such a reason, and to lighten the
notation, we will omit the explicit time dependence of both $\Omega$ and its boundaries and the solution
in the following sections.

%
It is convenient to split the potential $\Phi$ as the sum between
a mean flow potential (due to the free stream velocity) and the so called \emph{perturbation potential}
$\phi$ due to the presence of the wing, namely


\begin{subequations}
  \begin{alignat}{1}
    \Phi(\xb) &= \Vb_{\infty}\cdot\xb +
    \phi(\xb)\\
    \vb(\xb) = \nablab\Phi(\xb) &= \Vb_{\infty} +
    \nablab\phi(\xb).
  \end{alignat}
  \label{eq:potSplit}
\end{subequations}
%
The perturbation potential still satisfies Poisson equation

\begin{equation}\label{eq:laplacian_perturb}
\Delta\phi = 0 \qquad\qquad \text{in} \ \ \Omega,
\end{equation}
and in practice it is convenient to solve for $\phi$, and obtain the
total potential $\Phi$ from equation (\ref{eq:potSplit}).

Figure~\ref{fig:wholeDomain} presents a sketch of the explicit splitting
of the various parts of the boundary $\Gamma$. On the wing surface $\Gamma_{B}$,
we apply a non penetration condition which 
takes the form

\begin{equation}\label{eq:non_penetr_bc}
\phin := \nablab\phi\cdot \nb= -\nb\cdot\Vb_{\infty}
\qquad\qquad \text{on} \ \ \Gamma_{B},
\end{equation}
when expressed in terms of the perturbation potential and
$\Vb_{\infty}$. 

As the name suggests, the perturbation potential represents the flow potential component associated
with the presence of the wing in the otherwise undisturbed stream. Clearly, 
the influence of a body having finite dimensions on the flow field is expected to fade at infinite distance
from the body itself. Thus, a possible condition for the perturbation potential on the far field boundary $\Gamma_\infty$ is the
homogeneous Dirichlet condition


\begin{equation}\label{eq:far_field_bc}
\phi = 0 \qquad\qquad \text{on} \ \ \Gamma_{\infty}.
\end{equation}

A more accurate estimate of the expected asymptotic behavior of the perturbation potential will be 
obtained based on the boundary integral formulation of the problem described in Section \ref{sec::bound_int_form}.

\subsection{Wake modeling and quasi-potential flow}
\label{sec:quasi_pot}

The model described by Equation \eqref{eq:laplacian_perturb} and Boundary Conditions \eqref{eq:non_penetr_bc} and \eqref{eq:far_field_bc}
results unfortunately unable to predict that forces are exchanged between the ideal fluid and a body moving at constant speed in it.
This fact is known as D'Alembert's paradox \cite{GRIMBERG20081878}. However, it is possible to modify the potential flow model so as
to extend its limits and make it able to compute accurate approximations of the fluid dynamic forces exchanged by the fluid and
a streamlined body. This leads to the so called \emph{quasi-potential} flow theory. We refer the interested reader
to \cite{morinoGennaretti1992,GENNARETTI1996281} for a more detailed analysis such flow model, while in this section we will only report its main features.

\begin{figure}[htb!]
\centerline{
  \ifpdf
  \resizebox{0.7\textwidth}{!}{
    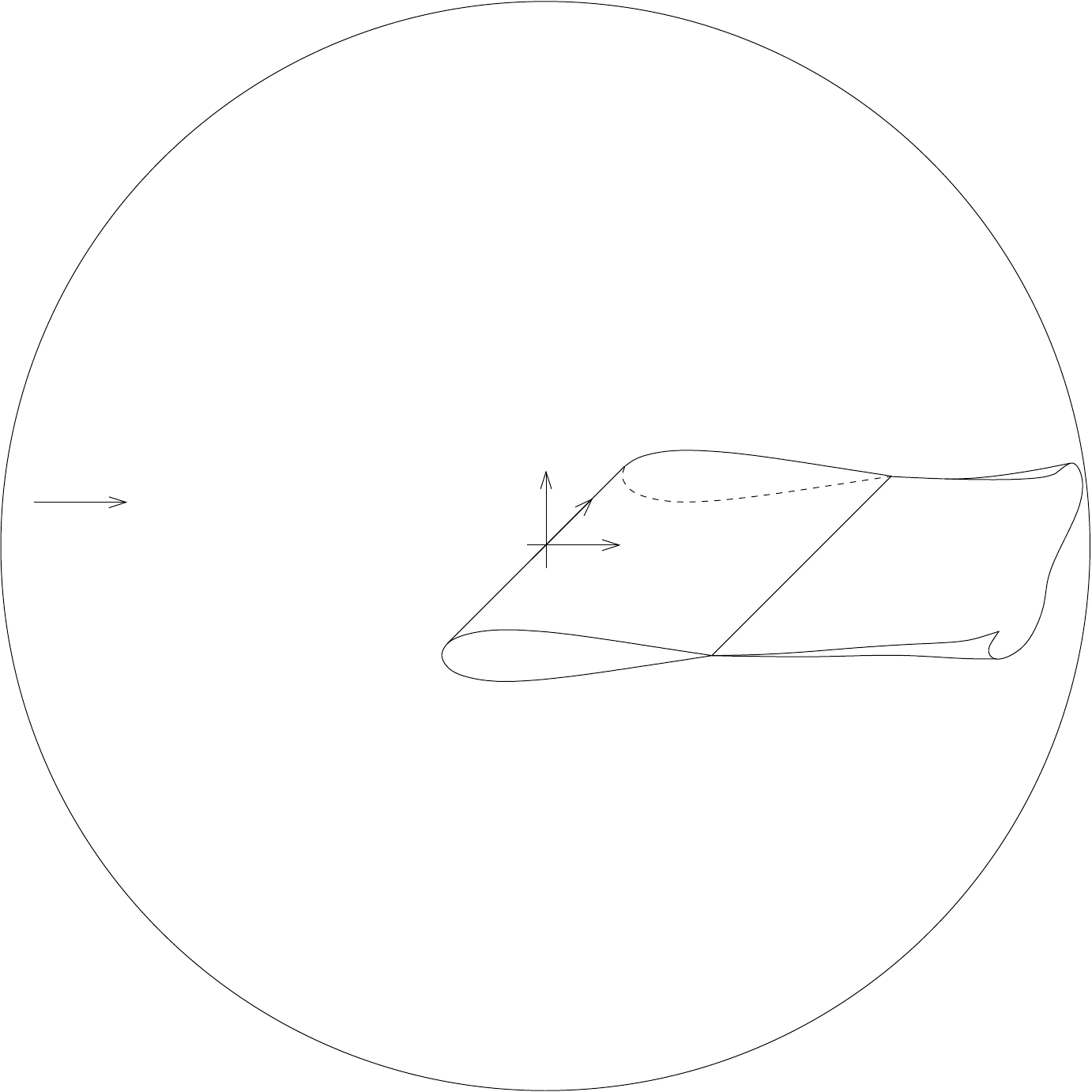
  }
  \else
  \resizebox{0.7\textwidth}{!}{
    \input{./figures/naca_domain_wake.pstex_t}
  }
  \fi
}       	
\caption{The computational domain including the presence of the wake sheet $\Gamma_W$ detaching from the trailing edge $\gamma_{TE}$\label{fig:wholeDomainWake}}
\end{figure}

In a quasi-potential model, the flow is assumed irrotational in all the flow domain except from
a free vortex sheet of null thickness detaching from the trailing edge. In this way,
the effects of the finite thickness vortical shear layer detaching from the wing or airfoil trailing edge are recovered,
in the potential flow theory framework, by means of a wake of null thickness.
The surface $\Gamma_W$ representing such a wake sheet (see Figure \ref{fig:wholeDomainWake}) is by all means an additional
boundary of the domain $\Omega$. As such, on each point $\xb\in\Gamma_W$ the limits of the potential as the surface is approached
from the top side $\phi^+(\xb)$ or from the bottom side $\phi^-(\xb)$ are in principle different.    

Enforcing the conservation of mass as in \cite{MORINO2001805} on a infinitely thin control volume sitting across $\Gamma_W$ yields

\begin{equation}\label{eq:mass_cons_wake}
\frac{\Pa \phi}{\Pa n}^+ = \frac{\Pa \phi}{\Pa n}^-  \qquad \text{on} \ \ \Gamma_W.
\end{equation} 

On the same infinitely thin control volume across $\Gamma_W$, under the assumption that pressure is continuous across the wake sheet,
the enforcement of the momentum conservation results in

\begin{equation}\label{eq:momentum_cons_wake}
\frac{D (\delta\phi)}{D t} = \frac{\Pa (\delta\phi)}{\Pa t} +\Vb_w\cdot\nablab(\delta\phi) = 0 \qquad \text{on}\ \  \Gamma_W,
\end{equation} 
in which

\begin{equation}\label{eq:delta_phi}
\delta\phi = \phi^+-\phi^-
\end{equation} 
is the potential jump across the wake surface and
$\Vb_w = \frac{(\nablab\phi^++\nablab\phi^-)}{2}+\Vb_\infty$ is the mean fluid velocity on $\Gamma_w$. 
Equation \eqref{eq:momentum_cons_wake} indicates that the potential jump remains constant following a
material point on the wake which moves with velocity $\Vb_w$. As a consequence, once the values
of $\delta\phi$ on the trailing edge points are known, they can be used to readily obtain the potential jump
field on the entire wake. Of course, this requires that on every point of the wake boundary, the latter is
aligned with the local wake velocity. Unfortunately, since the flow velocity is an unknown of our mathematical problem,
the shape of the wake boundary $\Gamma_w$ cannot be known \emph{a priori}. So, the values of $\delta\phi$ on the trailing
edge points and the effective shape of the wake are additional unknowns of the potential flow problem 
of the quasi-potential model governing equations. Section \ref{sec:kutta} will discuss
the supplementary condition added to the problem to allow for the computation of the additional unknown  $\delta\phi$.
Section \ref{sec:wake_relaxation} will provide details of the fixed point wake relaxation algorithm used to compute the
correct shape and position of the wake during the simulations.

\subsection{Kutta condition}\label{sec:kutta}

As mentioned, the wake vortex sheet is introduced to account
for the important effects of the thin shear layer detaching from the trailing edge of an airfoil or wing. In the idealization
of the quasi potential model, all the vorticity contained in such a shear layer is in fact carried within a null thickness
wake sheet. As a consequence, the velocity potential becomes discontinuous across the wake. 

A discontinuity in the potential is an acceptable feature in the simulations, but, of course, it should not translate
into a discontinuity in the pressure field. 
Equation \eqref{eq:momentum_cons_wake} is derived based on the assumption that the pressure is continuous across the wake. Thus,
aligning the wake with the local velocity field ensures that pressure is continuous on the two sides of the wake. We then only
need to impose pressure continuity on the points of the trailing edge, namely 

\begin{equation}\label{eq:kutta}
\delta p = p^+-p^- = 0 \qquad \text{on}\ \  \gamma_{TE},
\end{equation} 
where $p^+$ and $p^-$ are the pressure values computed on the leeward and windward side of the trailing edge, respectively.
The latter condition is the one added in this work to obtain the correct value of the additional trailing edge unknown $\delta \phi$. 
There are of course many forms available in the literature of the so called \emph{Kutta condition} needed to find the correct $\delta \phi$
value (see among others \cite{XU1998415}). 

As pointed out in Section \ref{sec:potential_flow}, the pressure value is obtained making use of Bernoulli's Equation \eqref{eq:bernoulli}.
Given the steady flow assumption used for all the simulations considered in this work, a perturbation potential based expression for pressure
reads

\begin{equation}\label{eq:pressure}
p(\xb,t) = \frac{1}{2}\rho\left(|\Vb_\infty|^2-|\nablab\phi+\Vb_\infty|^2\right) -\rho g z\qquad \text{in} \ \ \Omega,
\end{equation} 
where in the derivation of Equation~\eqref{eq:pressure} we have assumed that the only inertial force considered is due to a
gravity acceleration of module $g$ aligned with the $z$ axis. Kutta condition will then read

\begin{equation}\label{eq:kutta_vel}
|\nablab\phi^++\Vb_\infty|^2 = |\nablab\phi^-+\Vb_\infty|^2 \qquad \text{on}\ \  \gamma_{TE}.
\end{equation}

\section{Numerical approximation based on the Boundary Element Method}\label{sec:numerical_model}

In this work, we make use of the Boundary Element Method (BEM) for the
spatial discretization of the governing Laplace equation. In the context
of potential flow simulation, this is quite a common choice. We must
however remark that the Laplace equation for the velocity potential can be also discretized by means of
the Finite Element Method (FEM). At a first glance, it would appear that the most important advantage of BEM compared to FEM 
is the reduced number of unknowns associated with the codimension one grid. Unfortunately, in the practice such
an advantage is offset by the presence of a dense resolution matrix in the discretized algebraic system.
Yet, there are other advantages of BEM that made us favor it over FEM. In particular, in the context of the present
physical problem, the codimension one grids required by BEM are much 
easier to generate, deform and manage without significant quality drop.

For the BEM discretization we use the same formalism presented in~\cite{Giuliani2015,brebbia,pi-BEM}

\subsection{Boundary integral formulation}
\label{sec::bound_int_form}

To rewrite \eqref{eq:laplacian_perturb} as a Boundary Integral Equation (BIE) we make use of a \emph{fundamental solution} (or Green's function) of Laplace equation. More specifically, in this work we employ the free space Green's function

\begin{equation*}
G(\yb-\xb) = \frac{1}{4\pi|\yb-\xb|},
\end{equation*}
which is the distributional solution of the following problem

\begin{eqnarray*}
-\Delta G(\yb-\xb)&=&\delta(\xb)  \quad \text{in}\ \mathbb{R}^3\\
\lim_{|\yb|\rightarrow \infty} G(\yb-\xb)& = &0,
\end{eqnarray*}
where $\yb\in\mathbb{R}^3$ is a generic point, and
 $\xb\in\mathbb{R}^3$ is the center of the Dirac delta distribution $\delta(\xb)$.
If we select $\xb$ to be inside $\Omega$, using the deﬁning property of
the Dirac delta and the second Green identity, we obtain

\begin{equation}
\label{integral}
  \phi(\boldsymbol{x}) =\int_{\Gamma}G(\yb-\xb)\frac{\partial \phi}{\partial n}(\boldsymbol{ y})\d s_y
-\int_{\Gamma}\phi(\boldsymbol{ y}){\nablab G(\yb-\xb)\cdot\nb(\yb)}\d s_y,
\ \ \forall\boldsymbol{x}\in\Omega,
\end{equation}
where in this case $\yb\in\Gamma$ is a generic integration point on the domain boundary --- as indicated by the subscript
in differential $\d s_y$ --- and $\nb(\yb)$ is the outward unit normal vector to boundary $\Gamma$. The Green's function
gradient is $\nablab G=\frac{\yb-\yb}{-4\pi|\yb-\xb|^3}$.

From~\eqref{integral} we notice that if the solution and its normal derivative on the boundary $\Gamma$ are known then the potential $\phi$ can be computed in any point of the domain. 
Considering the trace of~\eqref{integral} we can write the boundary integral form of the original problem as 
\begin{eqnarray} 
  c(\boldsymbol{x})\phi(\boldsymbol{x}) &=&  \int_{\Gamma}G(\yb-\xb)\frac{\partial \phi}{\partial n}(\boldsymbol{y})\d s_y\\
&&-\int^{PV}_{\Gamma}\!\!\!\!\phi(\boldsymbol{ y}){\nablab G(\yb-\xb)\cdot\nb(\yb)}\d s_y \quad 
\forall\boldsymbol{x}\in\Gamma=\Pa\Omega\label{BoundIntFormul}
\end{eqnarray}
where $c(\boldsymbol{x})$ is obtained from the Cauchy Principal Value (CPV) evaluation of the integral involving the derivative of
the { Green's function} --- as indicated by the $\int^{PV}$ symbol. By a geometric perspective, the value of $c$ represents the fraction of
solid angle of the domain $\Omega$ seen from the boundary point $\boldsymbol{x}$. As we are interested in computing the solution on the
surface of the lifting body, the point $\boldsymbol{x}$ is here assumed to always lie on $\Gamma_B$.
We then need to introduce the boundary conditions of the Laplace problem described in Sections \ref{sec:pert_pot_prob} and \ref{sec:quasi_pot}
to obtain the values of the integrals appearing in Equation \eqref{BoundIntFormul} and formulate the discretized lifting body problem.
To this end, we split the domain boundary into its non overlapping regions illustrated in Figure \ref{fig:wholeDomainWake},
namely $\Gamma=\Gamma_B\cup\Gamma_\infty\cup\Gamma_{W+}\cup\Gamma_{W-}$ and analyze separately the contributions of each region to the integrals
in Equation \eqref{BoundIntFormul}.

\paragraph{Wake surface integrals}

Note that $\Gamma_{W+}$ and $\Gamma_{W-}$, are the top and bottom boundaries, respectively, of the potential flow region interfacing it with the
vortical wake region. As discussed in Section \ref{sec:quasi_pot}, in the framework of the quasi-potential model the vortical region past the trailing edge
is assumed to have zero thickness, and it is represented as a vortex sheet surface. Thus $\Gamma_{W+}$ and $\Gamma_{W-}$ will be characterized by the
same position and shape --- which coincides with  $\Gamma_{W}$ --- but by opposite orientation and normal vectors. Given these considerations, the contributions
of the wake boundary $\Gamma_{W}$ to the integrals appearing in Equation \eqref{BoundIntFormul} result in

\begin{eqnarray}
\int_{\Gamma_{W+}}G(\yb-\xb)\frac{\partial \phi}{\partial n}^+(\boldsymbol{y})\d s_y\!\!\!\! &+& \!\!\!\!\int_{\Gamma_{W-}}G(\yb-\xb)\frac{\partial \phi}{\partial n}^-(\boldsymbol{y})\d s_y  = 0\\\nonumber
\int^{PV}_{\Gamma_{W+}}\!\!\!\!\phi^+(\boldsymbol{ y}){\nablab G(\yb-\xb)\cdot\nb(\yb)}\d s_y \!\!\!\!&+& \!\!\!\!\int^{PV}_{\Gamma_{W-}}\!\!\!\!\phi^-(\boldsymbol{ y}){\nablab G(\yb-\xb)\cdot\nb(\yb)}\d s_y \\
&=& \!\!\!\!\int^{PV}_{\Gamma_W}\!\!\!\!\delta\phi(\boldsymbol{ y}){\nablab G(\yb-\xb)\cdot\nb(\yb)}\d s_y,
\end{eqnarray} 
where we have used Equations \eqref{eq:mass_cons_wake} and \eqref{eq:delta_phi}, and the opposite sign of the normal unit vectors on the two sides of the wake surface.

\paragraph{Far field surface integrals}

As for boundary $\Gamma_\infty$, we can assume, without loss of generality, that it can be represented as the surface of a sphere with
radius $R\rightarrow\infty$ (see Figure \ref{fig:wholeDomainWake}). Making use of spherical coordinates ($r,\varphi,\theta$) a generic point $\xb\in\Gamma_{\infty}$
can be expressed as $\xb=R\sin\varphi\cos\theta\hat{\ib}+R\sin\varphi\sin\theta\hat{\jb}+R\cos\varphi\hat{\kb}$ and the BIE integrals on $\Gamma_{\infty}$ read
\begin{eqnarray}
\label{eq:far_field}
&\displaystyle\int_{\Gamma_{\infty}}&\!\!\!\!\!\!\!\!\!G(\yb-\xb)\frac{\partial \phi}{\partial n}(\boldsymbol{y})\d s_y = \lim_{R\rightarrow\infty}\int_0^\pi\!\! \sin\varphi\d \varphi \int_0^{2\pi} \!\!\!\frac{1}{4\pi R} \frac{\partial \phi}{\partial n}(R,\varphi,\theta)R^2\d \theta \\\nonumber
&\displaystyle\int_{\Gamma_{\infty}}&\!\!\!\!\!\!\!\!\!\!\!\phi(\boldsymbol{ y}){\nablab G(\yb-\xb)\cdot\nb(\yb)}\d s_y  = \lim_{R\rightarrow\infty}\int_0^\pi\!\! \sin\varphi\d \varphi \int_0^{2\pi} \!\!\frac{R}{4\pi R^3}\phi(R,\varphi,\theta)R^2\d \theta,
\end{eqnarray}
where we have used the fact that as $R\rightarrow\infty$ $\xb-\yb\simeq \xb$ and $(\xb-\yb)\cdot \nb_{\Gamma_\infty}\simeq R$. 
Equations \eqref{eq:far_field} suggest that assuming $\phi=o(1)$ and $\frac{\Pa \phi}{\Pa n}=o(\frac{1}{R})$ as $R\rightarrow\infty$ results in null values for both
integrals on $\Gamma_\infty$. Such assumptions are by no means restrictive. Requiring that the perturbation potential is a function that fades at an infinite
distance from the body generating the perturbation is in fact consistent with its very definition. In light of this, the perturbation potential normal derivative must present an
asymptotic behavior proportional to the derivative of the perturbation potential, which makes the assumption $\frac{\Pa \phi}{\Pa n}=o(\frac{1}{R})$ also reasonable. For these
reasons, in this work the contribution of $\Gamma_\infty$ to the boundary integrals in Equation \eqref{BoundIntFormul} is considered null. The discretization of such boundary
will then not be necessary.

\paragraph{Body surface integrals}

On the body boundary $\Gamma_B$, making use of  Equation \eqref{eq:non_penetr_bc}, we can write

\begin{equation}
\int_{\Gamma_{B}}G(\yb-\xb)\frac{\partial \phi}{\partial n}(\boldsymbol{y})\d s_y = -\int_{\Gamma_{B}}G(\yb-\xb)(\nb(\yb)\cdot\Vb_\infty)\d s_y  
\end{equation}

The final boundary integral formulation of the problem, written for a generic point $\xb\in\Gamma_B$, reads

\begin{eqnarray} \label{BoundIntFormul-2}
  c(\boldsymbol{x})\phi(\boldsymbol{x})\!\!\!\! &+&\!\!\!\! \int^{PV}_{\Gamma}\!\!\!\!\phi(\boldsymbol{ y}){\nablab G(\yb-\xb)\cdot\nb(\yb)}\d s_y + \\ \nonumber
\!\!\!\!&&\!\!\!\!  \int^{PV}_{\Gamma_W}\!\!\!\!\delta\phi(\boldsymbol{ y}){\nablab G(\yb-\xb)\cdot\nb(\yb)}\d s_y = -\int_{\Gamma_{B}}\!\!\!\!G(\yb-\xb)(\nb(\yb)\cdot\Vb_\infty)\d s_y 
\end{eqnarray}

\subsection{Discretisation}
\label{sec:discretisation}
The numerical discretization of~\eqref{BoundIntFormul-2} leads to a real-valued Boundary Element Method (BEM). The resolution of a BEM requires the
discretization of the unknowns using functional spaces defined on a Lipschitz boundary. We address this problem introducing suitable discretizations
of the unknown $\phi$ on the body surface and of the unknown $\delta\phi$ on the wake surface.  Such discretizations are based on standard Lagrangian
finite element spaces defined on $\Gamma$. We then use the same functional space to describe the geometry. This setting is often referred to as
Isoparametric BEM.

We define the computational mesh as a quadrilateral decomposition ${\Gamma}_h$ of the boundary $\Gamma$. We require that two cells $K,\ K'$ of the mesh only
intersect on common edges or vertices, and that there exists a mapping from a reference cell $\hat K$ to $K$. The Jacobian must be uniformly bounded away from
zero for all cells $K$. 
{
To ease mesh generation, the simulation tool developed allows for the definition of a very coarse initial grid, which is then automatically refined on
the user prescribed geometry of the wing up to the desired level of refinement. Following~\cite{pi-BEM-repo, dealII92} an interface to CAD data
structures --- which are the most common tool to define arbitrary geometrical
descriptions~\cite{heltaiEtAlACMTOMS2021,BangerthHeisterHeltai-2016-b,ArndtBangerthDavydov-2017-a} --- is used
to specify the desired geometry of the wing. This feature has been previously employed in hydrodynamics simulations through
BEM, \cite{waveBem,Mola2014,MolaHeltaiDeSimone-2017-a}, and~\cite{pi-BEM} presents an example of a Laplacian problem with known analytic solution
solved for convergence testing purposes on the geometry of an aeronautics-like NACA wing shape.
As will be shown in Section~\ref{sec:convergence}, such CAD-interfacing feature has also been used in this work for the proper placing of the
degrees of freedom on the CAD surface whenever a high order finite element is selected for the mapping and the solution.  We point out that the
latter feature is ultimately what results in the possibility for the user to effectively select arbitrary order finite elements at the start
of each simulation. In fact, the generation of quality unstructured high-order meshes remains a significant obstacle in the adoption of
high-order finite elements. In the framework of the methodology here presented however, only a standard low order grid is required at
the start of the simulation, as possible additional high order mapping degrees of freedom are automatically generated on the CAD surfaces.   
}
 
If $\phi$ and $\frac{\Pa\phi}{\Pa n}$ lie in the spaces $V$, defined as

\begin{equation*}
\begin{aligned}
V &:= \left\{ \phi \in H^{\frac{1}{2}}(\Gamma)   \right\} ,
\end{aligned}
\end{equation*}
where $\Gamma=\partial\Omega$, then the integrals in Equation \eqref{BoundIntFormul-2} are bounded.
$H^{\frac{1}{2}}(\Gamma)$ is the space of traces on $\Gamma$ of functions in $H^{1}(\Omega)$. 
We construct the discretized space $V_h$ as conforming finite dimensional subspace of $V$

    \begin{align}\label{eq:shape_functs}
      V_h &:= \left\{ \phi_h \in C^0(\Gamma_h)  :  \phi_{h|K} \in \mathcal{Q}^{r}(K), \,K \in \Gamma_h   \right\}  \equiv \text{span}\{\psi_i\}_{i=1}^{N_V}  
    \end{align}
where  $\mathcal{Q}^{r}(K)$ is the space of polynomials of order $r$ in each coordinate direction.  In Equation \eqref{eq:shape_functs}
$\psi_i$ denotes the shape function associated with the $i$-th degree of freedom of the discretized space. For additional details on the
finite elements used in the BEM formulation we refer the interested reader to \cite{GiulianiPoliMIThesis}.
Following~\cite{pi-BEM} we use iso-parametric discretizations based on standard $Q_N$ Lagrangian finite elements, and
by collocating the support points of the geometry patches directly on the CAD surfaces. This work reports results obtained with bi-linear
elements ($r=1$), as well as of higher order bi-quadratic ($r=2$) or bi-cubic ($r=3$) elements. 

\begin{figure}[htb!]
\centerline{
  \ifpdf
  \resizebox{0.9\textwidth}{!}{
    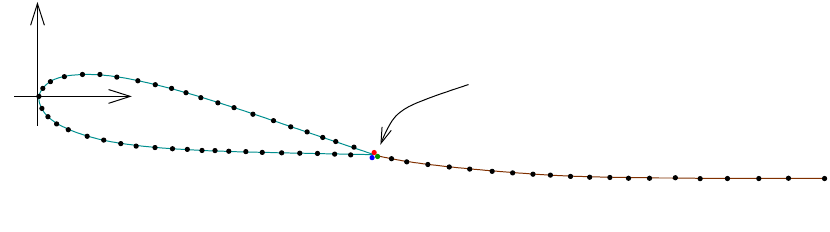
  }
  \else
  \resizebox{0.9\textwidth}{!}{
    \input{./figures/airfoil_algorithm.pstex_t}
  }
  \fi
}       	
\caption{A two dimensional sketch of the computational domain including details of the collocation nodes, coinciding with on the degrees of freedom (DOFs) of the
         discretization space. The image also indicates that in correspondence with the trailing edge $\gamma_{TE}$, a triple DOF is present to allow for the
         space $V_h$ functions to have different trailing edge values on the cells of the airfoil leeward side, windward side, and wake. \label{fig:airfoil_algo}}
\end{figure}

Since the unknowns $\phi$ and $\delta\phi$ are defined on complementary portions of the discretized domain, we can define a single function in the discretized
space $\phi_h \in V_h$ that approximates $\phi$ on the body surface and $\delta\phi$ on the wake surface. The generic element of the discretized space reads

\begin{equation}
\phi_h(\xb)=\sum_{j=1}^{N_V}\hat{\phi}_j \psi_j(\xb), 
\end{equation}
where the vector $\hat{\phib}\in\mathbb{R}^{N_V}$ contains the value of the discretized potential $\phi_h$ on each collocation point
located on the wing surface $\Gamma_B$, and of the potential jump $\delta\phi_h$ on each collocation point located on the wake.  $N_V$ then represents
the overall number of degrees of freedom of the discretized space $V_h$.
To allow for possible discontinuities of the unknown functions across the trailing edge $\gamma_{TE}$ we use a methodology inspired to the double nodes technique
presented in \cite{grilli2001}. More precisely, the use of a continuous Lagrangian Finite Element space is in principle not suited for the
representation of the solution at the trailing edge. On one hand in fact, the value of the potential obtained approaching the trailing edge from the windward
side of the airfoil is different from that obtained approaching the same point from the leeward side. On the other hand, the value obtained approaching the
same point from the wake side must be the difference between the two airfoil side values, as stated by Equation \eqref{eq:delta_phi}. To accommodate for
this, the computational mesh will be built in such a way that  in correspondence with the sharp trailing edge, the degrees of freedom of the discretization
space will be triple (see Figure \ref{fig:airfoil_algo}). This allows for the functions of
space $V_h$ to have different trailing edge values on the cells of the airfoil leeward side, windward side, and wake. Note that the increased number of unknowns
implies that additional equations must be introduced for the closure of the mathematical problem. As it will be shown in Section \ref{sec:num_kutta}, such equations
will be obtained making use of Equation \eqref{eq:delta_phi} and by imposing a discretized version of the Kutta condition \eqref{eq:kutta_vel}.   

The collocation method is a common resolution technique for a BEM since it does not require any additional integration of~\eqref{BoundIntFormul-2}.
We refer the interested reader to~\cite{pi-BEM} for a detailed analysis of the accuracy of this setting.
Collocating Equation \eqref{BoundIntFormul-2} results in the following linear system 
\begin{equation}
(C+N) \hat{\phib} - \bb = 0,
\label{LaplaceLinearSystem}
\end{equation}
where, for all the row indices $i$ corresponding to collocation points $\xb_i\in\Gamma_B$:
\begin{itemize}
\item{$C$ is a diagonal matrix with diagonal entries
    $C_{ii} = c(\xb_i)$}
\item{$N_{ij}=\sum_{k=1}^K \sum_q^{N_q} \dfrac{\partial G}{\partial n}(\xb_i
    - \xb_q) \psi_q^j J^k(\xb_q) $, where $J^k$ is the determinant of the first fundamental
    form for each panel $K$, and $\xb_q$ are the $N_q$ quadrature nodes
    used to compute the integral at the numerical level;}
\item{$b_{i}=\sum_{k=1}^K \sum_q^{N_q} G(\xb_i - \xb_q) (-\Vb_\infty\cdot\nb(\xb_q)) J^k(\xb_q) $;}
\end{itemize}

When the collocation point lies inside the current integration cell, we make use of bidimensional quadrature
formulas based on Duffy transformation as in \cite{Mousavi2010} to treat singular kernel integrals; in any other case,
we use standard Gauss quadrature rules.

It is important to point out that the linear system equation written on each row corresponding to a wake collocation point is
not the discretized BIE, but is instead the discretized version of Equation \eqref{eq:momentum_cons_wake}. This requires
the use of a structured grid in such a region, and that the grid is conformal along the trailing
edge $\gamma_{TE}$. An example of structured wake grid is depicted in Figure \ref{fig:meshRef}. With this kind of grids it
is possible to arrange the degrees of freedom of the wake in lines that sit on the pathlines detaching from the trailing edge.
The wake relaxation algorithm used to align such wake degrees of freedom lines with
the pathlines is described in Section \ref{sec:wake_relaxation}. Thus, for all the row indices $i$ corresponding to collocation points $\xb_i\in\Gamma_W$:
\begin{itemize}
\item computing integrals for the BEM matrices is not needed on the wake, as on $\Gamma_W$ Equation \eqref{eq:momentum_cons_wake}
      is exploited;
\item the system matrix {$(C+N)$} is used to impose that the value of $\delta\phi$ at each degree of freedom on the
      wake is equal to that of the previous degree of freedom on the same pathline. Thus, for each line the only non null
      entries will be a 1 on the principal diagonal, and a -1 in correspondence with the index of the previous degree of freedom
      along the pathline;
\item the right hand side entries $b_i$ are null in this region.
\end{itemize}

\subsection{Numerical enforcement of the nonlinear Kutta condition}
\label{sec:num_kutta}

System \eqref{LaplaceLinearSystem} has to be completed with the additional conditions on the trailing edge degrees of freedom
to close the problem and obtain the desired solution. As discussed in the previous section and illustrated in
Figure \ref{fig:airfoil_algo}, each collocation point on $\gamma_{TE}$ will correspond to three degrees of freedom in the
resolution system. The equations written on the row corresponding to each of these degrees of freedom in the system will be:

\begin{itemize}
\item \textbf{Wing leeward side:} on the row corresponding to this degree of freedom, we write the
      Boundary Integral Equation \eqref{BoundIntFormul-2}, centered at the collocation point on the trailing edge.  
\item \textbf{Wing windward side:} on the row corresponding to this degree of freedom, we write the discretized  Kutta condition \eqref{eq:kutta_vel}, namely
      \begin{equation}\label{eq:kutta_vel_2}
      |\nablab\phi^+_{lw}+\Vb_\infty|^2 = |\nablab\phi^-_{ww}+\Vb_\infty|^2,
      \end{equation}
      where the subscripts $lw$ and $ww$ indicate the local degree of freedom index on the airfoil leeward side and airfoil windward side,
      respectively.
\item \textbf{Wake degree of freedom:} in this degree of freedom, the unknown quantity is the potential jump $\delta\phi$ on the
      wake point attached to the trailing edge. As such, the equation written is
      \begin{equation}\label{eq:kutta_pot}
      (\delta\phi)_{wk} = \phi_{lw}-\phi_{ww},
      \end{equation}
      where the subscript $wk$ indicates the local degree of freedom index on the wake. 
\end{itemize} 

Introducing Equation \eqref{eq:kutta_pot} in the resolution system does not pose significant problems, as it translates
in adding a line in correspondence with the ${wk}$-th degree of freedom, in which only a unit diagonal term and two negative unit off diagonal terms
at columns ${lw}$ and ${ww}$ are present. On the other hand, the discretization of Equation \eqref{eq:kutta_vel_2} requires special
attention, as the gradients of the potential function are not a single valued function if evaluated in correspondence with
collocation points located on cells edges. We then resort to a $L_2$ projection strategy to obtain suitable approximations $\chi_x,\chi_y,\chi_z\in V_h$ of the
components of the surface gradient $\nabla_s\phi_h$ combined with the potential normal
derivative $\left(\frac{\partial \phi}{\partial n}\right)_h=-\Vb_\infty\cdot\nb$ and
the asymptotic velocity $\Vb_\infty$. The weak form of the velocity $x$ component computation reads
\begin{equation}
\int_{\Gamma}  v_h \chi_x \d s = \int_{\Gamma} v_h \left(\nabla_s \phi_h  + \left(-\Vb_\infty\cdot\nb\right) \nb
+\Vb_\infty\right)\cdot\hat{\ib}\d s  \quad\forall v_h \in V_h.
\label{gradientrecovery}
\end{equation}
From the discretization of Equation \eqref{gradientrecovery} we get the $N$ dimensional sparse system

\begin{equation}
M\hat{\chib}_x = \bb^{\chi_x},
\end{equation}
where

\begin{itemize}
\item $M:\ m_{ij} = \int_{\Gamma}\psi_i\psi_j\d s$ is a sparse mass matrix;
\item $\bb^{\chi_x}:\ b^{\chi_x}_i=\int_{\Gamma}\psi_i\left(\nabla_s \phi_h  -\Vb_\infty\cdot\nb
+\Vb_\infty\right)\cdot\hat{\ib}\d s$  is the right hand side vector;
\item $\hat{\chib}_x$ is the vector containing the values of the fluid velocity on each collocation point.
\end{itemize}
Note that the value of the local discretized velocity at each integration point is computed as

\begin{equation}\label{eq:velocity_on_quad_point}
\nabla_s \phi_h - \Vb_\infty\cdot \nb+\Vb_\infty = \sum_{j=1}^{N_{V}}\left((\nabla_s \phi)_j  - \Vb_\infty\cdot \nb+\Vb_\infty\right)\psi_j.
\end{equation}

Since all the standard Gauss quadrature points are located inside the grid cells, the velocity evaluated with Equation \eqref{eq:velocity_on_quad_point}
is single valued. Thus, $\chi_x(\xb)=\sum_{j=1}^{N_V}\hat{\chi_x}_j \psi_j(\xb)$ is the best possible approximation, in the $L_2$ norm sense,
of the $x$ velocity component in space $V_h$. 

A system similar to \eqref{gradientrecovery} --- with identical matrix $M$ --- is
solved for the $y$ and $z$ components of the velocity. Thus, the discretized version of Equation \eqref{eq:kutta_vel_2}
reads

\begin{equation}\label{eq:kutta_vel_3}
{\hat{\chi_x}_{lw}^2}+{\hat{\chi_y}_{lw}^2}+{\hat{\chi_y}_{lw}^2} = {\hat{\chi_x}_{ww}^2}+{\hat{\chi_y}_{ww}^2}+{\hat{\chi_y}_{ww}^2},
\end{equation}

\subsection{Newton iterations solution}

The final resolution nonlinear system is assembled combining Equations \eqref{LaplaceLinearSystem}, \eqref{eq:kutta_pot} and \eqref{eq:kutta_vel_3} as discussed
in Section \ref{sec:num_kutta}, where the nonlinearity is given by the quadratic term in Equation \eqref{eq:kutta_vel_3}. The final
resolution system is solved by means of the Newton method. The nonlinear solver implemented is typically able to converge to the solution in few iterations.

\subsection{Velocity computation on wake and wake relaxation}\label{sec:wake_relaxation}

After the resolution of the nonlinear system, the values of the velocity potential $\phi_h$ and of the wake potential jump $\delta\phi_h$ are available
at the body and wake collocation points, respectively. As discussed in Section \ref{sec:quasi_pot}, the position of the wake surface is an unknown of the
numerical problem. For such a reason, the values of $\phi_h$ and $\delta\phi_h$ computed through the resolution of the nonlinear system are used to obtain
the fluid velocity values at the wake collocation points, and consequently deform the wake so as to align it with the local velocity field. Once the
new wake position is computed, the nonlinear system is solved to obtain the new body surface potential and wake potential jump, leading to new wake
collocation point velocities. Such a fixed point strategy is repeated until convergence, which typically requires 5 to 10 iterations.

An important aspect of the procedure just outlined is represented by the computation of the fluid velocity in correspondence with the wake
collocation points. Once again, we point out that with our finite elements choice, the collocation points are often located on the
edges of the quadrilaterals composing the grid, where the potential gradients, and consequently the velocity, is discontinuous.  
We could in principle think of avoiding the problem, also in this case, making use of a gradient recovery strategy as in
Equation \eqref{gradientrecovery}. However, the unknown in the wake boundary is not $\phi$, but its jump $\delta\phi$ across the vortex sheet.
For such a reason, the result of the application of a gradient recovery strategy on such region would result in a function with
dubious physical interpretation. In this work, we then resort to the Hypersingular Boundary Integral Equation for the computation of the
fluid velocity on the wake. The Hypiersingular BIE is obtained taking the gradient of Equation \eqref{BoundIntFormul} with respect to $\xb$.
Such equation reads

\begin{eqnarray}
\nonumber
\ab(\xb)\phi(\xb)+\bar{\bar{C}}(\xb)\nablab\phi(\xb)\!\!\!\! & + & \!\!\!\!\int^{FP}_{\Gamma}\!\!\!\!\phi(\yb)(\nablab(\nablab G(\yb-\xb)))\cdot\nb(\yb)\d s_y\\
\label{eq:hypersingBIE}
& =& \int^{FP}_{\Gamma}\!\!\!\!\nablab G(\yb-\xb)(-\Vb_\infty\cdot\nb(\boldsymbol{y}))\d s_y
\end{eqnarray}
where the notation $\int^{FP}$ indicates the Hadamard's finite part evaluation of the integrals involving the Green's function gradients and
Hessian, $\ab$ is a free coefficient vector and $\bar{\bar{C}}$ is a free coefficient double tensor. Given the free space Green's function here used,
the double gradient term in the first integral takes the form
\begin{equation}
(\nablab(\nablab G(\yb-\xb)))\cdot\nb(\yb)=-\frac{1}{4\pi|\yb-\xb|^3}\left(\frac{(\yb-\xb)\cdot\nb}{|\yb-\xb|^2}(\yb-\xb)-\nb\right).
\end{equation}

Please note that, as pointed out in \cite{guiggianiEtAl1992}, such synthetic notation is not providing specific information on
the --- spherical --- shape of exclusion region in the finite parts evaluation. Such shape has a direct influence on the specific values
assumed by the free coefficients, as well as on the numerical quadrature rules developed Guiggiani et al. in \cite{guiggianiEtAl1992} for
the integrals evaluation, which we employ in the present work. We then refer the interested reader to \cite{guiggianiEtAl1992} for a
more detailed description of the equation and its derivation, which goes beyond the scope of the present paper.

After the solution of the potential problem, the collocation of Equation \eqref{eq:hypersingBIE} at a point $\xb\in\Gamma_W$ allows us for the
direct evaluation of the potential full gradient $\nablab\phi$, that can be in turn used to compute the local fluid
velocity $\Vb_s=\Vb_\infty+\nablab\phi_h$. Ideally, $\nablab\phi_h$ should be evaluated at each BEM collocation point on the wake surface $\Gamma_W$,
so that it can be readily used to align the wake geometry with the fluid velocity. However, in \cite{guiggianiEtAl1992} Equation \eqref{eq:hypersingBIE}
is derived under the assumptions that function $\frac{\Pa\phi}{\Pa n}$ is H\"older continuous at $\xb$, and that $\phi$ must be differentiable at the same point,
with H\"older continuous first derivatives. With our Lagrangian finite elements choice, the latter condition is clearly violated at all the collocation
points located on the cell edges and vertices. For such a reason, in this work we resort the weak form of the hypersingular BIE. The corresponding
Galerkin numerical problem will be that of finding the approximation of $\nablab\phi$, $\ub_h\in W_h$, such that 

\begin{eqnarray}
\label{eq:WeakHypersingBIE}
&&\forall \vb \in W_h  \\
\nonumber
\int_{\Gamma}\!\!\!\!\!\!\!&\!\!\!\!\!\!\!&\!\!\!\!\!\!\!\Big(\ab(\xb)\phi_h(\xb)+\bar{\bar{C}}(\xb)\ub_h(\xb)\Big)\cdot\vb(\xb)\d s_x  =  \\
\nonumber
& =& \!\!\!\!\!\int_{\Gamma}\!\!\!\underbrace{\left(\int^{FP}_{\Gamma}\!\!\!\!\!\!\!\!\Big(\nablab G(\yb-\xb)(-\Vb_\infty\!\!\cdot\!\nb(\boldsymbol{y}))
\!-\!\phi_h(\yb)(\nablab(\nablab G(\yb-\xb)))\!\cdot\!\nb(\yb)\Big)\!\d s_y\!\!\right)}_{\gb(\xb)}\!\!\cdot \vb(\xb)\!\d s_x.
\end{eqnarray}

Here, the scalar solution space $V_h$ is used to define the space $W_h$ for the vector valued solution function $\ub_h$ as

\begin{align}\label{eq:shape_functs_vect}
  W_h &:= \left\{ (u_h)_k\in V_h, k=1,2,3   \right\}  \equiv \text{span}\{\xi_i\}_{i=1}^{3N_V}  .
\end{align}
%


Note that in Equation \eqref{eq:WeakHypersingBIE} the function $\gb(\xb)$ indicated in the right hand side represents the perturbation potential
gradient computed on point $\xb$ of the external integration loop. In this regard, the main advantage of the Galerkin formulation used for the
gradient computation problem in Equation \eqref{eq:WeakHypersingBIE} is that, in the framework of the double integrals computation, the
Hypersingular BIE singularities are always located on the --- Gauss --- quadrature nodes. Since such points are always internal to the cells
of the computational grid, the H\"older continuous first derivatives requirement on $\phi$ is satisfied, resulting in a stable numerical algorithm.
In addition, the surface of each cell is smooth in correspondence with the quadrature nodes. On such a geometry, the values for the free
coefficients on the double sided wake surface are set to $\ab=\boldsymbol{0}$ and $\bar{\bar{C}}=\bar{\bar{I}}$, where $\bar{\bar{I}}$ is the identity double
tensor. Thus, the Galerkin formulation results in the following algebraic linear system

\begin{equation}
A\hat{\Ub} = \bb^{U},
\end{equation}
where

\begin{itemize}
\item $\hat{\Ub}\in\mathbb{R}^{3 N_V}$ is the solution vector containing all the components of the potential gradient solution on all the BEM collocation points;  
\item $A\in\mathbb{R}^{3 N_V \times 3 N_V}:\ A_{ij} = \int_{\Gamma}\xi_i(\xb)\xi_j(\xb)\d s$ is a sparse mass matrix;
\item $\bb^{U}\in\mathbb{R}^{3 N_V}:\ b^{U}_i=\int_{\Gamma}\gb(\xb)\cdot \xi(\xb)\!\d s_x$  is the right hand side vector.
\end{itemize}

Finally, the velocity on the BEM collocation points is obtained adding the asymptotic velocity to the gradient values contained in $\hat{\Ub}$.
The resulting velocity vectors are used to carry out the wake relaxation step. As previously mentioned, and as illustrated
in Figure \ref{fig:meshRef}, the computational grid is structured in the wake region. Thus, it is possible to arrange the wake collocation
points in lines oriented in the streamwise direction. Given a wake degree of freedom $iw$ in the BEM problem, and the degree of freedom $pw$ corresponding
to the previous collocation point on a pathline, the position of collocation point $iw$ will be computed as

\begin{equation}
\xb_{iw} = \xb_{pw}+\frac{\vb_i}{|\vb_i|}d \quad \forall \xb_i\in\Gamma_W
\label{eq:wake_relax}
\end{equation}
in which $d$ is a user prescribed wake cell length. The wake relaxation step in Equation \ref{eq:wake_relax} is carried out on all the wake
collocation points, leading to a modified wake geometry. The BEM nonlinear problem for the velocity potential is then solved making use of the
new wake shape, and the new fluid velocities are used to again modify the wake geometry. Such a fixed point iteration strategy is repeated until
convergence of the wake position and velocity potential, which in the test cases considered in Sections \ref{sec:convergence} and \ref{sec:results}
requires less than 10 iterations.

\section{Convergence tests}\label{sec:convergence}

This section presents and discusses the numerical results obtained with the lifting surface potential flow solver proposed. The numerical tests
have been designed to evaluate the accuracy and spatial convergence capability of both the nonlinear system solved to obtain the potential and
the wake relaxation algorithm.

\subsection{Test case setup and resolution}

\begin{figure}[htb!]
\centerline{
  \ifpdf
  \resizebox{0.7\textwidth}{!}{
    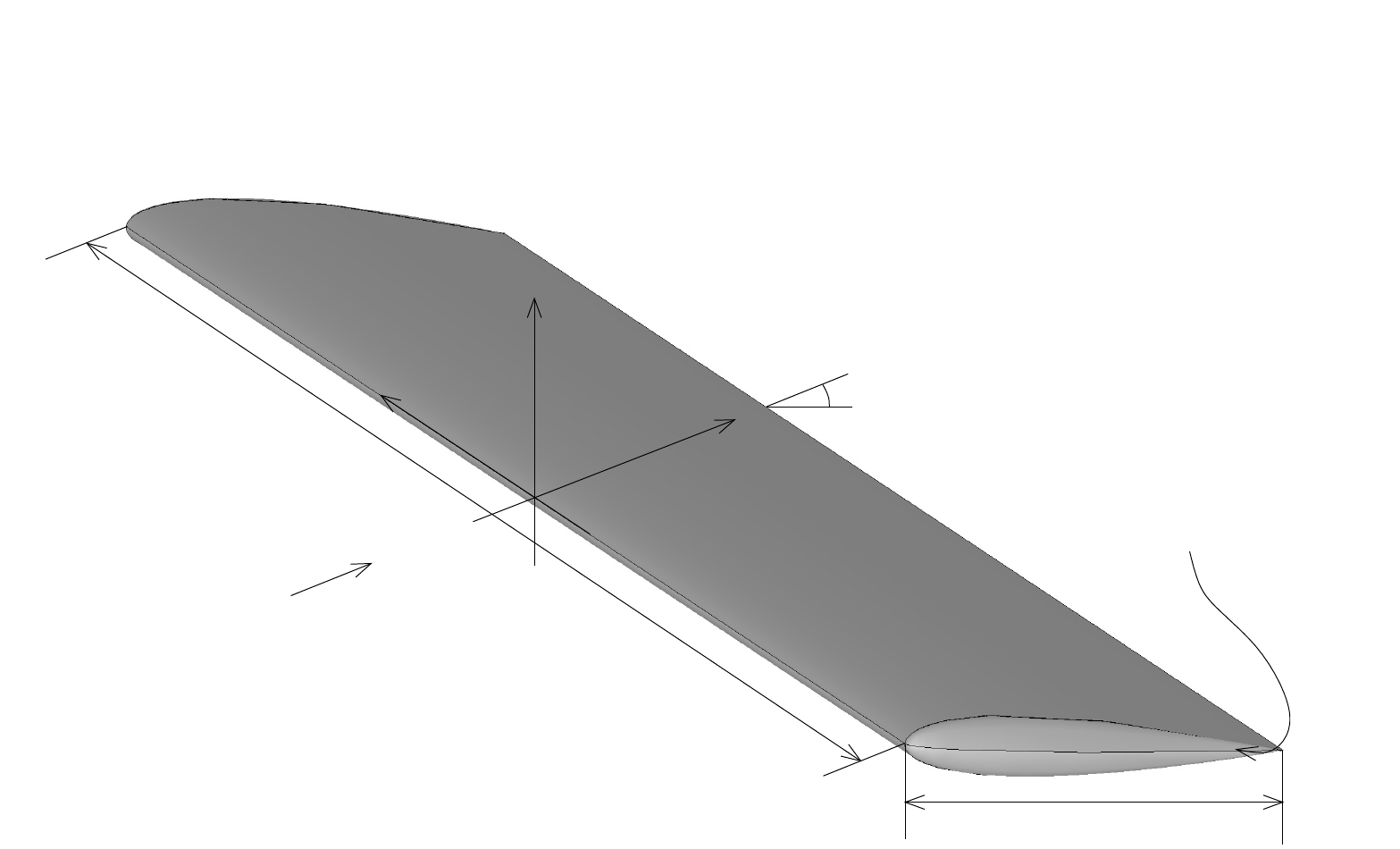
  }
  \else
  \resizebox{0.7\textwidth}{!}{
    \input{./figures/wing_setup.pstex_t}
  }
  \fi
}       	
\caption{A CAD rendering of the wing used in the first numerical tests, with relevant length measures.\label{fig:wingSetup}}
\end{figure}

The initial test case considered is the potential flow past a rectangular and constant-section wing. The airfoil section
considered is a NACA 0012 with a chord $c=1\ $m,
and the wing is generated translating such a shape for a span $b=4\ $m along a direction perpendicular to it. In the simulations, we make use of the right handed
reference frame described in Sec. \ref{sec:pert_pot_prob}  in which axis $x$ is aligned with the inflow velocity vector $\Vb_{\infty}=V_{\infty}\hat{\ib}$,
axis $z$ is aligned with the gravity acceleration vector $\gb=-g\hat{\kb}$, with $g=-9.81 m/s^2$ and axis $y$ normal to both $x$ and $z$. In such
framework, the NACA 0012
airfoil sections are parallel to the $xz$ plane. The wing span direction is then aligned with the $y$ axis, and the center of the reference
frame is located at the leading edge of the section corresponding with the symmetry plane of the wing, as illustrated in Figure \ref{fig:wingSetup}.
As the image shows, the angle of attack between the the chord $c$ and the $x$ axis is $\alpha=8.5^\circ$. In this particular
test the wing tip has been capped with a solid of revolution smoothly connected to the remaining part of the solid. As will be seen, this will help
showcasing the ability of the solver developed, to automatically refine the computational mesh on curved CAD surfaces. In fact, an extremely coarse initial
mesh composed by three rectangular cells (one covering the leeward side of the wing, one covering the windward side of the wing, and one covering the initially
flat wake --- having length 4$c$) is automatically refined to obtain the desired computational grid for each simulation. As
discussed in~\cite{pi-BEM}, the automated mesh generation process consists of an initial phase in which the aspect ratio of each cell is limited to a
user defined maximum value, followed by user specified numbers of uniform refinement cycles and adaptive refinements based on the local underlying
CAD surface curvature. 

\begin{figure}[htb!]
\begin{center}
\includegraphics[width=.48\textwidth]{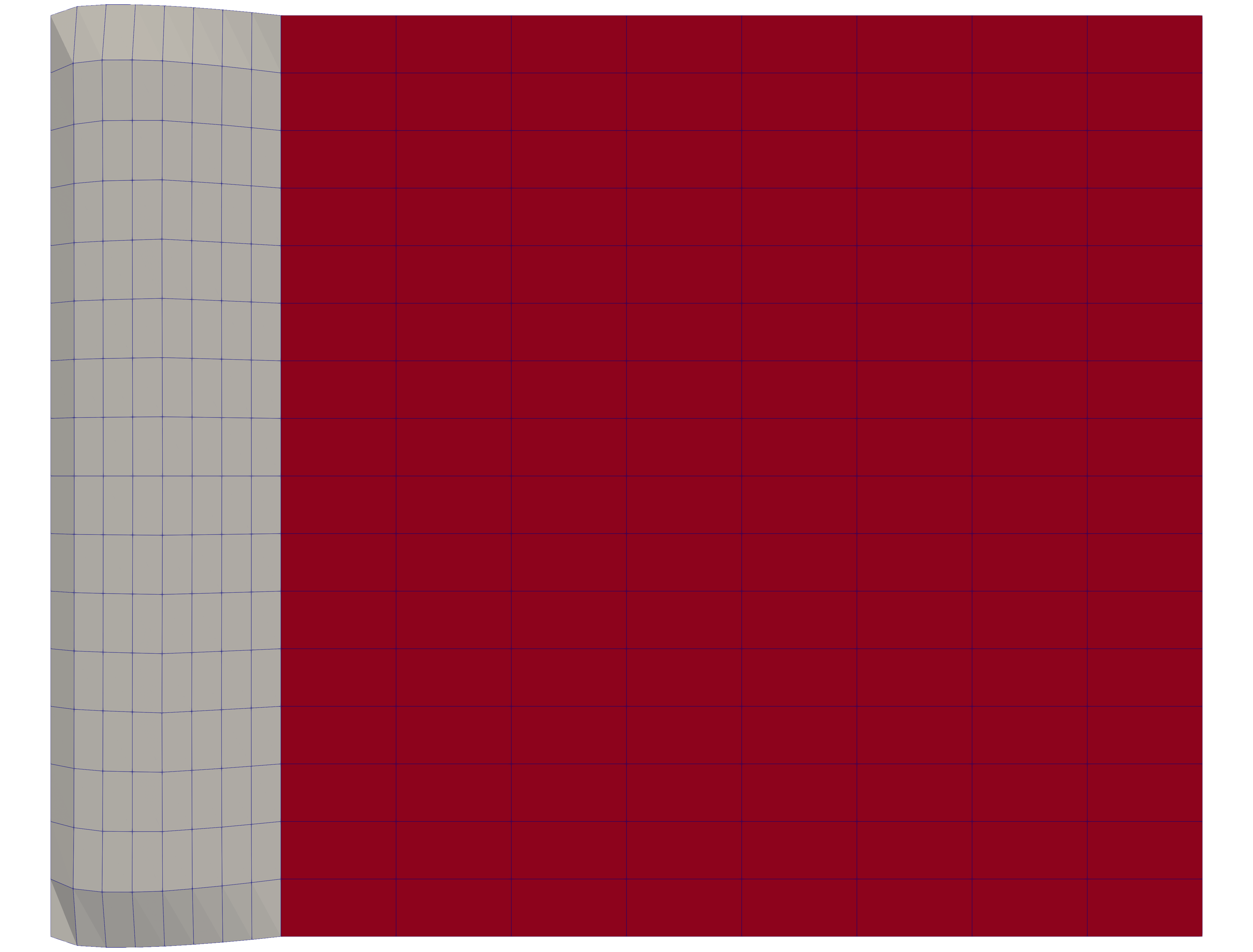}
\includegraphics[width=.48\textwidth]{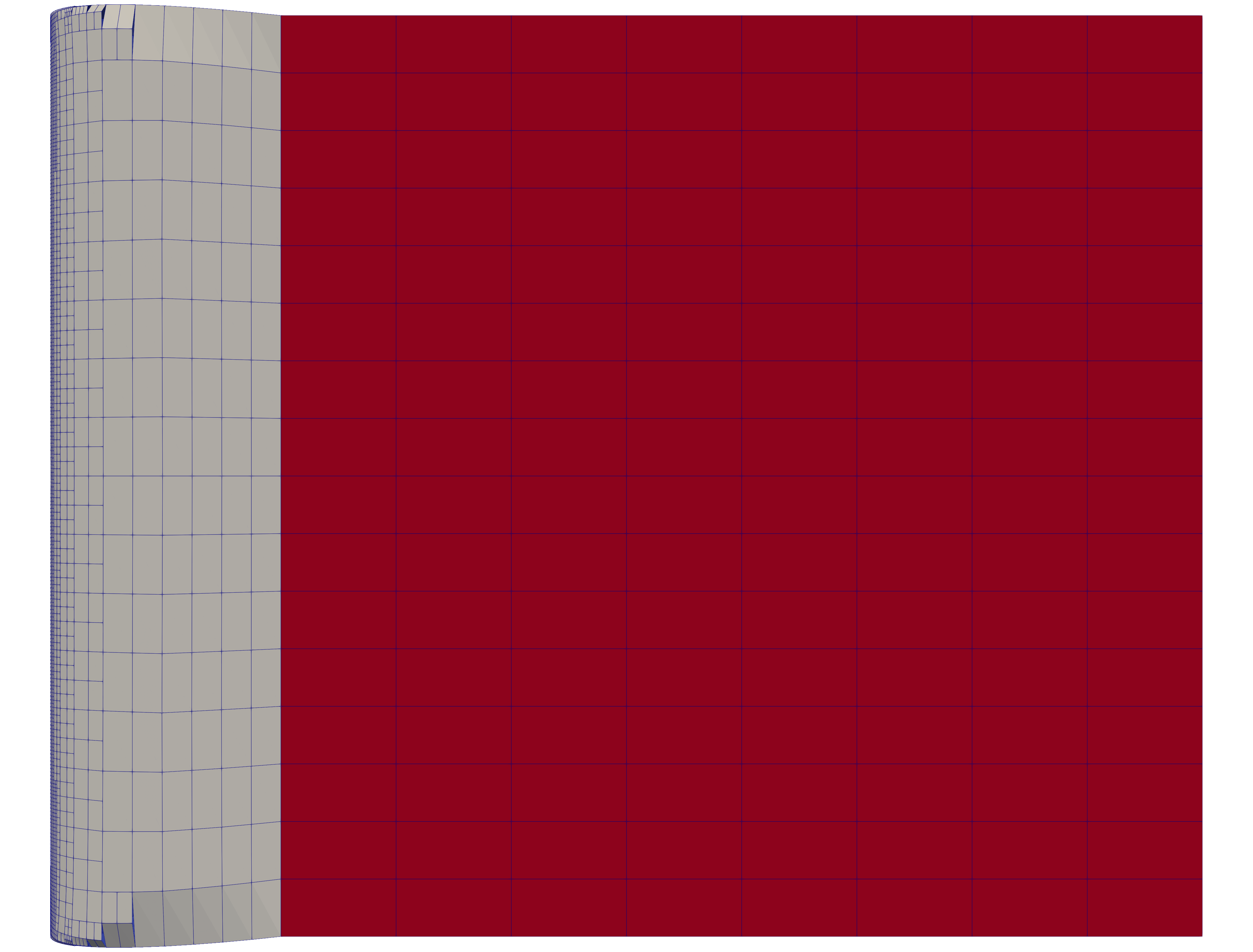}
\end{center}
\caption{A top view of the automatically generated initial mesh on the wing and wake grid. On the left, the --- structured --- grid obtained setting maximum aspect
         ratio 2.5, and 3 uniform refinement cycles. On the right, the --- unstructured --- grid obtained setting maximum aspect
         ratio 2.5, 3 uniform refinement cycles, and 4 adaptive refinement cycles based on CAD curvature. Note that, regardless of the wing grid type, the
         wake mesh is structured. The interface between the two regions, occurring at the trailing edge, is conformal.  
         \label{fig:meshRef}}
\end{figure}

Examples of the initial grids obtained with
the automated mesh generation algorithm are presented in Figure \ref{fig:meshRef}. The grid shown on the left, in particular, accounts for a total of 384 cells.
Such a grid results in 459, 1683 and 3675 degrees of freedom when bi-linear, bi-quadratic and bi-cubic finite elements are selected, respectively.
The additional collocation points required by the high order finite elements are automatically placed on the curved CAD surfaces at the moment of the degrees
of freedom distribution. Thus, the switch from low to high order finite elements does not require the generation of specific computational grids. 
For reference, the solution of the nonlinear problem with
459 degrees of freedom requires few seconds --- including jacobian matrix assembling --- on a single processor of the Intel Quad Core i7-7700HQ 2.80GHz,
32 GB RAM laptop used
to produce all the results shown in this work. After this, the velocity computation combined with the wake relaxation step, requires up to 10 seconds. Thus,
the execution time required by 10 wake relaxation cycles is approximately 2 minutes. Scaling up to the largest grids considered in this work, which feature up to
5500 degrees of freedom, the resolution of the nonlinear problem requires approximately a minute, while the computation of the velocity on the wake nodes
can take as much as 10 minutes. In such case, the required computation time for the full wake relaxation can reach more than an hour. Things can be made even
worse when high order finite elements are considered, as they typically demand higher quadrature order to obtain accurate solutions. Clearly, the
computational cost highlighted, can be significantly reduced through parallelization of both the nonlinear system resolution and the wake
velocity computation algorithms. For such a reason, a parallel version of the solver is currently being implemented.

\subsection{Nonlinear system resolution accuracy and convergence}

The first result presented aims at confirming that the algorithm developed is correctly solving the nonlinear system described in Section \ref{sec:num_kutta}.
Thus, the wake geometry is kept flat in these first numerical tests presented.
Clearly, no analytic solution is available for the realistic --- and non analytic --- wing geometry selected. For such a reason, the convergence of the resolution
algorithm developed is assessed observing the numerical solutions behavior when increasingly refined grids or higher order discretizations are considered.

\subsubsection{Potential jump at trailing edge}

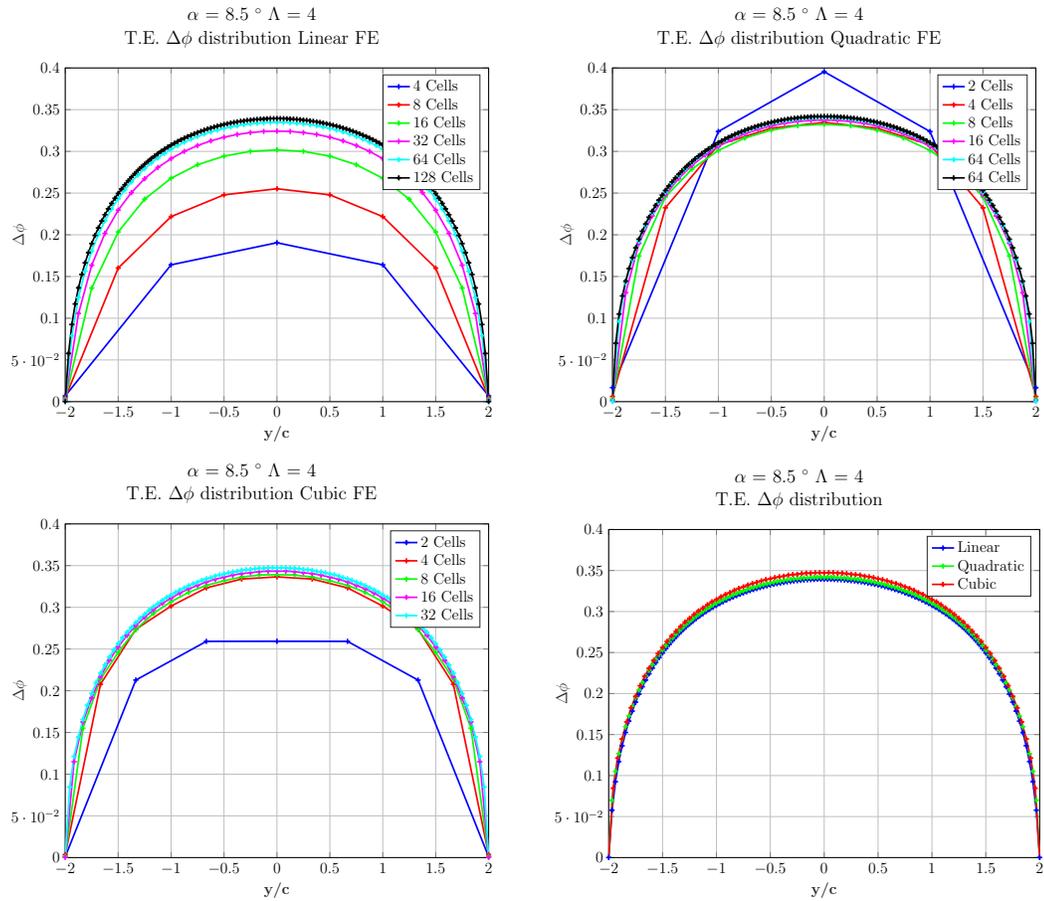
\begin{figure}[htb!]
\begin{tabular}{c c}
\resizebox{0.5\linewidth}{!}{
\centerline{
%
\definecolor{mycolor1}{rgb}{1.00000,0.00000,1.00000}%
\definecolor{mycolor2}{rgb}{0.00000,1.00000,1.00000}%
\definecolor{mycolor3}{rgb}{0.00000,0.00000,0.00000}%
\begin{tikzpicture}

\begin{axis}[%
width=4.521in,
height=3.566in,
at={(0.758in,0.481in)},
scale only axis,
xmin=-2,
xmax=2,
xlabel style={font=\color{white!15!black}},
xlabel={\textbf{y/c}},
ymin=0,
ymax=0.4,
ylabel style={font=\color{white!15!black}},
ylabel={\textbf{$\Delta\phi$}},
axis background/.style={fill=white},
xmajorgrids,
ymajorgrids,
legend style={legend cell align=left, align=left, draw=white!15!black}
]
\addplot [color=blue, line width=1.3pt, mark=+, mark options={solid, blue}]
  table[row sep=crcr]{%
  -2.0000   0.0069 \\
  -1.0000   0.1640 \\
        0   0.1905 \\
   1.0000   0.1640 \\
   2.0000   0.0069 \\
};
\addlegendentry{4 Cells}

\addplot [color=red, line width=1.3pt, mark=+, mark options={solid, red}]
  table[row sep=crcr]{%
  -2.0000   0.0047 \\
  -1.5000   0.1601 \\
  -1.0000   0.2218 \\
  -0.5000   0.2479 \\
        0   0.2552 \\
   0.5000   0.2479 \\
   1.0000   0.2218 \\
   1.5000   0.1601 \\
   2.0000   0.0047 \\
};
\addlegendentry{8 Cells}

\addplot [color=green, line width=1.3pt, mark=+, mark options={solid, green}]
  table[row sep=crcr]{%
  -2.0000   0.0032 \\
  -1.7500   0.1362 \\
  -1.5000   0.2033 \\
  -1.2500   0.2427 \\
  -1.0000   0.2678 \\
  -0.7500   0.2840 \\
  -0.5000   0.2943 \\
  -0.2500   0.3000 \\
        0   0.3018 \\
   0.2500   0.3000 \\
   0.5000   0.2943 \\
   0.7500   0.2840 \\
   1.0000   0.2678 \\
   1.2500   0.2427 \\
   1.5000   0.2033 \\
   1.7500   0.1362 \\
   2.0000   0.0032 \\
};
\addlegendentry{16 Cells}

\addplot [color=mycolor1, line width=1.3pt, mark=+, mark options={solid, mycolor1}]
  table[row sep=crcr]{%
  -2.0000   0.0019 \\
  -1.8750   0.1057 \\
  -1.7500   0.1636 \\
  -1.6250   0.2017 \\
  -1.5000   0.2295 \\
  -1.3750   0.2506 \\
  -1.2500   0.2673 \\
  -1.1250   0.2806 \\
  -1.0000   0.2913 \\
  -0.8750   0.3000 \\
  -0.7500   0.3070 \\
  -0.6250   0.3126 \\
  -0.5000   0.3170 \\
  -0.3750   0.3203 \\
  -0.2500   0.3225 \\
  -0.1250   0.3239 \\
        0   0.3243 \\
   0.1250   0.3239 \\
   0.2500   0.3225 \\
   0.3750   0.3203 \\
   0.5000   0.3170 \\
   0.6250   0.3126 \\
   0.7500   0.3070 \\
   0.8750   0.3000 \\
   1.0000   0.2913 \\
   1.1250   0.2806 \\
   1.2500   0.2673 \\
   1.3750   0.2506 \\
   1.5000   0.2295 \\
   1.6250   0.2017 \\
   1.7500   0.1636 \\
   1.8750   0.1057 \\
   2.0000   0.0019 \\
};
\addlegendentry{32 Cells}

\addplot [color=mycolor2, line width=1.3pt, mark=+, mark options={solid, mycolor2}]
  table[row sep=crcr]{%
  -2.0000   0.0008 \\
  -1.9375   0.0790 \\
  -1.8750   0.1247 \\
  -1.8125   0.1563 \\
  -1.7500   0.1806 \\
  -1.6875   0.2003 \\
  -1.6250   0.2168 \\
  -1.5625   0.2309 \\
  -1.5000   0.2431 \\
  -1.4375   0.2538 \\
  -1.3750   0.2633 \\
  -1.3125   0.2717 \\
  -1.2500   0.2793 \\
  -1.1875   0.2860 \\
  -1.1250   0.2921 \\
  -1.0625   0.2976 \\
  -1.0000   0.3025 \\
  -0.9375   0.3070 \\
  -0.8750   0.3110 \\
  -0.8125   0.3146 \\
  -0.7500   0.3178 \\
  -0.6875   0.3207 \\
  -0.6250   0.3233 \\
  -0.5625   0.3255 \\
  -0.5000   0.3275 \\
  -0.4375   0.3293 \\
  -0.3750   0.3307 \\
  -0.3125   0.3320 \\
  -0.2500   0.3330 \\
  -0.1875   0.3337 \\
  -0.1250   0.3343 \\
  -0.0625   0.3346 \\
        0   0.3347 \\
   0.0625   0.3346 \\
   0.1250   0.3343 \\
   0.1875   0.3337 \\
   0.2500   0.3330 \\
   0.3125   0.3320 \\
   0.3750   0.3307 \\
   0.4375   0.3293 \\
   0.5000   0.3275 \\
   0.5625   0.3255 \\
   0.6250   0.3233 \\
   0.6875   0.3207 \\
   0.7500   0.3178 \\
   0.8125   0.3146 \\
   0.8750   0.3110 \\
   0.9375   0.3070 \\
   1.0000   0.3025 \\
   1.0625   0.2976 \\
   1.1250   0.2921 \\
   1.1875   0.2860 \\
   1.2500   0.2793 \\
   1.3125   0.2717 \\
   1.3750   0.2633 \\
   1.4375   0.2538 \\
   1.5000   0.2431 \\
   1.5625   0.2309 \\
   1.6250   0.2168 \\
   1.6875   0.2003 \\
   1.7500   0.1806 \\
   1.8125   0.1563 \\
   1.8750   0.1247 \\
   1.9375   0.0790 \\
   2.0000   0.0008 \\
};
\addlegendentry{64 Cells}

\addplot [color=mycolor3, line width=1.3pt, mark=+, mark options={solid, mycolor3}]
  table[row sep=crcr]{%
  -2.0000   0.0001 \\
  -1.9688   0.0578 \\
  -1.9375   0.0925 \\
  -1.9062   0.1171 \\
  -1.8750   0.1364 \\
  -1.8438   0.1525 \\
  -1.8125   0.1664 \\
  -1.7812   0.1786 \\
  -1.7500   0.1895 \\
  -1.7188   0.1993 \\
  -1.6875   0.2083 \\
  -1.6562   0.2165 \\
  -1.6250   0.2241 \\
  -1.5938   0.2311 \\
  -1.5625   0.2377 \\
  -1.5312   0.2438 \\
  -1.5000   0.2495 \\
  -1.4688   0.2549 \\
  -1.4375   0.2599 \\
  -1.4062   0.2647 \\
  -1.3750   0.2691 \\
  -1.3438   0.2734 \\
  -1.3125   0.2774 \\
  -1.2812   0.2812 \\
  -1.2500   0.2847 \\
  -1.2188   0.2881 \\
  -1.1875   0.2914 \\
  -1.1562   0.2944 \\
  -1.1250   0.2973 \\
  -1.0938   0.3001 \\
  -1.0625   0.3027 \\
  -1.0312   0.3052 \\
  -1.0000   0.3076 \\
  -0.9688   0.3098 \\
  -0.9375   0.3119 \\
  -0.9062   0.3140 \\
  -0.8750   0.3159 \\
  -0.8438   0.3177 \\
  -0.8125   0.3194 \\
  -0.7812   0.3211 \\
  -0.7500   0.3226 \\
  -0.7188   0.3241 \\
  -0.6875   0.3255 \\
  -0.6562   0.3268 \\
  -0.6250   0.3280 \\
  -0.5938   0.3292 \\
  -0.5625   0.3303 \\
  -0.5312   0.3313 \\
  -0.5000   0.3323 \\
  -0.4688   0.3331 \\
  -0.4375   0.3340 \\
  -0.4062   0.3347 \\
  -0.3750   0.3354 \\
  -0.3438   0.3361 \\
  -0.3125   0.3366 \\
  -0.2812   0.3372 \\
  -0.2500   0.3376 \\
  -0.2188   0.3380 \\
  -0.1875   0.3384 \\
  -0.1562   0.3387 \\
  -0.1250   0.3389 \\
  -0.0938   0.3391 \\
  -0.0625   0.3393 \\
  -0.0312   0.3393 \\
        0   0.3394 \\
   0.0312   0.3393 \\
   0.0625   0.3393 \\
   0.0938   0.3391 \\
   0.1250   0.3389 \\
   0.1562   0.3387 \\
   0.1875   0.3384 \\
   0.2188   0.3380 \\
   0.2500   0.3376 \\
   0.2812   0.3372 \\
   0.3125   0.3366 \\
   0.3438   0.3361 \\
   0.3750   0.3354 \\
   0.4062   0.3347 \\
   0.4375   0.3340 \\
   0.4688   0.3331 \\
   0.5000   0.3323 \\
   0.5312   0.3313 \\
   0.5625   0.3303 \\
   0.5938   0.3292 \\
   0.6250   0.3280 \\
   0.6562   0.3268 \\
   0.6875   0.3255 \\
   0.7188   0.3241 \\
   0.7500   0.3226 \\
   0.7812   0.3211 \\
   0.8125   0.3194 \\
   0.8438   0.3177 \\
   0.8750   0.3159 \\
   0.9062   0.3140 \\
   0.9375   0.3119 \\
   0.9688   0.3098 \\
   1.0000   0.3076 \\
   1.0312   0.3052 \\
   1.0625   0.3027 \\
   1.0938   0.3001 \\
   1.1250   0.2973 \\
   1.1562   0.2944 \\
   1.1875   0.2914 \\
   1.2188   0.2881 \\
   1.2500   0.2847 \\
   1.2812   0.2812 \\
   1.3125   0.2774 \\
   1.3438   0.2734 \\
   1.3750   0.2691 \\
   1.4062   0.2647 \\
   1.4375   0.2599 \\
   1.4688   0.2549 \\
   1.5000   0.2495 \\
   1.5312   0.2438 \\
   1.5625   0.2377 \\
   1.5938   0.2311 \\
   1.6250   0.2241 \\
   1.6562   0.2165 \\
   1.6875   0.2083 \\
   1.7188   0.1993 \\
   1.7500   0.1895 \\
   1.7812   0.1786 \\
   1.8125   0.1664 \\
   1.8438   0.1525 \\
   1.8750   0.1364 \\
   1.9062   0.1171 \\
   1.9375   0.0925 \\
   1.9688   0.0578 \\
   2.0000   0.0001 \\
};
\addlegendentry{128 Cells}

\end{axis}

\node[align=center,font=\large, yshift=2em] (title) 
    at (current bounding box.north)
    {$\alpha$ = 8.5 ${^\circ}$ $\Lambda$ = 4\\ T.E. $\Delta\phi$ distribution Linear FE};

\end{tikzpicture}%
}    
}
&
\resizebox{0.5\linewidth}{!}{
\centerline{
%
%
\definecolor{mycolor1}{rgb}{1.00000,0.00000,1.00000}%
\definecolor{mycolor2}{rgb}{0.00000,1.00000,1.00000}%
\definecolor{mycolor3}{rgb}{0.00000,0.00000,0.00000}%
\begin{tikzpicture}

\begin{axis}[%
width=4.521in,
height=3.566in,
at={(0.758in,0.481in)},
scale only axis,
xmin=-2,
xmax=2,
xlabel style={font=\color{white!15!black}},
xlabel={\textbf{y/c}},
ymin=0,
ymax=0.4,
ylabel style={font=\color{white!15!black}},
ylabel={\textbf{$\Delta\phi$}},
axis background/.style={fill=white},
xmajorgrids,
ymajorgrids,
legend style={legend cell align=left, align=left, draw=white!15!black}
]
\addplot [color=blue, line width=1.3pt, mark=+, mark options={solid, blue}]
  table[row sep=crcr]{%
 -2.0000   0.0168 \\
  -1.0000   0.3238 \\
        0   0.3955 \\
   1.0000   0.3238 \\
   2.0000   0.0168 \\
};
\addlegendentry{2 Cells}

\addplot [color=red, line width=1.3pt, mark=+, mark options={solid, red}]
  table[row sep=crcr]{%
  -2.0000   0.0062 \\
  -1.5000   0.2324 \\
  -1.0000   0.3077 \\
  -0.5000   0.3278 \\
        0   0.3347 \\
   0.5000   0.3278 \\
   1.0000   0.3077 \\
   1.5000   0.2324 \\
   2.0000   0.0062 \\
};
\addlegendentry{4 Cells}

\addplot [color=green, line width=1.3pt, mark=+, mark options={solid, green}]
  table[row sep=crcr]{%
  -2.0000   0.0030 \\
  -1.7500   0.1746 \\
  -1.5000   0.2455 \\
  -1.2500   0.2787 \\
  -1.0000   0.3012 \\
  -0.7500   0.3160 \\
  -0.5000   0.3255 \\
  -0.2500   0.3308 \\
        0   0.3325 \\
   0.2500   0.3308 \\
   0.5000   0.3255 \\
   0.7500   0.3160 \\
   1.0000   0.3012 \\
   1.2500   0.2787 \\
   1.5000   0.2455 \\
   1.7500   0.1746 \\
   2.0000   0.0030 \\
};
\addlegendentry{8 Cells}

\addplot [color=mycolor1, line width=1.3pt, mark=+, mark options={solid, mycolor1}]
  table[row sep=crcr]{%
  -2.0000   0.0013 \\
  -1.8750   0.1305 \\
  -1.7500   0.1898 \\
  -1.6250   0.2227 \\
  -1.5000   0.2479 \\
  -1.3750   0.2675 \\
  -1.2500   0.2831 \\
  -1.1250   0.2956 \\
  -1.0000   0.3058 \\
  -0.8750   0.3141 \\
  -0.7500   0.3209 \\
  -0.6250   0.3263 \\
  -0.5000   0.3305 \\
  -0.3750   0.3336 \\
  -0.2500   0.3358 \\
  -0.1250   0.3371 \\
        0   0.3375 \\
   0.1250   0.3371 \\
   0.2500   0.3358 \\
   0.3750   0.3336 \\
   0.5000   0.3305 \\
   0.6250   0.3263 \\
   0.7500   0.3209 \\
   0.8750   0.3141 \\
   1.0000   0.3058 \\
   1.1250   0.2956 \\
   1.2500   0.2831 \\
   1.3750   0.2675 \\
   1.5000   0.2479 \\
   1.6250   0.2227 \\
   1.7500   0.1898 \\
   1.8750   0.1305 \\
   2.0000   0.0013 \\
};
\addlegendentry{16 Cells}

\addplot [color=mycolor2, line width=1.3pt, mark=+, mark options={solid, mycolor2}]
  table[row sep=crcr]{%
  -2.0000   0.0003 \\
  -1.9375   0.0961 \\
  -1.8750   0.1424 \\
  -1.8125   0.1701 \\
  -1.7500   0.1924 \\
  -1.6875   0.2107 \\
  -1.6250   0.2263 \\
  -1.5625   0.2396 \\
  -1.5000   0.2513 \\
  -1.4375   0.2616 \\
  -1.3750   0.2708 \\
  -1.3125   0.2789 \\
  -1.2500   0.2863 \\
  -1.1875   0.2928 \\
  -1.1250   0.2987 \\
  -1.0625   0.3041 \\
  -1.0000   0.3089 \\
  -0.9375   0.3133 \\
  -0.8750   0.3172 \\
  -0.8125   0.3207 \\
  -0.7500   0.3239 \\
  -0.6875   0.3267 \\
  -0.6250   0.3293 \\
  -0.5625   0.3315 \\
  -0.5000   0.3335 \\
  -0.4375   0.3352 \\
  -0.3750   0.3366 \\
  -0.3125   0.3378 \\
  -0.2500   0.3388 \\
  -0.1875   0.3396 \\
  -0.1250   0.3401 \\
  -0.0625   0.3405 \\
        0   0.3406 \\
   0.0625   0.3405 \\
   0.1250   0.3401 \\
   0.1875   0.3396 \\
   0.2500   0.3388 \\
   0.3125   0.3378 \\
   0.3750   0.3366 \\
   0.4375   0.3352 \\
   0.5000   0.3335 \\
   0.5625   0.3315 \\
   0.6250   0.3293 \\
   0.6875   0.3267 \\
   0.7500   0.3239 \\
   0.8125   0.3207 \\
   0.8750   0.3172 \\
   0.9375   0.3133 \\
   1.0000   0.3089 \\
   1.0625   0.3041 \\
   1.1250   0.2987 \\
   1.1875   0.2928 \\
   1.2500   0.2863 \\
   1.3125   0.2789 \\
   1.3750   0.2708 \\
   1.4375   0.2616 \\
   1.5000   0.2513 \\
   1.5625   0.2396 \\
   1.6250   0.2263 \\
   1.6875   0.2107 \\
   1.7500   0.1924 \\
   1.8125   0.1701 \\
   1.8750   0.1424 \\
   1.9375   0.0961 \\
   2.0000   0.0003 \\
};
\addlegendentry{64 Cells}

\addplot [color=mycolor3, line width=1.3pt, mark=+, mark options={solid, mycolor3}]
  table[row sep=crcr]{%
  -2.0000  -0.0004 \\
  -1.9688   0.0700 \\
  -1.9375   0.1051 \\
  -1.9062   0.1267 \\
  -1.8750   0.1445 \\
  -1.8438   0.1596 \\
  -1.8125   0.1727 \\
  -1.7812   0.1844 \\
  -1.7500   0.1948 \\
  -1.7188   0.2043 \\
  -1.6875   0.2130 \\
  -1.6562   0.2210 \\
  -1.6250   0.2283 \\
  -1.5938   0.2352 \\
  -1.5625   0.2416 \\
  -1.5312   0.2476 \\
  -1.5000   0.2532 \\
  -1.4688   0.2584 \\
  -1.4375   0.2634 \\
  -1.4062   0.2681 \\
  -1.3750   0.2725 \\
  -1.3438   0.2766 \\
  -1.3125   0.2806 \\
  -1.2812   0.2843 \\
  -1.2500   0.2878 \\
  -1.2188   0.2912 \\
  -1.1875   0.2944 \\
  -1.1562   0.2974 \\
  -1.1250   0.3003 \\
  -1.0938   0.3030 \\
  -1.0625   0.3056 \\
  -1.0312   0.3081 \\
  -1.0000   0.3104 \\
  -0.9688   0.3126 \\
  -0.9375   0.3147 \\
  -0.9062   0.3167 \\
  -0.8750   0.3187 \\
  -0.8438   0.3205 \\
  -0.8125   0.3222 \\
  -0.7812   0.3238 \\
  -0.7500   0.3253 \\
  -0.7188   0.3268 \\
  -0.6875   0.3282 \\
  -0.6562   0.3295 \\
  -0.6250   0.3307 \\
  -0.5938   0.3319 \\
  -0.5625   0.3329 \\
  -0.5312   0.3340 \\
  -0.5000   0.3349 \\
  -0.4688   0.3358 \\
  -0.4375   0.3366 \\
  -0.4062   0.3374 \\
  -0.3750   0.3381 \\
  -0.3438   0.3387 \\
  -0.3125   0.3393 \\
  -0.2812   0.3398 \\
  -0.2500   0.3402 \\
  -0.2188   0.3407 \\
  -0.1875   0.3410 \\
  -0.1562   0.3413 \\
  -0.1250   0.3415 \\
  -0.0938   0.3417 \\
  -0.0625   0.3419 \\
  -0.0312   0.3419 \\
        0   0.3420 \\
   0.0312   0.3419 \\
   0.0625   0.3419 \\
   0.0938   0.3417 \\
   0.1250   0.3415 \\
   0.1562   0.3413 \\
   0.1875   0.3410 \\
   0.2188   0.3407 \\
   0.2500   0.3402 \\
   0.2812   0.3398 \\
   0.3125   0.3393 \\
   0.3438   0.3387 \\
   0.3750   0.3381 \\
   0.4062   0.3374 \\
   0.4375   0.3366 \\
   0.4688   0.3358 \\
   0.5000   0.3349 \\
   0.5312   0.3340 \\
   0.5625   0.3329 \\
   0.5938   0.3319 \\
   0.6250   0.3307 \\
   0.6562   0.3295 \\
   0.6875   0.3282 \\
   0.7188   0.3268 \\
   0.7500   0.3253 \\
   0.7812   0.3238 \\
   0.8125   0.3222 \\
   0.8438   0.3205 \\
   0.8750   0.3187 \\
   0.9062   0.3167 \\
   0.9375   0.3147 \\
   0.9688   0.3126 \\
   1.0000   0.3104 \\
   1.0312   0.3081 \\
   1.0625   0.3056 \\
   1.0938   0.3030 \\
   1.1250   0.3003 \\
   1.1562   0.2974 \\
   1.1875   0.2944 \\
   1.2188   0.2912 \\
   1.2500   0.2878 \\
   1.2812   0.2843 \\
   1.3125   0.2806 \\
   1.3438   0.2766 \\
   1.3750   0.2725 \\
   1.4062   0.2681 \\
   1.4375   0.2634 \\
   1.4688   0.2584 \\
   1.5000   0.2532 \\
   1.5312   0.2476 \\
   1.5625   0.2416 \\
   1.5938   0.2352 \\
   1.6250   0.2283 \\
   1.6562   0.2210 \\
   1.6875   0.2130 \\
   1.7188   0.2043 \\
   1.7500   0.1948 \\
   1.7812   0.1844 \\
   1.8125   0.1727 \\
   1.8438   0.1596 \\
   1.8750   0.1445 \\
   1.9062   0.1267 \\
   1.9375   0.1051 \\
   1.9688   0.0700 \\
   2.0000  -0.0004 \\
};
\addlegendentry{64 Cells}

\end{axis}

\node[align=center,font=\large, yshift=2em] (title) 
    at (current bounding box.north)
    {$\alpha$ = 8.5 ${^\circ}$ $\Lambda$ = 4\\ T.E. $\Delta\phi$ distribution Quadratic FE};
    
\end{tikzpicture}%
}    
}
\\
\resizebox{0.5\linewidth}{!}{
\centerline{
%
%
\definecolor{mycolor1}{rgb}{1.00000,0.00000,1.00000}%
\definecolor{mycolor2}{rgb}{0.00000,1.00000,1.00000}%

\begin{tikzpicture}

\begin{axis}[%
width=4.521in,
height=3.566in,
at={(0.758in,0.481in)},
scale only axis,
xmin=-2,
xmax=2,
xlabel style={font=\color{white!15!black}},
xlabel={\textbf{y/c}},
ymin=0,
ymax=0.4,
ylabel style={font=\color{white!15!black}},
ylabel={\textbf{$\Delta\phi$}},
axis background/.style={fill=white},
xmajorgrids,
ymajorgrids,
legend style={legend cell align=left, align=left, draw=white!15!black}
]
\addplot [color=blue, line width=1.3pt, mark=+, mark options={solid, blue}]
  table[row sep=crcr]{%
  -2.0000   0.0010 \\
  -1.3333   0.2128 \\
  -0.6667   0.2592 \\
        0   0.2593 \\
   0.6667   0.2592 \\
   1.3333   0.2128 \\
   2.0000   0.0010 \\
};
\addlegendentry{2 Cells}

\addplot [color=red, line width=1.3pt, mark=+, mark options={solid, red}]
  table[row sep=crcr]{%
  -2.0000   0.0036 \\
  -1.6667   0.2078 \\
  -1.3333   0.2737 \\
  -1.0000   0.3015 \\
  -0.6667   0.3233 \\
  -0.3333   0.3339 \\
        0   0.3365 \\
   0.3333   0.3339 \\
   0.6667   0.3233 \\
   1.0000   0.3015 \\
   1.3333   0.2737 \\
   1.6667   0.2078 \\
   2.0000   0.0036 \\
};
\addlegendentry{4 Cells}

\addplot [color=green, line width=1.3pt, mark=+, mark options={solid, green}]
  table[row sep=crcr]{%
  -2.0000   0.0014 \\
  -1.8333   0.1550 \\
  -1.6667   0.2139 \\
  -1.5000   0.2461 \\
  -1.3333   0.2733 \\
  -1.1667   0.2927 \\
  -1.0000   0.3067 \\
  -0.8333   0.3180 \\
  -0.6667   0.3263 \\
  -0.5000   0.3319 \\
  -0.3333   0.3363 \\
  -0.1667   0.3387 \\
        0   0.3392 \\
   0.1667   0.3387 \\
   0.3333   0.3363 \\
   0.5000   0.3319 \\
   0.6667   0.3263 \\
   0.8333   0.3180 \\
   1.0000   0.3067 \\
   1.1667   0.2927 \\
   1.3333   0.2733 \\
   1.5000   0.2461 \\
   1.6667   0.2139 \\
   1.8333   0.1550 \\
   2.0000   0.0014 \\
};
\addlegendentry{8 Cells}

\addplot [color=mycolor1, line width=1.3pt, mark=+, mark options={solid, mycolor1}]
  table[row sep=crcr]{%
  -2.0000   0.0003 \\
  -1.9167   0.1149 \\
  -1.8333   0.1626 \\
  -1.7500   0.1913 \\
  -1.6667   0.2167 \\
  -1.5833   0.2365 \\
  -1.5000   0.2524 \\
  -1.4167   0.2664 \\
  -1.3333   0.2782 \\
  -1.2500   0.2880 \\
  -1.1667   0.2970 \\
  -1.0833   0.3047 \\
  -1.0000   0.3110 \\
  -0.9167   0.3171 \\
  -0.8333   0.3222 \\
  -0.7500   0.3263 \\
  -0.6667   0.3304 \\
  -0.5833   0.3336 \\
  -0.5000   0.3361 \\
  -0.4167   0.3386 \\
  -0.3333   0.3404 \\
  -0.2500   0.3415 \\
  -0.1667   0.3428 \\
  -0.0833   0.3434 \\
        0   0.3433 \\
   0.0833   0.3434 \\
   0.1667   0.3428 \\
   0.2500   0.3415 \\
   0.3333   0.3404 \\
   0.4167   0.3386 \\
   0.5000   0.3361 \\
   0.5833   0.3336 \\
   0.6667   0.3304 \\
   0.7500   0.3263 \\
   0.8333   0.3222 \\
   0.9167   0.3171 \\
   1.0000   0.3110 \\
   1.0833   0.3047 \\
   1.1667   0.2970 \\
   1.2500   0.2880 \\
   1.3333   0.2782 \\
   1.4167   0.2664 \\
   1.5000   0.2524 \\
   1.5833   0.2365 \\
   1.6667   0.2167 \\
   1.7500   0.1913 \\
   1.8333   0.1626 \\
   1.9167   0.1149 \\
   2.0000   0.0003 \\
};
\addlegendentry{16 Cells}

\addplot [color=mycolor2, line width=1.3pt, mark=+, mark options={solid, mycolor2}]
  table[row sep=crcr]{%
  -2.0000  -0.0004 \\
  -1.9583   0.0846 \\
  -1.9167   0.1213 \\
  -1.8750   0.1445 \\
  -1.8333   0.1654 \\
  -1.7917   0.1825 \\
  -1.7500   0.1967 \\
  -1.7083   0.2096 \\
  -1.6667   0.2210 \\
  -1.6250   0.2309 \\
  -1.5833   0.2404 \\
  -1.5417   0.2489 \\
  -1.5000   0.2563 \\
  -1.4583   0.2636 \\
  -1.4167   0.2702 \\
  -1.3750   0.2760 \\
  -1.3333   0.2818 \\
  -1.2917   0.2872 \\
  -1.2500   0.2917 \\
  -1.2083   0.2965 \\
  -1.1667   0.3008 \\
  -1.1250   0.3044 \\
  -1.0833   0.3084 \\
  -1.0417   0.3119 \\
  -1.0000   0.3147 \\
  -0.9583   0.3180 \\
  -0.9167   0.3209 \\
  -0.8750   0.3232 \\
  -0.8333   0.3260 \\
  -0.7917   0.3283 \\
  -0.7500   0.3301 \\
  -0.7083   0.3324 \\
  -0.6667   0.3342 \\
  -0.6250   0.3356 \\
  -0.5833   0.3375 \\
  -0.5417   0.3390 \\
  -0.5000   0.3399 \\
  -0.4583   0.3414 \\
  -0.4167   0.3425 \\
  -0.3750   0.3431 \\
  -0.3333   0.3444 \\
  -0.2917   0.3451 \\
  -0.2500   0.3454 \\
  -0.2083   0.3463 \\
  -0.1667   0.3468 \\
  -0.1250   0.3467 \\
  -0.0833   0.3474 \\
  -0.0417   0.3475 \\
        0   0.3472 \\
   0.0417   0.3475 \\
   0.0833   0.3474 \\
   0.1250   0.3467 \\
   0.1667   0.3468 \\
   0.2083   0.3463 \\
   0.2500   0.3454 \\
   0.2917   0.3451 \\
   0.3333   0.3444 \\
   0.3750   0.3431 \\
   0.4167   0.3425 \\
   0.4583   0.3414 \\
   0.5000   0.3399 \\
   0.5417   0.3390 \\
   0.5833   0.3375 \\
   0.6250   0.3356 \\
   0.6667   0.3342 \\
   0.7083   0.3324 \\
   0.7500   0.3301 \\
   0.7917   0.3283 \\
   0.8333   0.3260 \\
   0.8750   0.3232 \\
   0.9167   0.3209 \\
   0.9583   0.3180 \\
   1.0000   0.3147 \\
   1.0417   0.3119 \\
   1.0833   0.3084 \\
   1.1250   0.3044 \\
   1.1667   0.3008 \\
   1.2083   0.2965 \\
   1.2500   0.2917 \\
   1.2917   0.2872 \\
   1.3333   0.2818 \\
   1.3750   0.2760 \\
   1.4167   0.2702 \\
   1.4583   0.2636 \\
   1.5000   0.2563 \\
   1.5417   0.2489 \\
   1.5833   0.2404 \\
   1.6250   0.2309 \\
   1.6667   0.2210 \\
   1.7083   0.2096 \\
   1.7500   0.1967 \\
   1.7917   0.1825 \\
   1.8333   0.1654 \\
   1.8750   0.1445 \\
   1.9167   0.1213 \\
   1.9583   0.0846 \\
   2.0000  -0.0004 \\
};
\addlegendentry{32 Cells}

\end{axis}

\node[align=center,font=\large, yshift=2em] (title) 
    at (current bounding box.north)
    {$\alpha$ = 8.5 ${^\circ}$ $\Lambda$ = 4\\ T.E. $\Delta\phi$ distribution Cubic FE};
    
\end{tikzpicture}%
}    
}
&
\resizebox{0.5\linewidth}{!}{
\centerline{
%
%
\begin{tikzpicture}

\begin{axis}[%
width=4.602in,
height=3.506in,
at={(0.772in,0.473in)},
scale only axis,
xmin=-2,
xmax=2,
xlabel style={font=\color{white!15!black}},
xlabel={\textbf{y/c}},
ymin=0,
ymax=0.4,
ylabel style={font=\color{white!15!black}},
ylabel={\textbf{$\Delta\phi$}},
axis background/.style={fill=white},
xmajorgrids,
ymajorgrids,
legend style={legend cell align=left, align=left, draw=white!15!black}
]
\addplot [color=blue, line width=1.3pt, mark=+, mark options={solid, blue}]
  table[row sep=crcr]{%
  -2.0000   0.0001 \\
  -1.9688   0.0578 \\
  -1.9375   0.0925 \\
  -1.9062   0.1171 \\
  -1.8750   0.1364 \\
  -1.8438   0.1525 \\
  -1.8125   0.1664 \\
  -1.7812   0.1786 \\
  -1.7500   0.1895 \\
  -1.7188   0.1993 \\
  -1.6875   0.2083 \\
  -1.6562   0.2165 \\
  -1.6250   0.2241 \\
  -1.5938   0.2311 \\
  -1.5625   0.2377 \\
  -1.5312   0.2438 \\
  -1.5000   0.2495 \\
  -1.4688   0.2549 \\
  -1.4375   0.2599 \\
  -1.4062   0.2647 \\
  -1.3750   0.2691 \\
  -1.3438   0.2734 \\
  -1.3125   0.2774 \\
  -1.2812   0.2812 \\
  -1.2500   0.2847 \\
  -1.2188   0.2881 \\
  -1.1875   0.2914 \\
  -1.1562   0.2944 \\
  -1.1250   0.2973 \\
  -1.0938   0.3001 \\
  -1.0625   0.3027 \\
  -1.0312   0.3052 \\
  -1.0000   0.3076 \\
  -0.9688   0.3098 \\
  -0.9375   0.3119 \\
  -0.9062   0.3140 \\
  -0.8750   0.3159 \\
  -0.8438   0.3177 \\
  -0.8125   0.3194 \\
  -0.7812   0.3211 \\
  -0.7500   0.3226 \\
  -0.7188   0.3241 \\
  -0.6875   0.3255 \\
  -0.6562   0.3268 \\
  -0.6250   0.3280 \\
  -0.5938   0.3292 \\
  -0.5625   0.3303 \\
  -0.5312   0.3313 \\
  -0.5000   0.3323 \\
  -0.4688   0.3331 \\
  -0.4375   0.3340 \\
  -0.4062   0.3347 \\
  -0.3750   0.3354 \\
  -0.3438   0.3361 \\
  -0.3125   0.3366 \\
  -0.2812   0.3372 \\
  -0.2500   0.3376 \\
  -0.2188   0.3380 \\
  -0.1875   0.3384 \\
  -0.1562   0.3387 \\
  -0.1250   0.3389 \\
  -0.0938   0.3391 \\
  -0.0625   0.3393 \\
  -0.0312   0.3393 \\
        0   0.3394 \\
   0.0312   0.3393 \\
   0.0625   0.3393 \\
   0.0938   0.3391 \\
   0.1250   0.3389 \\
   0.1562   0.3387 \\
   0.1875   0.3384 \\
   0.2188   0.3380 \\
   0.2500   0.3376 \\
   0.2812   0.3372 \\
   0.3125   0.3366 \\
   0.3438   0.3361 \\
   0.3750   0.3354 \\
   0.4062   0.3347 \\
   0.4375   0.3340 \\
   0.4688   0.3331 \\
   0.5000   0.3323 \\
   0.5312   0.3313 \\
   0.5625   0.3303 \\
   0.5938   0.3292 \\
   0.6250   0.3280 \\
   0.6562   0.3268 \\
   0.6875   0.3255 \\
   0.7188   0.3241 \\
   0.7500   0.3226 \\
   0.7812   0.3211 \\
   0.8125   0.3194 \\
   0.8438   0.3177 \\
   0.8750   0.3159 \\
   0.9062   0.3140 \\
   0.9375   0.3119 \\
   0.9688   0.3098 \\
   1.0000   0.3076 \\
   1.0312   0.3052 \\
   1.0625   0.3027 \\
   1.0938   0.3001 \\
   1.1250   0.2973 \\
   1.1562   0.2944 \\
   1.1875   0.2914 \\
   1.2188   0.2881 \\
   1.2500   0.2847 \\
   1.2812   0.2812 \\
   1.3125   0.2774 \\
   1.3438   0.2734 \\
   1.3750   0.2691 \\
   1.4062   0.2647 \\
   1.4375   0.2599 \\
   1.4688   0.2549 \\
   1.5000   0.2495 \\
   1.5312   0.2438 \\
   1.5625   0.2377 \\
   1.5938   0.2311 \\
   1.6250   0.2241 \\
   1.6562   0.2165 \\
   1.6875   0.2083 \\
   1.7188   0.1993 \\
   1.7500   0.1895 \\
   1.7812   0.1786 \\
   1.8125   0.1664 \\
   1.8438   0.1525 \\
   1.8750   0.1364 \\
   1.9062   0.1171 \\
   1.9375   0.0925 \\
   1.9688   0.0578 \\
   2.0000   0.0001 \\
};
\addlegendentry{Linear}

\addplot [color=green, line width=1.3pt, mark=+, mark options={solid, green}]
  table[row sep=crcr]{%
  -2.0000  -0.0004 \\
  -1.9688   0.0700 \\
  -1.9375   0.1051 \\
  -1.9062   0.1267 \\
  -1.8750   0.1445 \\
  -1.8438   0.1596 \\
  -1.8125   0.1727 \\
  -1.7812   0.1844 \\
  -1.7500   0.1948 \\
  -1.7188   0.2043 \\
  -1.6875   0.2130 \\
  -1.6562   0.2210 \\
  -1.6250   0.2283 \\
  -1.5938   0.2352 \\
  -1.5625   0.2416 \\
  -1.5312   0.2476 \\
  -1.5000   0.2532 \\
  -1.4688   0.2584 \\
  -1.4375   0.2634 \\
  -1.4062   0.2681 \\
  -1.3750   0.2725 \\
  -1.3438   0.2766 \\
  -1.3125   0.2806 \\
  -1.2812   0.2843 \\
  -1.2500   0.2878 \\
  -1.2188   0.2912 \\
  -1.1875   0.2944 \\
  -1.1562   0.2974 \\
  -1.1250   0.3003 \\
  -1.0938   0.3030 \\
  -1.0625   0.3056 \\
  -1.0312   0.3081 \\
  -1.0000   0.3104 \\
  -0.9688   0.3126 \\
  -0.9375   0.3147 \\
  -0.9062   0.3167 \\
  -0.8750   0.3187 \\
  -0.8438   0.3205 \\
  -0.8125   0.3222 \\
  -0.7812   0.3238 \\
  -0.7500   0.3253 \\
  -0.7188   0.3268 \\
  -0.6875   0.3282 \\
  -0.6562   0.3295 \\
  -0.6250   0.3307 \\
  -0.5938   0.3319 \\
  -0.5625   0.3329 \\
  -0.5312   0.3340 \\
  -0.5000   0.3349 \\
  -0.4688   0.3358 \\
  -0.4375   0.3366 \\
  -0.4062   0.3374 \\
  -0.3750   0.3381 \\
  -0.3438   0.3387 \\
  -0.3125   0.3393 \\
  -0.2812   0.3398 \\
  -0.2500   0.3402 \\
  -0.2188   0.3407 \\
  -0.1875   0.3410 \\
  -0.1562   0.3413 \\
  -0.1250   0.3415 \\
  -0.0938   0.3417 \\
  -0.0625   0.3419 \\
  -0.0312   0.3419 \\
        0   0.3420 \\
   0.0312   0.3419 \\
   0.0625   0.3419 \\
   0.0938   0.3417 \\
   0.1250   0.3415 \\
   0.1562   0.3413 \\
   0.1875   0.3410 \\
   0.2188   0.3407 \\
   0.2500   0.3402 \\
   0.2812   0.3398 \\
   0.3125   0.3393 \\
   0.3438   0.3387 \\
   0.3750   0.3381 \\
   0.4062   0.3374 \\
   0.4375   0.3366 \\
   0.4688   0.3358 \\
   0.5000   0.3349 \\
   0.5312   0.3340 \\
   0.5625   0.3329 \\
   0.5938   0.3319 \\
   0.6250   0.3307 \\
   0.6562   0.3295 \\
   0.6875   0.3282 \\
   0.7188   0.3268 \\
   0.7500   0.3253 \\
   0.7812   0.3238 \\
   0.8125   0.3222 \\
   0.8438   0.3205 \\
   0.8750   0.3187 \\
   0.9062   0.3167 \\
   0.9375   0.3147 \\
   0.9688   0.3126 \\
   1.0000   0.3104 \\
   1.0312   0.3081 \\
   1.0625   0.3056 \\
   1.0938   0.3030 \\
   1.1250   0.3003 \\
   1.1562   0.2974 \\
   1.1875   0.2944 \\
   1.2188   0.2912 \\
   1.2500   0.2878 \\
   1.2812   0.2843 \\
   1.3125   0.2806 \\
   1.3438   0.2766 \\
   1.3750   0.2725 \\
   1.4062   0.2681 \\
   1.4375   0.2634 \\
   1.4688   0.2584 \\
   1.5000   0.2532 \\
   1.5312   0.2476 \\
   1.5625   0.2416 \\
   1.5938   0.2352 \\
   1.6250   0.2283 \\
   1.6562   0.2210 \\
   1.6875   0.2130 \\
   1.7188   0.2043 \\
   1.7500   0.1948 \\
   1.7812   0.1844 \\
   1.8125   0.1727 \\
   1.8438   0.1596 \\
   1.8750   0.1445 \\
   1.9062   0.1267 \\
   1.9375   0.1051 \\
   1.9688   0.0700 \\
   2.0000  -0.0004 \\
};
\addlegendentry{Quadratic}

\addplot [color=red, line width=1.3pt, mark=+, mark options={solid, red}]
  table[row sep=crcr]{%
  -2.0000  -0.0004 \\
  -1.9583   0.0846 \\
  -1.9167   0.1213 \\
  -1.8750   0.1445 \\
  -1.8333   0.1654 \\
  -1.7917   0.1825 \\
  -1.7500   0.1967 \\
  -1.7083   0.2096 \\
  -1.6667   0.2210 \\
  -1.6250   0.2309 \\
  -1.5833   0.2404 \\
  -1.5417   0.2489 \\
  -1.5000   0.2563 \\
  -1.4583   0.2636 \\
  -1.4167   0.2702 \\
  -1.3750   0.2760 \\
  -1.3333   0.2818 \\
  -1.2917   0.2872 \\
  -1.2500   0.2917 \\
  -1.2083   0.2965 \\
  -1.1667   0.3008 \\
  -1.1250   0.3044 \\
  -1.0833   0.3084 \\
  -1.0417   0.3119 \\
  -1.0000   0.3147 \\
  -0.9583   0.3180 \\
  -0.9167   0.3209 \\
  -0.8750   0.3232 \\
  -0.8333   0.3260 \\
  -0.7917   0.3283 \\
  -0.7500   0.3301 \\
  -0.7083   0.3324 \\
  -0.6667   0.3342 \\
  -0.6250   0.3356 \\
  -0.5833   0.3375 \\
  -0.5417   0.3390 \\
  -0.5000   0.3399 \\
  -0.4583   0.3414 \\
  -0.4167   0.3425 \\
  -0.3750   0.3431 \\
  -0.3333   0.3444 \\
  -0.2917   0.3451 \\
  -0.2500   0.3454 \\
  -0.2083   0.3463 \\
  -0.1667   0.3468 \\
  -0.1250   0.3467 \\
  -0.0833   0.3474 \\
  -0.0417   0.3475 \\
        0   0.3472 \\
   0.0417   0.3475 \\
   0.0833   0.3474 \\
   0.1250   0.3467 \\
   0.1667   0.3468 \\
   0.2083   0.3463 \\
   0.2500   0.3454 \\
   0.2917   0.3451 \\
   0.3333   0.3444 \\
   0.3750   0.3431 \\
   0.4167   0.3425 \\
   0.4583   0.3414 \\
   0.5000   0.3399 \\
   0.5417   0.3390 \\
   0.5833   0.3375 \\
   0.6250   0.3356 \\
   0.6667   0.3342 \\
   0.7083   0.3324 \\
   0.7500   0.3301 \\
   0.7917   0.3283 \\
   0.8333   0.3260 \\
   0.8750   0.3232 \\
   0.9167   0.3209 \\
   0.9583   0.3180 \\
   1.0000   0.3147 \\
   1.0417   0.3119 \\
   1.0833   0.3084 \\
   1.1250   0.3044 \\
   1.1667   0.3008 \\
   1.2083   0.2965 \\
   1.2500   0.2917 \\
   1.2917   0.2872 \\
   1.3333   0.2818 \\
   1.3750   0.2760 \\
   1.4167   0.2702 \\
   1.4583   0.2636 \\
   1.5000   0.2563 \\
   1.5417   0.2489 \\
   1.5833   0.2404 \\
   1.6250   0.2309 \\
   1.6667   0.2210 \\
   1.7083   0.2096 \\
   1.7500   0.1967 \\
   1.7917   0.1825 \\
   1.8333   0.1654 \\
   1.8750   0.1445 \\
   1.9167   0.1213 \\
   1.9583   0.0846 \\
   2.0000  -0.0004 \\
};
\addlegendentry{Cubic}

\end{axis}

\node[align=center,font=\large, yshift=2em] (title) 
    at (current bounding box.north)
    {$\alpha$ = 8.5 ${^\circ}$ $\Lambda$ = 4\\ T.E. $\Delta\phi$ distribution};
    
\end{tikzpicture}%
}    
}
\end{tabular}  	
\caption{Analysis concerning the distribution of $\delta\phi$ on the wake surface. 
The graphs refer to a plane section normal to the x-direction. 
This section is located downstream of the wing trailing edge, at $x/c$ = 1.05.
In three of these four graphs, a mesh convergence analysis was conducted, 
while in the fourth (the one on the bottom right), the various degrees 
of approximation of the solution were compared in terms of the quality of the results.\label{fig:DeltaPhiDistribution}}
\end{figure}

Figure \ref{fig:DeltaPhiDistribution} presents the values of the potential jump $\delta \phi$ obtained at the wing trailing edge, as a function of the
span wise coordinate $y$. The automated grid generation algorithm is here set up to produce structured computational grids (as in \ref{fig:meshRef}, on the left),
which are progressively refined doubling at each additional cycle the number of nodes in both the span wise and chord wise directions. As the main focus
of the present numerical test is that of evaluating the solution at the trailing edge, such structured grids appear better suited to evaluate the spatial
convergence of the solver in such region.

The top left plot in Figure \ref{fig:DeltaPhiDistribution} refers to the convergence results obtained making use of bi-linear finite elements. In this
case, the $\delta \phi$ curves clearly
suggest that the solution converges as the number of equispaced nodes on the trailing edge is increased. This is the case also when bi-quadratic and bi-cubic finite
elements are considered, as shown in the top right and bottom left diagram, respectively. The limit solutions obtained with finite elements of different
degrees are also compared in the bottom right plot of Figure \ref{fig:DeltaPhiDistribution}. We point out that, because a structured mesh is used in the
present test, the bi-linear and bi-quadratic solutions in the bottom right plot feature the same number --- and location --- of the system
degrees of freedom. In fact, given any structured quadrilateral grid, the degrees of freedom obtained using bi-quadratic elements or using bi-linear finite
elements on a grid with an additional uniform refinement step are the same. For such a reason, the bottom right plot of Figure \ref{fig:DeltaPhiDistribution}
suggests that the best solution obtained with bi-linear and bi-quadratic finite elements, when using 128 of degrees of freedom, appear rather similar.
Also the solutions obtained with bi-cubic elements appear similar --- although with these elements it is not possible to have a comparison with exactly
the same degrees of freedom number and distribution used for the lower degree discretizations. However, the plots in Figure \ref{fig:DeltaPhiDistribution} also
suggest that the high order solutions seem, as expected, to converge faster to the limit solution, as the cubic and quadratic curves obtained on the coarser
grids are much closer to the limit solution than their linear counterparts.

\subsubsection{Pressure coefficient on wing sections}

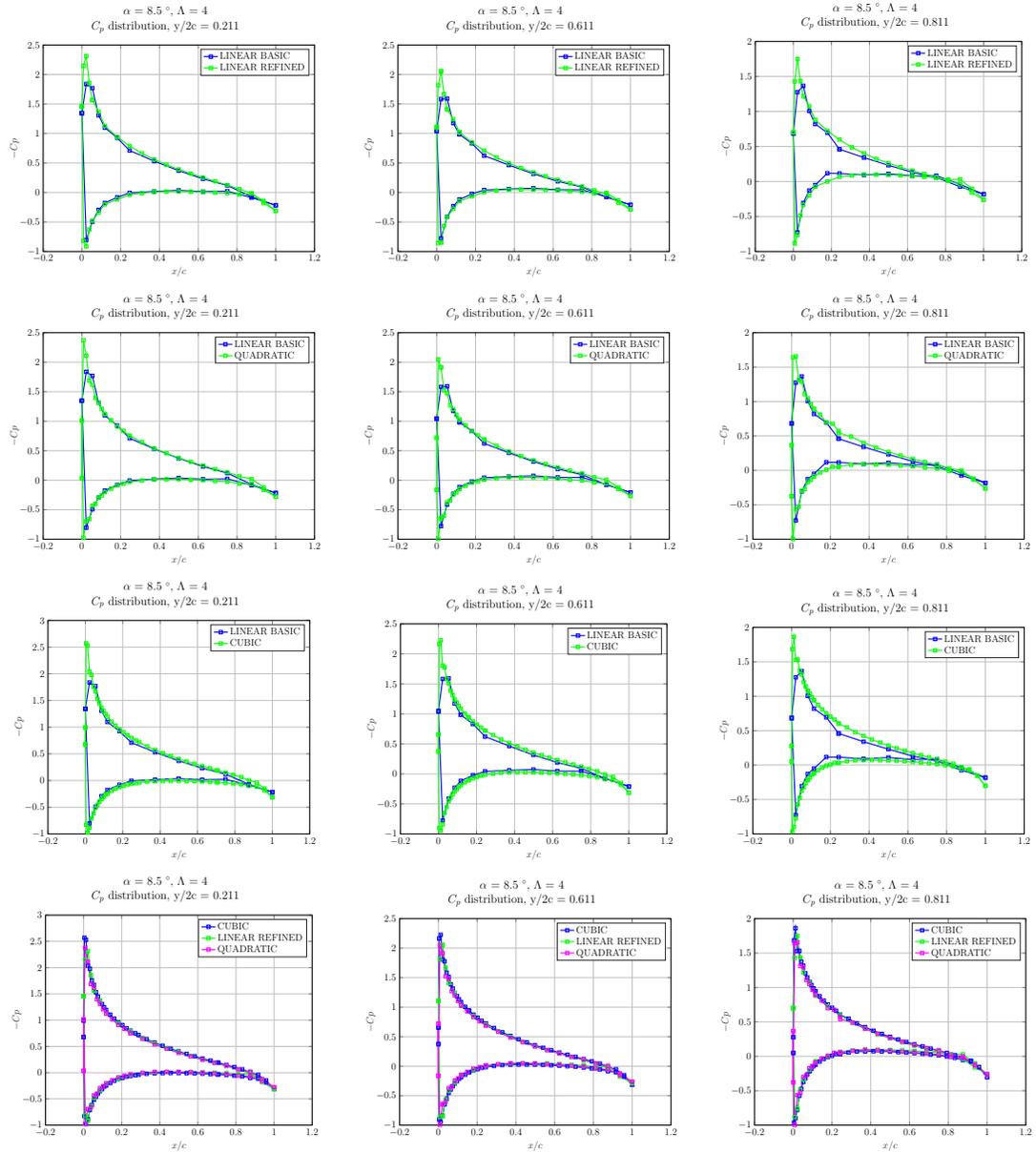
\begin{figure}
\begin{tabular}{l c r}
\resizebox{0.33\linewidth}{!}{
\centerline{
%
%
\begin{tikzpicture}

\begin{axis}[%
width=4.602in,
height=3.506in,
at={(0.772in,0.473in)},
scale only axis,
xmin=-0.2,
xmax=1.2,
xlabel style={font=\color{white!15!black}},
xlabel={\textbf{$x/c$}},
ymin=-1,
ymax=2.5,
ylabel style={font=\color{white!15!black}},
ylabel={\textbf{$-Cp$}},
axis background/.style={fill=white},
xmajorgrids,
ymajorgrids,
legend style={legend cell align=left, align=left, draw=white!15!black}
]
\addplot [color=blue, line width=1.3pt, mark=square, mark options={solid, blue}]
  table[row sep=crcr]{%
-2.10067467573863e-18	1.3463546\\
-2.10067452078647e-18	1.3463546\\
0.022500835002742	1.8375067\\
0.0540064126417008	1.7690562\\
0.0858467048895266	1.3095248\\
0.118125304877199	1.0983894\\
0.181915310111255	0.92643726\\
0.245969396727604	0.7105338\\
0.372539037906482	0.53161079\\
0.498963363571903	0.37001854\\
0.624564409614299	0.23177417\\
0.749947102206648	0.11913677\\
0.87510234979763	-0.072056696\\
1.00000001776646	-0.21733509\\
};
\addlegendentry{LINEAR BASIC}

\addplot [color=blue, line width=1.3pt, mark=square, mark options={solid, blue}, forget plot]
  table[row sep=crcr]{%
-2.10067467573863e-18	1.3463546\\
-2.10067452078647e-18	1.3463546\\
0.0225009205583868	-0.8045851\\
0.0540065905073232	-0.49306706\\
0.0858469719470585	-0.29502469\\
0.118125665562518	-0.17636052\\
0.181915867664454	-0.080585092\\
0.245970137197054	-0.0064670793\\
0.372540118370416	0.018355535\\
0.498964822744851	0.036668353\\
0.624565509484518	0.015847767\\
0.749947836688085	0.022108328\\
0.875102704426079	-0.086268716\\
1.00000001776646	-0.21733509\\
};
\addplot [color=green, line width=1.3pt, mark=square, mark options={solid, green}]
  table[row sep=crcr]{%
-2.10067467573863e-18	1.4598382\\
-1.24028267280665e-18	1.4598382\\
0.00802630116672502	2.1460655\\
0.0225008352983609	2.3131003\\
0.0379779687597454	1.8561957\\
0.0540064086856373	1.5658172\\
0.0540064126417008	1.5658172\\
0.0858467048895266	1.3703508\\
0.118125304877199	1.1283838\\
0.181915310111255	0.94024694\\
0.245969396284175	0.78853732\\
0.309257753006218	0.66285861\\
0.372539036945721	0.56038421\\
0.435777248266454	0.47060874\\
0.498963362537238	0.3916581\\
0.561791460283731	0.32008207\\
0.624564409023061	0.25389338\\
0.687282677104425	0.19028479\\
0.749947102206648	0.12739223\\
0.812554831581903	0.057264138\\
0.87510234979763	-0.015627971\\
0.93758148876842	-0.17929213\\
1.00000001776646	-0.3090359\\
};
\addlegendentry{LINEAR REFINED}

\addplot [color=green, line width=1.3pt, mark=square, mark options={solid, green}, forget plot]
  table[row sep=crcr]{%
-2.10067467573863e-18	1.4598382\\
-1.24028267280665e-18	1.4598382\\
0.00802633698660474	-0.81932634\\
0.0225009205583868	-0.91294062\\
0.0379780992901577	-0.62716496\\
0.0540065905073232	-0.4787676\\
0.0858469719470585	-0.33804163\\
0.118125665562518	-0.20082845\\
0.181915867664454	-0.11522083\\
0.245970138231719	-0.03836308\\
0.309258645384469	-0.009267698\\
0.372540149518986	0.0090931179\\
0.43577853141086	0.015944762\\
0.498964793074375	0.017096039\\
0.561792719324196	0.013949339\\
0.624565450143566	0.0081310673\\
0.687283604548716	-0.0010691882\\
0.749947836688085	-0.0094432402\\
0.812555401914088	-0.034114793\\
0.875102766723219	-0.015350613\\
0.937581715749889	-0.14156334\\
1.00000001776646	-0.3090359\\
};
\end{axis}

\node[align=center,font=\large, yshift=2em] (title) 
    at (current bounding box.north)
    {$\alpha$ = 8.5 ${^\circ}$, $\Lambda$ = 4\\ $C_p$ distribution, y/2c = 0.211};
\end{tikzpicture}}    
}
&
\resizebox{0.33\linewidth}{!}{
\centerline{
%
%
\begin{tikzpicture}

\begin{axis}[%
width=4.602in,
height=3.506in,
at={(0.772in,0.473in)},
scale only axis,
xmin=-0.2,
xmax=1.2,
xlabel style={font=\color{white!15!black}},
xlabel={\textbf{$x/c$}},
ymin=-1,
ymax=2.5,
ylabel style={font=\color{white!15!black}},
ylabel={\textbf{$-Cp$}},
axis background/.style={fill=white},
xmajorgrids,
ymajorgrids,
legend style={legend cell align=left, align=left, draw=white!15!black}
]
\addplot [color=blue, line width=1.3pt, mark=square, mark options={solid, blue}]
  table[row sep=crcr]{%
-2.02434943799647e-18	1.0447187\\
-1.24722131053235e-18	1.0447187\\
0.0222577352381907	1.5824107\\
0.0535441837511689	1.5922148\\
0.0851987611350909	1.1753397\\
0.117295679333204	0.98458451\\
0.117295679776632	0.98458451\\
0.180778165251204	0.83228636\\
0.244532380110597	0.62362581\\
0.37142817016392	0.46469167\\
0.498155411565668	0.31603265\\
0.62399847619308	0.19009674\\
0.749610722062309	0.091619998\\
0.874942343553201	-0.077700965\\
1.00000001776646	-0.20898511\\
};
\addlegendentry{LINEAR BASIC}

\addplot [color=blue, line width=1.3pt, mark=square, mark options={solid, blue}, forget plot]
  table[row sep=crcr]{%
-2.02434943799647e-18	1.0447187\\
-1.24722131053235e-18	1.0447187\\
0.0222580534023887	-0.77773809\\
0.0535448636671292	-0.41318652\\
0.0851997924719858	-0.23089908\\
0.117297068618011	-0.11669878\\
0.180780253794332	-0.023440197\\
0.244535155223405	0.041797101\\
0.371430974212348	0.059766065\\
0.498158217198355	0.072879225\\
0.624000519202426	0.046826061\\
0.749612045428836	0.046273548\\
0.874943053288412	-0.072692491\\
1.00000001776646	-0.20898511\\
};
\addplot [color=green, line width=1.3pt, mark=square, mark options={solid, green}]
  table[row sep=crcr]{%
-2.02434943799647e-18	1.1072608\\
-1.24722130905426e-18	1.1072608\\
0.00791694185315451	1.8170447\\
0.021631210322021	2.048521\\
0.0222577349425719	2.059083\\
0.037623779053015	1.6656756\\
0.0535441837511689	1.4123451\\
0.0535441874116135	1.4123449\\
0.0851987611350909	1.240484\\
0.117295679333204	1.0216686\\
0.117295679776632	1.0216686\\
0.180778165251204	0.85008788\\
0.244532377302218	0.71023911\\
0.307986327775636	0.59273106\\
0.371428170119577	0.49713197\\
0.43482098314434	0.41351634\\
0.498155410235383	0.34060201\\
0.561106099440338	0.27460366\\
0.623998475601842	0.21434057\\
0.686833105022484	0.1568047\\
0.749610720879834	0.10033423\\
0.81230864516365	0.037687205\\
0.874942343553201	-0.027540985\\
0.93750352261668	-0.1726097\\
1.00000001776646	-0.28707844\\
};
\addlegendentry{LINEAR REFINED}

\addplot [color=green, line width=1.3pt, mark=square, mark options={solid, green}, forget plot]
  table[row sep=crcr]{%
-2.02434943799647e-18	1.1072608\\
-1.24722130905426e-18	1.1072608\\
0.00791707384118923	-0.85768664\\
0.0216251717040897	-0.84868509\\
0.0222580527089917	-0.84828395\\
0.0376242751192349	-0.5654335\\
0.0535448636671292	-0.41579044\\
0.0851997924719858	-0.27690333\\
0.117297068618011	-0.14476109\\
0.180780253794332	-0.064885795\\
0.244535158475212	0.0076795914\\
0.307989098524936	0.033125121\\
0.371430974212348	0.048334397\\
0.43482378253549	0.052320223\\
0.498158217198355	0.050691228\\
0.561108570904566	0.045100786\\
0.624000578543378	0.036723182\\
0.686834816527729	0.024994286\\
0.749612045428836	0.013297787\\
0.812309686505364	-0.013772331\\
0.874943053288412	-0.005267804\\
0.937503948149521	-0.13057314\\
1.00000001776646	-0.28707844\\
};
\end{axis}

\node[align=center,font=\large, yshift=2em] (title) 
    at (current bounding box.north)
    {$\alpha$ = 8.5 ${^\circ}$, $\Lambda$ = 4\\ $C_p$ distribution, y/2c = 0.611};
\end{tikzpicture}}    
}
&
\resizebox{0.33\linewidth}{!}{
\centerline{
%
%
\begin{tikzpicture}

\begin{axis}[%
width=4.521in,
height=3.566in,
at={(0.758in,0.481in)},
scale only axis,
xmin=-0.2,
xmax=1.2,
xlabel style={font=\color{white!15!black}},
xlabel={\textbf{$x/c$}},
ymin=-1,
ymax=2,
ylabel style={font=\color{white!15!black}},
ylabel={\textbf{$-Cp$}},
axis background/.style={fill=white},
xmajorgrids,
ymajorgrids,
legend style={legend cell align=left, align=left, draw=white!15!black}
]
\addplot [color=blue, line width=1.3pt, mark=square, mark options={solid, blue}]
  table[row sep=crcr]{%
-1.8300675764247e-18	0.68493921\\
-9.83552885949522e-19	0.68493921\\
0.0216156951629845	1.275207\\
0.052387891779524	1.3651484\\
0.0836701736134814	1.0066783\\
0.115425638221391	0.82033533\\
0.178974315685522	0.69537598\\
0.242773592786953	0.46006137\\
0.242773592934762	0.46006134\\
0.371212033595351	0.34002617\\
0.499484888757955	0.231409\\
0.625262262706455	0.13044263\\
0.750825881283067	0.058112629\\
0.875586667142171	-0.074513435\\
1.00000001776646	-0.18223496\\
};
\addlegendentry{LINEAR BASIC}

\addplot [color=blue, line width=1.3pt, mark=square, mark options={solid, blue}, forget plot]
  table[row sep=crcr]{%
-1.8300675764247e-18	0.68493921\\
-9.83552885949522e-19	0.68493921\\
0.0216155712553421	-0.72666514\\
0.0523877121014522	-0.30608544\\
0.0836699947670298	-0.12858351\\
0.115425519940142	-0.050673705\\
0.178976110328969	0.11711863\\
0.242777238911956	0.11495374\\
0.371214640435728	0.09251745\\
0.499486384214262	0.10992701\\
0.625263308885059	0.080960035\\
0.750826672592786	0.083449341\\
0.875586997959648	-0.040171716\\
1.00000001776646	-0.18223496\\
};
\addplot [color=green, line width=1.3pt, mark=square, mark options={solid, green}]
  table[row sep=crcr]{%
-1.8300675764247e-18	0.69961536\\
-8.17025388750979e-19	0.69961536\\
0.00762118183428853	1.4286324\\
0.0216156801451465	1.7494626\\
0.0367163411875874	1.4345329\\
0.052387891779524	1.2196158\\
0.0836701736134814	1.0755891\\
0.115425638221391	0.8799746\\
0.178974315685522	0.72479415\\
0.242773584913838	0.59754503\\
0.306997107259851	0.48703215\\
0.371212031703391	0.40178329\\
0.435375487926569	0.32146326\\
0.499484827347671	0.26222357\\
0.562398376058854	0.20486018\\
0.625262259454648	0.15417184\\
0.68806843081763	0.10684429\\
0.75082593515507	0.061058957\\
0.813244173284663	0.010045909\\
0.875586604845031	-0.039834339\\
0.937832364541668	-0.15963557\\
1.00000001776646	-0.2603251\\
};
\addlegendentry{LINEAR REFINED}

\addplot [color=green, line width=1.3pt, mark=square, mark options={solid, green}, forget plot]
  table[row sep=crcr]{%
-1.8300675764247e-18	0.69961536\\
-8.17025388750979e-19	0.69961536\\
0.00762112051832905	-0.8832413\\
0.0216155581373077	-0.76263648\\
0.0367161815675636	-0.48516461\\
0.0523877121014522	-0.33379793\\
0.0836699947670298	-0.19796844\\
0.115425519940142	-0.072435357\\
0.178976110328969	-0.0010881792\\
0.242777236894912	0.067036182\\
0.307000230711601	0.085709527\\
0.371214643391917	0.10032669\\
0.43537756642936	0.098087557\\
0.499486384214262	0.093494236\\
0.562399694920442	0.087495267\\
0.625263313319341	0.077358484\\
0.688069300353065	0.065888859\\
0.750826679983257	0.054340661\\
0.813244725201233	0.028187478\\
0.875586997959648	0.031608168\\
0.937832562562223	-0.10279376\\
1.00000001776646	-0.2603251\\
};
\end{axis}

\node[align=center,font=\large, yshift=2em] (title) 
    at (current bounding box.north)
    {$\alpha$ = 8.5 ${^\circ}$, $\Lambda$ = 4\\ $C_p$ distribution, y/2c = 0.811};
\end{tikzpicture}}    
}
\\
\resizebox{0.33\linewidth}{!}{
\centerline{
%
%
\begin{tikzpicture}

\begin{axis}[%
width=4.602in,
height=3.506in,
at={(0.772in,0.473in)},
scale only axis,
xmin=-0.2,
xmax=1.2,
xlabel style={font=\color{white!15!black}},
xlabel={\textbf{$x/c$}},
ymin=-1,
ymax=2.5,
ylabel style={font=\color{white!15!black}},
ylabel={\textbf{$-Cp$}},
axis background/.style={fill=white},
xmajorgrids,
ymajorgrids,
legend style={legend cell align=left, align=left, draw=white!15!black}
]
\addplot [color=blue, line width=1.3pt, mark=square, mark options={solid, blue}]
  table[row sep=crcr]{%
-2.10067467573863e-18	1.3463546\\
-2.10067452078647e-18	1.3463546\\
0.022500835002742	1.8375067\\
0.0540064126417008	1.7690562\\
0.0858467048895266	1.3095248\\
0.118125304877199	1.0983894\\
0.181915310111255	0.92643726\\
0.245969396727604	0.7105338\\
0.372539037906482	0.53161079\\
0.498963363571903	0.37001854\\
0.624564409614299	0.23177417\\
0.749947102206648	0.11913677\\
0.87510234979763	-0.072056696\\
1.00000001776646	-0.21733509\\
};
\addlegendentry{LINEAR BASIC}

\addplot [color=blue, line width=1.3pt, mark=square, mark options={solid, blue}, forget plot]
  table[row sep=crcr]{%
-2.10067467573863e-18	1.3463546\\
-2.10067452078647e-18	1.3463546\\
0.0225009205583868	-0.8045851\\
0.0540065905073232	-0.49306706\\
0.0858469719470585	-0.29502469\\
0.118125665562518	-0.17636052\\
0.181915867664454	-0.080585092\\
0.245970137197054	-0.0064670793\\
0.372540118370416	0.018355535\\
0.498964822744851	0.036668353\\
0.624565509484518	0.015847767\\
0.749947836688085	0.022108328\\
0.875102704426079	-0.086268716\\
1.00000001776646	-0.21733509\\
};
\addplot [color=green, line width=1.3pt, mark=square, mark options={solid, green}]
  table[row sep=crcr]{%
-2.10067467573863e-18	0.036717705\\
-1.24028267280665e-18	1.0140986\\
0.00802647451397611	2.3725977\\
0.0225008352983609	2.1128097\\
0.0379779684641266	1.687606\\
0.0540064086856373	1.6154777\\
0.0698462598637024	1.3953323\\
0.0858467048895266	1.3121706\\
0.101949780581038	1.2095926\\
0.118125304285961	1.1254386\\
0.118125304877199	1.1253982\\
0.149964858434306	1.0159061\\
0.181915310111255	0.91010666\\
0.213925407485033	0.83552575\\
0.245969395988557	0.7534489\\
0.245969396284175	0.7533409\\
0.309257753006218	0.64402771\\
0.372539036945721	0.53834796\\
0.435777248266454	0.45549548\\
0.498963362537238	0.37747881\\
0.561791460283731	0.30999428\\
0.624564409023061	0.24771619\\
0.687282677104425	0.18639621\\
0.749947102206648	0.13422219\\
0.812554831581903	0.063259907\\
0.87510234979763	0.019904509\\
0.93758148876842	-0.12129729\\
1.00000001776646	-0.2801615\\
};
\addlegendentry{QUADRATIC}

\addplot [color=green, line width=1.3pt, mark=square, mark options={solid, green}, forget plot]
  table[row sep=crcr]{%
-2.10067467573863e-18	0.036717705\\
-1.24028267280665e-18	1.0140986\\
0.00802651036116693	-0.97565609\\
0.0225009205583868	-0.69419378\\
0.0379780992901577	-0.65719736\\
0.0540065905073232	-0.43300971\\
0.0698464883133735	-0.40116966\\
0.0858469719470585	-0.29656789\\
0.101950098464154	-0.26007688\\
0.118125666744993	-0.20275322\\
0.149965312872185	-0.14389247\\
0.181915867664454	-0.087783165\\
0.213926045916388	-0.059768148\\
0.245970138231719	-0.030141246\\
0.309258645384469	-0.0050371056\\
0.372540149518986	0.015182843\\
0.43577853141086	0.014717768\\
0.498964793074375	0.016934801\\
0.561792719324196	0.0073482529\\
0.624565450143566	0.0013127846\\
0.687283604548716	-0.015242588\\
0.749947836688085	-0.024871675\\
0.812555401914088	-0.057220351\\
0.875102766723219	-0.066563718\\
0.937581715749889	-0.16339183\\
1.00000001776646	-0.2801615\\
};
\end{axis}

\node[align=center,font=\large, yshift=2em] (title) 
    at (current bounding box.north)
    {$\alpha$ = 8.5 ${^\circ}$, $\Lambda$ = 4\\ $C_p$ distribution, y/2c = 0.211};
\end{tikzpicture}}    
}
&
\resizebox{0.33\linewidth}{!}{
\centerline{
%
%
\begin{tikzpicture}

\begin{axis}[%
width=4.602in,
height=3.506in,
at={(0.772in,0.473in)},
scale only axis,
xmin=-0.2,
xmax=1.2,
xlabel style={font=\color{white!15!black}},
xlabel={\textbf{$x/c$}},
ymin=-1,
ymax=2.5,
ylabel style={font=\color{white!15!black}},
ylabel={\textbf{$-Cp$}},
axis background/.style={fill=white},
xmajorgrids,
ymajorgrids,
legend style={legend cell align=left, align=left, draw=white!15!black}
]
\addplot [color=blue, line width=1.3pt, mark=square, mark options={solid, blue}]
  table[row sep=crcr]{%
-2.02434943799647e-18	1.0447187\\
-1.24722131053235e-18	1.0447187\\
0.0222577352381907	1.5824107\\
0.0535441837511689	1.5922148\\
0.0851987611350909	1.1753397\\
0.117295679333204	0.98458451\\
0.117295679776632	0.98458451\\
0.180778165251204	0.83228636\\
0.244532380110597	0.62362581\\
0.37142817016392	0.46469167\\
0.498155411565668	0.31603265\\
0.62399847619308	0.19009674\\
0.749610722062309	0.091619998\\
0.874942343553201	-0.077700965\\
1.00000001776646	-0.20898511\\
};
\addlegendentry{LINEAR BASIC}

\addplot [color=blue, line width=1.3pt, mark=square, mark options={solid, blue}, forget plot]
  table[row sep=crcr]{%
-2.02434943799647e-18	1.0447187\\
-1.24722131053235e-18	1.0447187\\
0.0222580534023887	-0.77773809\\
0.0535448636671292	-0.41318652\\
0.0851997924719858	-0.23089908\\
0.117297068618011	-0.11669878\\
0.180780253794332	-0.023440197\\
0.244535155223405	0.041797101\\
0.371430974212348	0.059766065\\
0.498158217198355	0.072879225\\
0.624000519202426	0.046826061\\
0.749612045428836	0.046273548\\
0.874943053288412	-0.072692491\\
1.00000001776646	-0.20898511\\
};
\addplot [color=green, line width=1.3pt, mark=square, mark options={solid, green}]
  table[row sep=crcr]{%
-2.02434943799647e-18	0.71826309\\
-1.24722130905426e-18	-0.16234292\\
0.00792250348087713	2.0466621\\
0.0216616442372479	1.918523\\
0.0222577349425719	1.9129419\\
0.037623857343937	1.5233558\\
0.0535441874116135	1.4751093\\
0.0692905640088142	1.2718035\\
0.0851987608394721	1.2005014\\
0.101210587317448	1.1060445\\
0.117295678150729	1.0304499\\
0.117295678741966	1.0303353\\
0.148980421080693	0.92922091\\
0.180778164216538	0.83336627\\
0.212637465244411	0.7616663\\
0.244532377302218	0.69070745\\
0.244532405642409	0.68998605\\
0.307986327775636	0.58330077\\
0.371428170119577	0.48484957\\
0.43482098314434	0.40726018\\
0.498155410235383	0.33474118\\
0.561106099440338	0.27274638\\
0.623998475601842	0.2155547\\
0.686833105022484	0.15979376\\
0.749610720879834	0.11234608\\
0.81230864516365	0.047948029\\
0.874942343553201	0.0075814654\\
0.93750352261668	-0.121348\\
1.00000001776646	-0.26613307\\
};
\addlegendentry{QUADRATIC}

\addplot [color=green, line width=1.3pt, mark=square, mark options={solid, green}, forget plot]
  table[row sep=crcr]{%
-2.02434943799647e-18	0.71826309\\
-1.24722130905426e-18	-0.16234292\\
0.00792263853941618	-0.99438637\\
0.0216559485803761	-0.65170497\\
0.0222580527089917	-0.63669246\\
0.037624358322658	-0.60188627\\
0.0535448642583668	-0.3741847\\
0.0692914216916997	-0.34536976\\
0.0851997930632234	-0.24207138\\
0.101211799903776	-0.20738861\\
0.117297058727852	-0.15188091\\
0.117297069652676	-0.1519822\\
0.148982167936896	-0.095866509\\
0.180780254828998	-0.042281218\\
0.212639906422036	-0.016766028\\
0.244535158475212	0.012036184\\
0.244535159657687	0.01157608\\
0.307989098524936	0.033080928\\
0.371430974212348	0.049800787\\
0.43482378253549	0.046861839\\
0.498158217198355	0.046314884\\
0.561108570904566	0.034449816\\
0.624000578543378	0.026054965\\
0.686834816527729	0.0072571211\\
0.749612045428836	-0.0049986853\\
0.812309686505364	-0.039559159\\
0.874943053288412	-0.053601913\\
0.937503948149521	-0.15064526\\
1.00000001776646	-0.26613307\\
};
\end{axis}

\node[align=center,font=\large, yshift=2em] (title) 
    at (current bounding box.north)
    {$\alpha$ = 8.5 ${^\circ}$, $\Lambda$ = 4\\ $C_p$ distribution, y/2c = 0.611};
\end{tikzpicture}}    
}
&
\resizebox{0.33\linewidth}{!}{
\centerline{
%
%
\begin{tikzpicture}

\begin{axis}[%
width=4.602in,
height=3.506in,
at={(0.772in,0.473in)},
scale only axis,
xmin=-0.2,
xmax=1.2,
xlabel style={font=\color{white!15!black}},
xlabel={\textbf{$x/c$}},
ymin=-1,
ymax=2,
ylabel style={font=\color{white!15!black}},
ylabel={\textbf{$-Cp$}},
axis background/.style={fill=white},
xmajorgrids,
ymajorgrids,
legend style={legend cell align=left, align=left, draw=white!15!black}
]
\addplot [color=blue, line width=1.3pt, mark=square, mark options={solid, blue}]
  table[row sep=crcr]{%
-1.8300675764247e-18	0.68493921\\
-9.83552885949522e-19	0.68493921\\
0.0216156951629845	1.275207\\
0.052387891779524	1.3651484\\
0.0836701736134814	1.0066783\\
0.115425638221391	0.82033533\\
0.178974315685522	0.69537598\\
0.242773592786953	0.46006137\\
0.242773592934762	0.46006134\\
0.371212033595351	0.34002617\\
0.499484888757955	0.231409\\
0.625262262706455	0.13044263\\
0.750825881283067	0.058112629\\
0.875586667142171	-0.074513435\\
1.00000001776646	-0.18223496\\
};
\addlegendentry{LINEAR BASIC}

\addplot [color=blue, line width=1.3pt, mark=square, mark options={solid, blue}, forget plot]
  table[row sep=crcr]{%
-1.8300675764247e-18	0.68493921\\
-9.83552885949522e-19	0.68493921\\
0.0216155712553421	-0.72666514\\
0.0523877121014522	-0.30608544\\
0.0836699947670298	-0.12858351\\
0.115425519940142	-0.050673705\\
0.178976110328969	0.11711863\\
0.242777238911956	0.11495374\\
0.371214640435728	0.09251745\\
0.499486384214262	0.10992701\\
0.625263308885059	0.080960035\\
0.750826672592786	0.083449341\\
0.875586997959648	-0.040171716\\
1.00000001776646	-0.18223496\\
};
\addplot [color=green, line width=1.3pt, mark=square, mark options={solid, green}]
  table[row sep=crcr]{%
-1.8300675764247e-18	0.36683193\\
-8.17025388750979e-19	-0.37652409\\
0.00762587522916001	1.6427321\\
0.0216156801451465	1.6563702\\
0.0367164841166687	1.3141203\\
0.0523878768747048	1.2903523\\
0.0679438165567362	1.1080149\\
0.0836701550590769	1.0489005\\
0.099507816943524	0.96475548\\
0.115425594617615	0.8936277\\
0.115425609546365	0.89570719\\
0.147145272838322	0.80942917\\
0.178974300461153	0.69992328\\
0.210859384162065	0.67793906\\
0.242773584913838	0.57264483\\
0.24277358652667	0.54137844\\
0.30699710740766	0.49127528\\
0.371212031703391	0.40065923\\
0.435375487926569	0.33123535\\
0.499484827347671	0.26715472\\
0.562398376058854	0.21390681\\
0.625262259454648	0.16588049\\
0.68806843081763	0.11983029\\
0.75082593515507	0.082076438\\
0.813244173284663	0.02823429\\
0.875586604845031	8.0749451e-05\\
0.937832366019762	-0.12003584\\
1.00000001776646	-0.26250795\\
};
\addlegendentry{QUADRATIC}

\addplot [color=green, line width=1.3pt, mark=square, mark options={solid, green}, forget plot]
  table[row sep=crcr]{%
-1.8300675764247e-18	0.36683193\\
-8.17025388750979e-19	-0.37652409\\
0.00762581470180202	-0.99838495\\
0.0216155581373077	-0.56289697\\
0.036716320634094	-0.52793169\\
0.0523876969945267	-0.29793435\\
0.0679436303744736	-0.27222234\\
0.0836699898350224	-0.17167231\\
0.0995076731489009	-0.1382522\\
0.115425482648936	-0.089084223\\
0.115425488091741	-0.088131212\\
0.14714611684868	-0.030985828\\
0.178976114776323	0.0072704963\\
0.210862122729364	0.051156536\\
0.242777224800684	0.047870021\\
0.242777236894912	0.058793563\\
0.307000230711601	0.080374181\\
0.371214643391917	0.094828166\\
0.43537756642936	0.086647063\\
0.499486384214262	0.08379633\\
0.562399694920442	0.069619186\\
0.625263313319341	0.060141865\\
0.688069300353065	0.040529653\\
0.750826679983257	0.028392183\\
0.813244725201233	-0.0082420707\\
0.875586997959648	-0.020341894\\
0.937832562562223	-0.13161916\\
1.00000001776646	-0.26250795\\
};
\end{axis}

\node[align=center,font=\large, yshift=2em] (title) 
    at (current bounding box.north)
    {$\alpha$ = 8.5 ${^\circ}$, $\Lambda$ = 4\\ $C_p$ distribution, y/2c = 0.811};
\end{tikzpicture}}    
}
\\
\resizebox{0.33\linewidth}{!}{
\centerline{
%
%
\begin{tikzpicture}

\begin{axis}[%
width=4.44in,
height=3.625in,
at={(0.745in,0.489in)},
scale only axis,
xmin=-0.2,
xmax=1.2,
xlabel style={font=\color{white!15!black}},
xlabel={\textbf{$x/c$}},
ymin=-1,
ymax=3,
ylabel style={font=\color{white!15!black}},
ylabel={\textbf{$-Cp$}},
axis background/.style={fill=white},
xmajorgrids,
ymajorgrids,
legend style={legend cell align=left, align=left, draw=white!15!black}
]
\addplot [color=blue, line width=1.3pt, mark=square, mark options={solid, blue}]
  table[row sep=crcr]{%
-2.10067467573863e-18	1.3463546\\
-2.10067452078647e-18	1.3463546\\
0.022500835002742	1.8375067\\
0.0540064126417008	1.7690562\\
0.0858467048895266	1.3095248\\
0.118125304877199	1.0983894\\
0.181915310111255	0.92643726\\
0.245969396727604	0.7105338\\
0.372539037906482	0.53161079\\
0.498963363571903	0.37001854\\
0.624564409614299	0.23177417\\
0.749947102206648	0.11913677\\
0.87510234979763	-0.072056696\\
1.00000001776646	-0.21733509\\
};
\addlegendentry{LINEAR BASIC}

\addplot [color=blue, line width=1.3pt, mark=square, mark options={solid, blue}, forget plot]
  table[row sep=crcr]{%
-2.10067467573863e-18	1.3463546\\
-2.10067452078647e-18	1.3463546\\
0.0225009205583868	-0.8045851\\
0.0540065905073232	-0.49306706\\
0.0858469719470585	-0.29502469\\
0.118125665562518	-0.17636052\\
0.181915867664454	-0.080585092\\
0.245970137197054	-0.0064670793\\
0.372540118370416	0.018355535\\
0.498964822744851	0.036668353\\
0.624565509484518	0.015847767\\
0.749947836688085	0.022108328\\
0.875102704426079	-0.086268716\\
1.00000001776646	-0.21733509\\
};
\addplot [color=green, line width=1.3pt, mark=square, mark options={solid, green}]
  table[row sep=crcr]{%
-2.10067452078647e-18	0.99438006\\
-1.88228638582612e-18	0.67665339\\
0.00415745706256359	2.5688586\\
0.0124023321454434	2.5293353\\
0.0225008352983609	2.037878\\
0.0327421115943696	1.9779305\\
0.0432713961686896	1.7577556\\
0.0540064126417008	1.6708629\\
0.0645461089431045	1.5305665\\
0.075163926258632	1.4484965\\
0.0858467048895266	1.3652699\\
0.0965733645374024	1.3057446\\
0.107334193841775	1.2459959\\
0.118125304285961	1.1936649\\
0.139337921226713	1.1051929\\
0.160603936010295	1.0294514\\
0.181915309963445	0.96028548\\
0.203251304786666	0.90589374\\
0.224603284959509	0.84694284\\
0.245969425659032	0.80128413\\
0.270951174523791	0.75217587\\
0.288162576786297	0.7183392\\
0.33035221166873	0.64004529\\
0.372539036664883	0.57480943\\
0.414704268235396	0.51275921\\
0.456844581118166	0.45722273\\
0.498963361946	0.4055348\\
0.540855414081315	0.35799524\\
0.582721471876939	0.3129541\\
0.624564409023061	0.27009472\\
0.666382801629308	0.22900268\\
0.702323885514941	0.19402736\\
0.708176672470286	0.18833165\\
0.7499471615476	0.14863952\\
0.79169229572529	0.10609136\\
0.833410862525808	0.062784694\\
0.875102409138581	0.0050342479\\
0.916762856420234	-0.043032736\\
0.958389079193298	-0.14945902\\
1.00000001776646	-0.31565344\\
};
\addlegendentry{CUBIC}

\addplot [color=green, line width=1.3pt, mark=square, mark options={solid, green}, forget plot]
  table[row sep=crcr]{%
-2.10067452078647e-18	0.99438006\\
-1.88228638582612e-18	0.67665339\\
0.00415747746351858	-0.8277986\\
0.0124023837050943	-0.97464365\\
0.0225009188759739	-0.88725352\\
0.0327422267254169	-0.70996028\\
0.0432715430435515	-0.61484593\\
0.0540065905073232	-0.51366335\\
0.0645463248171926	-0.45755571\\
0.0751641707964003	-0.3982659\\
0.0858469788701695	-0.34963945\\
0.0965736701831975	-0.31233677\\
0.107334525328225	-0.26894519\\
0.118125666744993	-0.24395533\\
0.139338344171045	-0.19717513\\
0.160604416633501	-0.15355468\\
0.181915867664454	-0.12352144\\
0.20325189717262	-0.10003515\\
0.22460394134366	-0.074197888\\
0.245970138231719	-0.06115276\\
0.270840202328811	-0.048446424\\
0.288163415534036	-0.039594177\\
0.330353196467016	-0.020607136\\
0.372540149518986	-0.013994689\\
0.414705453636153	-0.0069990898\\
0.456845877817938	-0.0054265885\\
0.498964822744851	-0.005914703\\
0.540856724574902	-0.0087178107\\
0.58272264026283	-0.012847953\\
0.624565450143566	-0.018607693\\
0.666383815275177	-0.025334435\\
0.702424920789416	-0.032731742\\
0.708177512714931	-0.033912245\\
0.749947839644273	-0.043142859\\
0.791692883136026	-0.056659851\\
0.833411352948819	-0.070989475\\
0.875102766723219	-0.10037051\\
0.91676309846275	-0.11439648\\
0.958389181247267	-0.18403201\\
1.00000001776646	-0.31565344\\
};
\end{axis}

\node[align=center,font=\large, yshift=2em] (title) 
    at (current bounding box.north)
    {$\alpha$ = 8.5 ${^\circ}$, $\Lambda$ = 4\\ $C_p$ distribution, y/2c = 0.211};
\end{tikzpicture}}    
}
&
\resizebox{0.33\linewidth}{!}{
\centerline{
%
%
\begin{tikzpicture}

\begin{axis}[%
width=4.521in,
height=3.566in,
at={(0.758in,0.481in)},
scale only axis,
xmin=-0.2,
xmax=1.2,
xlabel style={font=\color{white!15!black}},
xlabel={\textbf{$x/c$}},
ymin=-1,
ymax=2.5,
ylabel style={font=\color{white!15!black}},
ylabel={\textbf{$-Cp$}},
axis background/.style={fill=white},
xmajorgrids,
ymajorgrids,
legend style={legend cell align=left, align=left, draw=white!15!black}
]
\addplot [color=blue, line width=1.3pt, mark=square, mark options={solid, blue}]
  table[row sep=crcr]{%
-2.02434943799647e-18	1.0447187\\
-1.24722131053235e-18	1.0447187\\
0.0222577352381907	1.5824107\\
0.0535441837511689	1.5922148\\
0.0851987611350909	1.1753397\\
0.117295679333204	0.98458451\\
0.117295679776632	0.98458451\\
0.180778165251204	0.83228636\\
0.244532380110597	0.62362581\\
0.37142817016392	0.46469167\\
0.498155411565668	0.31603265\\
0.62399847619308	0.19009674\\
0.749610722062309	0.091619998\\
0.874942343553201	-0.077700965\\
1.00000001776646	-0.20898511\\
};
\addlegendentry{LINEAR BASIC}

\addplot [color=blue, line width=1.3pt, mark=square, mark options={solid, blue}, forget plot]
  table[row sep=crcr]{%
-2.02434943799647e-18	1.0447187\\
-1.24722131053235e-18	1.0447187\\
0.0222580534023887	-0.77773809\\
0.0535448636671292	-0.41318652\\
0.0851997924719858	-0.23089908\\
0.117297068618011	-0.11669878\\
0.180780253794332	-0.023440197\\
0.244535155223405	0.041797101\\
0.371430974212348	0.059766065\\
0.498158217198355	0.072879225\\
0.624000519202426	0.046826061\\
0.749612045428836	0.046273548\\
0.874943053288412	-0.072692491\\
1.00000001776646	-0.20898511\\
};
\addplot [color=green, line width=1.3pt, mark=square, mark options={solid, green}]
  table[row sep=crcr]{%
-1.14213807692232e-18	0.65621138\\
-9.00182759915995e-19	0.37439698\\
0.00409915120658076	2.1641905\\
0.0122539238823407	2.2258725\\
0.0222577349425719	1.8049899\\
0.0324245843031084	1.7769955\\
0.0428809726444962	1.5844892\\
0.0535441871159947	1.5125289\\
0.0640214169077506	1.3869342\\
0.0745773459333662	1.3147616\\
0.0851987608394721	1.2402364\\
0.0958644849257839	1.1871358\\
0.10656477265969	1.1334383\\
0.117295678150729	1.0858679\\
0.138404899644579	1.0056758\\
0.156312530331545	0.94695622\\
0.159568392679639	0.93627512\\
0.18077816392092	0.87361723\\
0.202013450893539	0.82220846\\
0.223265470155025	0.76919907\\
0.244532375676315	0.72515655\\
0.244532376267552	0.7236799\\
0.286836489337739	0.64750582\\
0.329134900924747	0.57508403\\
0.371428170119577	0.51349974\\
0.413697122327436	0.45595461\\
0.45593837235865	0.4042204\\
0.498155409939764	0.35633692\\
0.54012945602504	0.31232253\\
0.582076224298323	0.27090472\\
0.623998475601842	0.2315651\\
0.665894805773255	0.19406715\\
0.707765086521932	0.15704507\\
0.749610722062309	0.12105686\\
0.791416427729308	0.082286686\\
0.833193887244142	0.043075569\\
0.874942343553201	-0.010725223\\
0.916658007270825	-0.052915532\\
0.958337596455102	-0.15310237\\
1.00000001776646	-0.31396174\\
};
\addlegendentry{CUBIC}

\addplot [color=green, line width=1.3pt, mark=square, mark options={solid, green}, forget plot]
  table[row sep=crcr]{%
-1.14213807692232e-18	0.65621138\\
-9.00182759915995e-19	0.37439698\\
0.00409922801405183	-0.90321326\\
0.0122541180827971	-0.94546849\\
0.0222580527089917	-0.84152615\\
0.0324250217720263	-0.65002596\\
0.0428815270818934	-0.55367446\\
0.0535448603023033	-0.4516069\\
0.0640222146110596	-0.39665622\\
0.0745782566575152	-0.33818033\\
0.0851997861401124	-0.29093486\\
0.0958656360681177	-0.2544843\\
0.106566040607035	-0.2124259\\
0.117297069652676	-0.18827038\\
0.138406502367362	-0.14385739\\
0.156275732160386	-0.10838422\\
0.159570254559709	-0.10182549\\
0.180780255863664	-0.074527323\\
0.202015770367388	-0.051268533\\
0.223268009443464	-0.028543672\\
0.244535159657687	-0.016292769\\
0.286839243558578	0.004215423\\
0.329137659663322	0.019813737\\
0.371430944541873	0.025097584\\
0.413699902488846	0.029555002\\
0.45594120534161	0.029271344\\
0.498158217198355	0.026891135\\
0.540132045456626	0.022306906\\
0.582078533250478	0.016397014\\
0.624000578543378	0.0088933362\\
0.66589667756479	0.00033328607\\
0.707766754594174	-0.010108377\\
0.749612045428836	-0.021341622\\
0.791417608498584	-0.036946569\\
0.833194756254551	-0.053708084\\
0.874943053288412	-0.085467763\\
0.916658454279872	-0.10386589\\
0.958337804698634	-0.17828794\\
1.00000001776646	-0.31396174\\
};
\end{axis}

\node[align=center,font=\large, yshift=2em] (title) 
    at (current bounding box.north)
    {$\alpha$ = 8.5 ${^\circ}$, $\Lambda$ = 4\\ $C_p$ distribution, y/2c = 0.611};
\end{tikzpicture}}    
}
&
\resizebox{0.33\linewidth}{!}{
\centerline{
%
%
\begin{tikzpicture}

\begin{axis}[%
width=4.602in,
height=3.506in,
at={(0.772in,0.473in)},
scale only axis,
xmin=-0.2,
xmax=1.2,
xlabel style={font=\color{white!15!black}},
xlabel={\textbf{$x/c$}},
ymin=-1,
ymax=2,
ylabel style={font=\color{white!15!black}},
ylabel={\textbf{$-Cp$}},
axis background/.style={fill=white},
xmajorgrids,
ymajorgrids,
legend style={legend cell align=left, align=left, draw=white!15!black}
]
\addplot [color=blue, line width=1.3pt, mark=square, mark options={solid, blue}]
  table[row sep=crcr]{%
-1.8300675764247e-18	0.68493921\\
-9.83552885949522e-19	0.68493921\\
0.0216156951629845	1.275207\\
0.052387891779524	1.3651484\\
0.0836701736134814	1.0066783\\
0.115425638221391	0.82033533\\
0.178974315685522	0.69537598\\
0.242773592786953	0.46006137\\
0.242773592934762	0.46006134\\
0.371212033595351	0.34002617\\
0.499484888757955	0.231409\\
0.625262262706455	0.13044263\\
0.750825881283067	0.058112629\\
0.875586667142171	-0.074513435\\
1.00000001776646	-0.18223496\\
};
\addlegendentry{LINEAR BASIC}

\addplot [color=blue, line width=1.3pt, mark=square, mark options={solid, blue}, forget plot]
  table[row sep=crcr]{%
-1.8300675764247e-18	0.68493921\\
-9.83552885949522e-19	0.68493921\\
0.0216155712553421	-0.72666514\\
0.0523877121014522	-0.30608544\\
0.0836699947670298	-0.12858351\\
0.115425519940142	-0.050673705\\
0.178976110328969	0.11711863\\
0.242777238911956	0.11495374\\
0.371214640435728	0.09251745\\
0.499486384214262	0.10992701\\
0.625263308885059	0.080960035\\
0.750826672592786	0.083449341\\
0.875586997959648	-0.040171716\\
1.00000001776646	-0.18223496\\
};
\addplot [color=green, line width=1.3pt, mark=square, mark options={solid, green}]
  table[row sep=crcr]{%
-1.54963118995704e-18	0.27725598\\
-8.18760392295643e-19	0.050570723\\
0.00392701810251808	1.683823\\
0.0118442665829781	1.8616087\\
0.02161567683683	1.5279721\\
0.0316037413049256	1.5323106\\
0.041889270345907	1.3722105\\
0.0523878704058813	1.3155442\\
0.062737225328494	1.2070961\\
0.0731691518516966	1.1455313\\
0.0836701513638417	1.0806531\\
0.0849760738031609	1.0749229\\
0.0942190294122309	1.0343511\\
0.104805520422297	0.98717183\\
0.115425593582949	0.94503903\\
0.115425600234373	0.94368833\\
0.136558476761236	0.87247592\\
0.157744220509282	0.80973178\\
0.178974295583442	0.75373411\\
0.20022751538331	0.70550197\\
0.221494575843468	0.66106158\\
0.242773583866101	0.60630143\\
0.242773584309529	0.60238242\\
0.285590235855046	0.54779869\\
0.328402874894375	0.47980261\\
0.371212031112153	0.42504907\\
0.413994302468823	0.37223929\\
0.456750730146968	0.32627368\\
0.499484826756434	0.28434667\\
0.541433255942272	0.24528919\\
0.583358064905119	0.20983364\\
0.625262317613125	0.17608741\\
0.667138723424506	0.14419271\\
0.708992818345403	0.11313544\\
0.750825875814118	0.08308211\\
0.792446480935992	0.050622411\\
0.834033637606829	0.018023754\\
0.875586606323125	-0.028062705\\
0.917093432404952	-0.063333198\\
0.958558266919661	-0.15377553\\
1.00000001776646	-0.30251953\\
};
\addlegendentry{CUBIC}

\addplot [color=green, line width=1.3pt, mark=square, mark options={solid, green}, forget plot]
  table[row sep=crcr]{%
-1.54963118995704e-18	0.27725598\\
-8.18760392295643e-19	0.050570723\\
0.00392698310241853	-0.96461976\\
0.0118441849778742	-0.89918256\\
0.0216155522488625	-0.77713251\\
0.0316035935890773	-0.57254791\\
0.0418891030807539	-0.47496161\\
0.0523876918907786	-0.37278166\\
0.0627370415134146	-0.31906211\\
0.0731689711406572	-0.26170692\\
0.0836699792406074	-0.21687232\\
0.0849714104316008	-0.21237987\\
0.0942188726148191	-0.18067032\\
0.104805373172876	-0.14172307\\
0.115425473941253	-0.11810376\\
0.115425488104812	-0.11850182\\
0.136559022474748	-0.075128399\\
0.157745396319784	-0.036952913\\
0.178976100612764	-0.0095137721\\
0.200229943317113	0.0082446421\\
0.221497620300509	0.033464823\\
0.242777227017825	0.034172844\\
0.242777228052491	0.035677571\\
0.28559351571825	0.055495754\\
0.328405816007735	0.070014924\\
0.371214613721441	0.070820905\\
0.413996528466602	0.074265525\\
0.456752599076327	0.071732879\\
0.499486384214262	0.068782106\\
0.541434606172042	0.061834406\\
0.583359293377541	0.055074498\\
0.625263311841247	0.046879616\\
0.667139618151478	0.037408415\\
0.708993684510187	0.026260639\\
0.750826620642305	0.014236098\\
0.792447165095704	-0.0026290566\\
0.834034107293575	-0.020926602\\
0.875586997959648	-0.055488463\\
0.917093708101339	-0.076441534\\
0.958558402145052	-0.15676838\\
1.00000001776646	-0.30251953\\
};
\end{axis}

\node[align=center,font=\large, yshift=2em] (title) 
    at (current bounding box.north)
    {$\alpha$ = 8.5 ${^\circ}$, $\Lambda$ = 4\\ $C_p$ distribution, y/2c = 0.811};
\end{tikzpicture}}    
}
\\
\resizebox{0.33\linewidth}{!}{
\centerline{
%
%
\definecolor{mycolor1}{rgb}{1.00000,0.00000,1.00000}%
\begin{tikzpicture}

\begin{axis}[%
width=4.521in,
height=3.566in,
at={(0.758in,0.481in)},
scale only axis,
xmin=-0.2,
xmax=1.2,
xlabel style={font=\color{white!15!black}},
xlabel={\textbf{$x/c$}},
ymin=-1,
ymax=3,
ylabel style={font=\color{white!15!black}},
ylabel={\textbf{$-Cp$}},
axis background/.style={fill=white},
xmajorgrids,
ymajorgrids,
legend style={legend cell align=left, align=left, draw=white!15!black}
]
\addplot [color=blue, line width=1.3pt, mark=square, mark options={solid, blue}]
  table[row sep=crcr]{%
-2.10067452078647e-18	0.99438006\\
-1.88228638582612e-18	0.67665339\\
0.00415745706256359	2.5688586\\
0.0124023321454434	2.5293353\\
0.0225008352983609	2.037878\\
0.0327421115943696	1.9779305\\
0.0432713961686896	1.7577556\\
0.0540064126417008	1.6708629\\
0.0645461089431045	1.5305665\\
0.075163926258632	1.4484965\\
0.0858467048895266	1.3652699\\
0.0965733645374024	1.3057446\\
0.107334193841775	1.2459959\\
0.118125304285961	1.1936649\\
0.139337921226713	1.1051929\\
0.160603936010295	1.0294514\\
0.181915309963445	0.96028548\\
0.203251304786666	0.90589374\\
0.224603284959509	0.84694284\\
0.245969425659032	0.80128413\\
0.270951174523791	0.75217587\\
0.288162576786297	0.7183392\\
0.33035221166873	0.64004529\\
0.372539036664883	0.57480943\\
0.414704268235396	0.51275921\\
0.456844581118166	0.45722273\\
0.498963361946	0.4055348\\
0.540855414081315	0.35799524\\
0.582721471876939	0.3129541\\
0.624564409023061	0.27009472\\
0.666382801629308	0.22900268\\
0.702323885514941	0.19402736\\
0.708176672470286	0.18833165\\
0.7499471615476	0.14863952\\
0.79169229572529	0.10609136\\
0.833410862525808	0.062784694\\
0.875102409138581	0.0050342479\\
0.916762856420234	-0.043032736\\
0.958389079193298	-0.14945902\\
1.00000001776646	-0.31565344\\
};
\addlegendentry{CUBIC}

\addplot [color=blue, line width=1.3pt, mark=square, mark options={solid, blue}, forget plot]
  table[row sep=crcr]{%
-2.10067452078647e-18	0.99438006\\
-1.88228638582612e-18	0.67665339\\
0.00415747746351858	-0.8277986\\
0.0124023837050943	-0.97464365\\
0.0225009188759739	-0.88725352\\
0.0327422267254169	-0.70996028\\
0.0432715430435515	-0.61484593\\
0.0540065905073232	-0.51366335\\
0.0645463248171926	-0.45755571\\
0.0751641707964003	-0.3982659\\
0.0858469788701695	-0.34963945\\
0.0965736701831975	-0.31233677\\
0.107334525328225	-0.26894519\\
0.118125666744993	-0.24395533\\
0.139338344171045	-0.19717513\\
0.160604416633501	-0.15355468\\
0.181915867664454	-0.12352144\\
0.20325189717262	-0.10003515\\
0.22460394134366	-0.074197888\\
0.245970138231719	-0.06115276\\
0.270840202328811	-0.048446424\\
0.288163415534036	-0.039594177\\
0.330353196467016	-0.020607136\\
0.372540149518986	-0.013994689\\
0.414705453636153	-0.0069990898\\
0.456845877817938	-0.0054265885\\
0.498964822744851	-0.005914703\\
0.540856724574902	-0.0087178107\\
0.58272264026283	-0.012847953\\
0.624565450143566	-0.018607693\\
0.666383815275177	-0.025334435\\
0.702424920789416	-0.032731742\\
0.708177512714931	-0.033912245\\
0.749947839644273	-0.043142859\\
0.791692883136026	-0.056659851\\
0.833411352948819	-0.070989475\\
0.875102766723219	-0.10037051\\
0.91676309846275	-0.11439648\\
0.958389181247267	-0.18403201\\
1.00000001776646	-0.31565344\\
};
\addplot [color=green, line width=1.3pt, mark=square, mark options={solid, green}]
  table[row sep=crcr]{%
-2.10067467573863e-18	1.4598382\\
-1.24028267280665e-18	1.4598382\\
0.00802630116672502	2.1460655\\
0.0225008352983609	2.3131003\\
0.0379779687597454	1.8561957\\
0.0540064086856373	1.5658172\\
0.0540064126417008	1.5658172\\
0.0858467048895266	1.3703508\\
0.118125304877199	1.1283838\\
0.181915310111255	0.94024694\\
0.245969396284175	0.78853732\\
0.309257753006218	0.66285861\\
0.372539036945721	0.56038421\\
0.435777248266454	0.47060874\\
0.498963362537238	0.3916581\\
0.561791460283731	0.32008207\\
0.624564409023061	0.25389338\\
0.687282677104425	0.19028479\\
0.749947102206648	0.12739223\\
0.812554831581903	0.057264138\\
0.87510234979763	-0.015627971\\
0.93758148876842	-0.17929213\\
1.00000001776646	-0.3090359\\
};
\addlegendentry{LINEAR REFINED}

\addplot [color=green, line width=1.3pt, mark=square, mark options={solid, green}, forget plot]
  table[row sep=crcr]{%
-2.10067467573863e-18	1.4598382\\
-1.24028267280665e-18	1.4598382\\
0.00802633698660474	-0.81932634\\
0.0225009205583868	-0.91294062\\
0.0379780992901577	-0.62716496\\
0.0540065905073232	-0.4787676\\
0.0858469719470585	-0.33804163\\
0.118125665562518	-0.20082845\\
0.181915867664454	-0.11522083\\
0.245970138231719	-0.03836308\\
0.309258645384469	-0.009267698\\
0.372540149518986	0.0090931179\\
0.43577853141086	0.015944762\\
0.498964793074375	0.017096039\\
0.561792719324196	0.013949339\\
0.624565450143566	0.0081310673\\
0.687283604548716	-0.0010691882\\
0.749947836688085	-0.0094432402\\
0.812555401914088	-0.034114793\\
0.875102766723219	-0.015350613\\
0.937581715749889	-0.14156334\\
1.00000001776646	-0.3090359\\
};
\addplot [color=mycolor1, line width=1.3pt, mark=square, mark options={solid, mycolor1}]
  table[row sep=crcr]{%
-2.10067467573863e-18	0.036717705\\
-1.24028267280665e-18	1.0140986\\
0.00802647451397611	2.3725977\\
0.0225008352983609	2.1128097\\
0.0379779684641266	1.687606\\
0.0540064086856373	1.6154777\\
0.0698462598637024	1.3953323\\
0.0858467048895266	1.3121706\\
0.101949780581038	1.2095926\\
0.118125304285961	1.1254386\\
0.118125304877199	1.1253982\\
0.149964858434306	1.0159061\\
0.181915310111255	0.91010666\\
0.213925407485033	0.83552575\\
0.245969395988557	0.7534489\\
0.245969396284175	0.7533409\\
0.309257753006218	0.64402771\\
0.372539036945721	0.53834796\\
0.435777248266454	0.45549548\\
0.498963362537238	0.37747881\\
0.561791460283731	0.30999428\\
0.624564409023061	0.24771619\\
0.687282677104425	0.18639621\\
0.749947102206648	0.13422219\\
0.812554831581903	0.063259907\\
0.87510234979763	0.019904509\\
0.93758148876842	-0.12129729\\
1.00000001776646	-0.2801615\\
};
\addlegendentry{QUADRATIC}

\addplot [color=mycolor1, line width=1.3pt, mark=square, mark options={solid, mycolor1}, forget plot]
  table[row sep=crcr]{%
-2.10067467573863e-18	0.036717705\\
-1.24028267280665e-18	1.0140986\\
0.00802651036116693	-0.97565609\\
0.0225009205583868	-0.69419378\\
0.0379780992901577	-0.65719736\\
0.0540065905073232	-0.43300971\\
0.0698464883133735	-0.40116966\\
0.0858469719470585	-0.29656789\\
0.101950098464154	-0.26007688\\
0.118125666744993	-0.20275322\\
0.149965312872185	-0.14389247\\
0.181915867664454	-0.087783165\\
0.213926045916388	-0.059768148\\
0.245970138231719	-0.030141246\\
0.309258645384469	-0.0050371056\\
0.372540149518986	0.015182843\\
0.43577853141086	0.014717768\\
0.498964793074375	0.016934801\\
0.561792719324196	0.0073482529\\
0.624565450143566	0.0013127846\\
0.687283604548716	-0.015242588\\
0.749947836688085	-0.024871675\\
0.812555401914088	-0.057220351\\
0.875102766723219	-0.066563718\\
0.937581715749889	-0.16339183\\
1.00000001776646	-0.2801615\\
};
\end{axis}

\node[align=center,font=\large, yshift=2em] (title) 
    at (current bounding box.north)
    {$\alpha$ = 8.5 ${^\circ}$, $\Lambda$ = 4\\ $C_p$ distribution, y/2c = 0.211};
\end{tikzpicture}}    
}
&
\resizebox{0.33\linewidth}{!}{
\centerline{
%
%
\definecolor{mycolor1}{rgb}{1.00000,0.00000,1.00000}%
\begin{tikzpicture}

\begin{axis}[%
width=4.602in,
height=3.506in,
at={(0.772in,0.473in)},
scale only axis,
xmin=-0.2,
xmax=1.2,
xlabel style={font=\color{white!15!black}},
xlabel={\textbf{$x/c$}},
ymin=-1,
ymax=2.5,
ylabel style={font=\color{white!15!black}},
ylabel={\textbf{$-Cp$}},
axis background/.style={fill=white},
xmajorgrids,
ymajorgrids,
legend style={legend cell align=left, align=left, draw=white!15!black}
]
\addplot [color=blue, line width=1.3pt, mark=square, mark options={solid, blue}]
  table[row sep=crcr]{%
-1.14213807692232e-18	0.65621138\\
-9.00182759915995e-19	0.37439698\\
0.00409915120658076	2.1641905\\
0.0122539238823407	2.2258725\\
0.0222577349425719	1.8049899\\
0.0324245843031084	1.7769955\\
0.0428809726444962	1.5844892\\
0.0535441871159947	1.5125289\\
0.0640214169077506	1.3869342\\
0.0745773459333662	1.3147616\\
0.0851987608394721	1.2402364\\
0.0958644849257839	1.1871358\\
0.10656477265969	1.1334383\\
0.117295678150729	1.0858679\\
0.138404899644579	1.0056758\\
0.156312530331545	0.94695622\\
0.159568392679639	0.93627512\\
0.18077816392092	0.87361723\\
0.202013450893539	0.82220846\\
0.223265470155025	0.76919907\\
0.244532375676315	0.72515655\\
0.244532376267552	0.7236799\\
0.286836489337739	0.64750582\\
0.329134900924747	0.57508403\\
0.371428170119577	0.51349974\\
0.413697122327436	0.45595461\\
0.45593837235865	0.4042204\\
0.498155409939764	0.35633692\\
0.54012945602504	0.31232253\\
0.582076224298323	0.27090472\\
0.623998475601842	0.2315651\\
0.665894805773255	0.19406715\\
0.707765086521932	0.15704507\\
0.749610722062309	0.12105686\\
0.791416427729308	0.082286686\\
0.833193887244142	0.043075569\\
0.874942343553201	-0.010725223\\
0.916658007270825	-0.052915532\\
0.958337596455102	-0.15310237\\
1.00000001776646	-0.31396174\\
};
\addlegendentry{CUBIC}

\addplot [color=blue, line width=1.3pt, mark=square, mark options={solid, blue}, forget plot]
  table[row sep=crcr]{%
-1.14213807692232e-18	0.65621138\\
-9.00182759915995e-19	0.37439698\\
0.00409922801405183	-0.90321326\\
0.0122541180827971	-0.94546849\\
0.0222580527089917	-0.84152615\\
0.0324250217720263	-0.65002596\\
0.0428815270818934	-0.55367446\\
0.0535448603023033	-0.4516069\\
0.0640222146110596	-0.39665622\\
0.0745782566575152	-0.33818033\\
0.0851997861401124	-0.29093486\\
0.0958656360681177	-0.2544843\\
0.106566040607035	-0.2124259\\
0.117297069652676	-0.18827038\\
0.138406502367362	-0.14385739\\
0.156275732160386	-0.10838422\\
0.159570254559709	-0.10182549\\
0.180780255863664	-0.074527323\\
0.202015770367388	-0.051268533\\
0.223268009443464	-0.028543672\\
0.244535159657687	-0.016292769\\
0.286839243558578	0.004215423\\
0.329137659663322	0.019813737\\
0.371430944541873	0.025097584\\
0.413699902488846	0.029555002\\
0.45594120534161	0.029271344\\
0.498158217198355	0.026891135\\
0.540132045456626	0.022306906\\
0.582078533250478	0.016397014\\
0.624000578543378	0.0088933362\\
0.66589667756479	0.00033328607\\
0.707766754594174	-0.010108377\\
0.749612045428836	-0.021341622\\
0.791417608498584	-0.036946569\\
0.833194756254551	-0.053708084\\
0.874943053288412	-0.085467763\\
0.916658454279872	-0.10386589\\
0.958337804698634	-0.17828794\\
1.00000001776646	-0.31396174\\
};
\addplot [color=green, line width=1.3pt, mark=square, mark options={solid, green}]
  table[row sep=crcr]{%
-2.02434943799647e-18	1.1072608\\
-1.24722130905426e-18	1.1072608\\
0.00791694185315451	1.8170447\\
0.021631210322021	2.048521\\
0.0222577349425719	2.059083\\
0.037623779053015	1.6656756\\
0.0535441837511689	1.4123451\\
0.0535441874116135	1.4123449\\
0.0851987611350909	1.240484\\
0.117295679333204	1.0216686\\
0.117295679776632	1.0216686\\
0.180778165251204	0.85008788\\
0.244532377302218	0.71023911\\
0.307986327775636	0.59273106\\
0.371428170119577	0.49713197\\
0.43482098314434	0.41351634\\
0.498155410235383	0.34060201\\
0.561106099440338	0.27460366\\
0.623998475601842	0.21434057\\
0.686833105022484	0.1568047\\
0.749610720879834	0.10033423\\
0.81230864516365	0.037687205\\
0.874942343553201	-0.027540985\\
0.93750352261668	-0.1726097\\
1.00000001776646	-0.28707844\\
};
\addlegendentry{LINEAR REFINED}

\addplot [color=green, line width=1.3pt, mark=square, mark options={solid, green}, forget plot]
  table[row sep=crcr]{%
-2.02434943799647e-18	1.1072608\\
-1.24722130905426e-18	1.1072608\\
0.00791707384118923	-0.85768664\\
0.0216251717040897	-0.84868509\\
0.0222580527089917	-0.84828395\\
0.0376242751192349	-0.5654335\\
0.0535448636671292	-0.41579044\\
0.0851997924719858	-0.27690333\\
0.117297068618011	-0.14476109\\
0.180780253794332	-0.064885795\\
0.244535158475212	0.0076795914\\
0.307989098524936	0.033125121\\
0.371430974212348	0.048334397\\
0.43482378253549	0.052320223\\
0.498158217198355	0.050691228\\
0.561108570904566	0.045100786\\
0.624000578543378	0.036723182\\
0.686834816527729	0.024994286\\
0.749612045428836	0.013297787\\
0.812309686505364	-0.013772331\\
0.874943053288412	-0.005267804\\
0.937503948149521	-0.13057314\\
1.00000001776646	-0.28707844\\
};
\addplot [color=mycolor1, line width=1.3pt, mark=square, mark options={solid, mycolor1}]
  table[row sep=crcr]{%
-2.02434943799647e-18	0.71826309\\
-1.24722130905426e-18	-0.16234292\\
0.00792250348087713	2.0466621\\
0.0216616442372479	1.918523\\
0.0222577349425719	1.9129419\\
0.037623857343937	1.5233558\\
0.0535441874116135	1.4751093\\
0.0692905640088142	1.2718035\\
0.0851987608394721	1.2005014\\
0.101210587317448	1.1060445\\
0.117295678150729	1.0304499\\
0.117295678741966	1.0303353\\
0.148980421080693	0.92922091\\
0.180778164216538	0.83336627\\
0.212637465244411	0.7616663\\
0.244532377302218	0.69070745\\
0.244532405642409	0.68998605\\
0.307986327775636	0.58330077\\
0.371428170119577	0.48484957\\
0.43482098314434	0.40726018\\
0.498155410235383	0.33474118\\
0.561106099440338	0.27274638\\
0.623998475601842	0.2155547\\
0.686833105022484	0.15979376\\
0.749610720879834	0.11234608\\
0.81230864516365	0.047948029\\
0.874942343553201	0.0075814654\\
0.93750352261668	-0.121348\\
1.00000001776646	-0.26613307\\
};
\addlegendentry{QUADRATIC}

\addplot [color=mycolor1, line width=1.3pt, mark=square, mark options={solid, mycolor1}, forget plot]
  table[row sep=crcr]{%
-2.02434943799647e-18	0.71826309\\
-1.24722130905426e-18	-0.16234292\\
0.00792263853941618	-0.99438637\\
0.0216559485803761	-0.65170497\\
0.0222580527089917	-0.63669246\\
0.037624358322658	-0.60188627\\
0.0535448642583668	-0.3741847\\
0.0692914216916997	-0.34536976\\
0.0851997930632234	-0.24207138\\
0.101211799903776	-0.20738861\\
0.117297058727852	-0.15188091\\
0.117297069652676	-0.1519822\\
0.148982167936896	-0.095866509\\
0.180780254828998	-0.042281218\\
0.212639906422036	-0.016766028\\
0.244535158475212	0.012036184\\
0.244535159657687	0.01157608\\
0.307989098524936	0.033080928\\
0.371430974212348	0.049800787\\
0.43482378253549	0.046861839\\
0.498158217198355	0.046314884\\
0.561108570904566	0.034449816\\
0.624000578543378	0.026054965\\
0.686834816527729	0.0072571211\\
0.749612045428836	-0.0049986853\\
0.812309686505364	-0.039559159\\
0.874943053288412	-0.053601913\\
0.937503948149521	-0.15064526\\
1.00000001776646	-0.26613307\\
};
\end{axis}

\node[align=center,font=\large, yshift=2em] (title) 
    at (current bounding box.north)
    {$\alpha$ = 8.5 ${^\circ}$, $\Lambda$ = 4\\ $C_p$ distribution, y/2c = 0.611};
\end{tikzpicture}%
}    
}
&
\resizebox{0.33\linewidth}{!}{
\centerline{
%
%
\definecolor{mycolor1}{rgb}{1.00000,0.00000,1.00000}%
\begin{tikzpicture}

\begin{axis}[%
width=4.602in,
height=3.506in,
at={(0.772in,0.473in)},
scale only axis,
xmin=-0.2,
xmax=1.2,
xlabel style={font=\color{white!15!black}},
xlabel={\textbf{$x/c$}},
ymin=-1,
ymax=2,
ylabel style={font=\color{white!15!black}},
ylabel={\textbf{$-Cp$}},
axis background/.style={fill=white},
xmajorgrids,
ymajorgrids,
legend style={legend cell align=left, align=left, draw=white!15!black}
]
\addplot [color=blue, line width=1.3pt, mark=square, mark options={solid, blue}]
  table[row sep=crcr]{%
-1.54963118995704e-18	0.27725598\\
-8.18760392295643e-19	0.050570723\\
0.00392701810251808	1.683823\\
0.0118442665829781	1.8616087\\
0.02161567683683	1.5279721\\
0.0316037413049256	1.5323106\\
0.041889270345907	1.3722105\\
0.0523878704058813	1.3155442\\
0.062737225328494	1.2070961\\
0.0731691518516966	1.1455313\\
0.0836701513638417	1.0806531\\
0.0849760738031609	1.0749229\\
0.0942190294122309	1.0343511\\
0.104805520422297	0.98717183\\
0.115425593582949	0.94503903\\
0.115425600234373	0.94368833\\
0.136558476761236	0.87247592\\
0.157744220509282	0.80973178\\
0.178974295583442	0.75373411\\
0.20022751538331	0.70550197\\
0.221494575843468	0.66106158\\
0.242773583866101	0.60630143\\
0.242773584309529	0.60238242\\
0.285590235855046	0.54779869\\
0.328402874894375	0.47980261\\
0.371212031112153	0.42504907\\
0.413994302468823	0.37223929\\
0.456750730146968	0.32627368\\
0.499484826756434	0.28434667\\
0.541433255942272	0.24528919\\
0.583358064905119	0.20983364\\
0.625262317613125	0.17608741\\
0.667138723424506	0.14419271\\
0.708992818345403	0.11313544\\
0.750825875814118	0.08308211\\
0.792446480935992	0.050622411\\
0.834033637606829	0.018023754\\
0.875586606323125	-0.028062705\\
0.917093432404952	-0.063333198\\
0.958558266919661	-0.15377553\\
1.00000001776646	-0.30251953\\
};
\addlegendentry{CUBIC}

\addplot [color=blue, line width=1.3pt, mark=square, mark options={solid, blue}, forget plot]
  table[row sep=crcr]{%
-1.54963118995704e-18	0.27725598\\
-8.18760392295643e-19	0.050570723\\
0.00392698310241853	-0.96461976\\
0.0118441849778742	-0.89918256\\
0.0216155522488625	-0.77713251\\
0.0316035935890773	-0.57254791\\
0.0418891030807539	-0.47496161\\
0.0523876918907786	-0.37278166\\
0.0627370415134146	-0.31906211\\
0.0731689711406572	-0.26170692\\
0.0836699792406074	-0.21687232\\
0.0849714104316008	-0.21237987\\
0.0942188726148191	-0.18067032\\
0.104805373172876	-0.14172307\\
0.115425473941253	-0.11810376\\
0.115425488104812	-0.11850182\\
0.136559022474748	-0.075128399\\
0.157745396319784	-0.036952913\\
0.178976100612764	-0.0095137721\\
0.200229943317113	0.0082446421\\
0.221497620300509	0.033464823\\
0.242777227017825	0.034172844\\
0.242777228052491	0.035677571\\
0.28559351571825	0.055495754\\
0.328405816007735	0.070014924\\
0.371214613721441	0.070820905\\
0.413996528466602	0.074265525\\
0.456752599076327	0.071732879\\
0.499486384214262	0.068782106\\
0.541434606172042	0.061834406\\
0.583359293377541	0.055074498\\
0.625263311841247	0.046879616\\
0.667139618151478	0.037408415\\
0.708993684510187	0.026260639\\
0.750826620642305	0.014236098\\
0.792447165095704	-0.0026290566\\
0.834034107293575	-0.020926602\\
0.875586997959648	-0.055488463\\
0.917093708101339	-0.076441534\\
0.958558402145052	-0.15676838\\
1.00000001776646	-0.30251953\\
};
\addplot [color=green, line width=1.3pt, mark=square, mark options={solid, green}]
  table[row sep=crcr]{%
-1.8300675764247e-18	0.69961536\\
-8.17025388750979e-19	0.69961536\\
0.00762118183428853	1.4286324\\
0.0216156801451465	1.7494626\\
0.0367163411875874	1.4345329\\
0.052387891779524	1.2196158\\
0.0836701736134814	1.0755891\\
0.115425638221391	0.8799746\\
0.178974315685522	0.72479415\\
0.242773584913838	0.59754503\\
0.306997107259851	0.48703215\\
0.371212031703391	0.40178329\\
0.435375487926569	0.32146326\\
0.499484827347671	0.26222357\\
0.562398376058854	0.20486018\\
0.625262259454648	0.15417184\\
0.68806843081763	0.10684429\\
0.75082593515507	0.061058957\\
0.813244173284663	0.010045909\\
0.875586604845031	-0.039834339\\
0.937832364541668	-0.15963557\\
1.00000001776646	-0.2603251\\
};
\addlegendentry{LINEAR REFINED}

\addplot [color=green, line width=1.3pt, mark=square, mark options={solid, green}, forget plot]
  table[row sep=crcr]{%
-1.8300675764247e-18	0.69961536\\
-8.17025388750979e-19	0.69961536\\
0.00762112051832905	-0.8832413\\
0.0216155581373077	-0.76263648\\
0.0367161815675636	-0.48516461\\
0.0523877121014522	-0.33379793\\
0.0836699947670298	-0.19796844\\
0.115425519940142	-0.072435357\\
0.178976110328969	-0.0010881792\\
0.242777236894912	0.067036182\\
0.307000230711601	0.085709527\\
0.371214643391917	0.10032669\\
0.43537756642936	0.098087557\\
0.499486384214262	0.093494236\\
0.562399694920442	0.087495267\\
0.625263313319341	0.077358484\\
0.688069300353065	0.065888859\\
0.750826679983257	0.054340661\\
0.813244725201233	0.028187478\\
0.875586997959648	0.031608168\\
0.937832562562223	-0.10279376\\
1.00000001776646	-0.2603251\\
};
\addplot [color=mycolor1, line width=1.3pt, mark=square, mark options={solid, mycolor1}]
  table[row sep=crcr]{%
-1.8300675764247e-18	0.36683193\\
-8.17025388750979e-19	-0.37652409\\
0.00762587522916001	1.6427321\\
0.0216156801451465	1.6563702\\
0.0367164841166687	1.3141203\\
0.0523878768747048	1.2903523\\
0.0679438165567362	1.1080149\\
0.0836701550590769	1.0489005\\
0.099507816943524	0.96475548\\
0.115425594617615	0.8936277\\
0.115425609546365	0.89570719\\
0.147145272838322	0.80942917\\
0.178974300461153	0.69992328\\
0.210859384162065	0.67793906\\
0.242773584913838	0.57264483\\
0.24277358652667	0.54137844\\
0.30699710740766	0.49127528\\
0.371212031703391	0.40065923\\
0.435375487926569	0.33123535\\
0.499484827347671	0.26715472\\
0.562398376058854	0.21390681\\
0.625262259454648	0.16588049\\
0.68806843081763	0.11983029\\
0.75082593515507	0.082076438\\
0.813244173284663	0.02823429\\
0.875586604845031	8.0749451e-05\\
0.937832366019762	-0.12003584\\
1.00000001776646	-0.26250795\\
};
\addlegendentry{QUADRATIC}

\addplot [color=mycolor1, line width=1.3pt, mark=square, mark options={solid, mycolor1}, forget plot]
  table[row sep=crcr]{%
-1.8300675764247e-18	0.36683193\\
-8.17025388750979e-19	-0.37652409\\
0.00762581470180202	-0.99838495\\
0.0216155581373077	-0.56289697\\
0.036716320634094	-0.52793169\\
0.0523876969945267	-0.29793435\\
0.0679436303744736	-0.27222234\\
0.0836699898350224	-0.17167231\\
0.0995076731489009	-0.1382522\\
0.115425482648936	-0.089084223\\
0.115425488091741	-0.088131212\\
0.14714611684868	-0.030985828\\
0.178976114776323	0.0072704963\\
0.210862122729364	0.051156536\\
0.242777224800684	0.047870021\\
0.242777236894912	0.058793563\\
0.307000230711601	0.080374181\\
0.371214643391917	0.094828166\\
0.43537756642936	0.086647063\\
0.499486384214262	0.08379633\\
0.562399694920442	0.069619186\\
0.625263313319341	0.060141865\\
0.688069300353065	0.040529653\\
0.750826679983257	0.028392183\\
0.813244725201233	-0.0082420707\\
0.875586997959648	-0.020341894\\
0.937832562562223	-0.13161916\\
1.00000001776646	-0.26250795\\
};
\end{axis}

\node[align=center,font=\large, yshift=2em] (title) 
    at (current bounding box.north)
    {$\alpha$ = 8.5 ${^\circ}$, $\Lambda$ = 4\\ $C_p$ distribution, y/2c = 0.811};
\end{tikzpicture}%
}    
}
\end{tabular}  	
\caption{Development of the pressure coefficient in three different
wing sections. All of the sections were obtained by using planes 
parallel to the chord of the profile and perpendicular to the $y$-axis.
These are located at $y/c$ = 0.211, 0.611 and 0.811.
The comparison was made by analysing four different case studies: starting 
with a base case, the aim is to observe how the solution changes in 
response to grid refinement or the use of higher degree finite elements, 
such as quadratic and cubic.\label{fig:GridSpacing}}
\end{figure}

The second numerical test is aimed at the characterization of the solver accuracy in the estimation of the pressure field on the wing surface. This is
relevant for the present application, in which the main goal of the simulations is the computation of the fluid dynamic forces on the wing. And, as well known,
such forces in the context of potential flow theory, are in fact obtained as the integral of the pressure stresses on the wing surface. The pressure field is
here evaluated by means of Equation \eqref{eq:pressure}, and Figure \ref{fig:GridSpacing} presents plots of the non dimensional pressure coefficient

\begin{equation}
c_p = \frac{p-p_\infty}{\frac{1}{2}\rho V_\infty^2}
\end{equation}
at vertical sections of the wing located at different span wise coordinate. The present convergence test starts from a moderately accurate base
solution obtained with bi-linear elements on a computational grid generated setting maximum aspect ratio 2.5, 3 uniform refinement cycles, and 2
adaptive refinement cycles based on CAD curvature. To evaluate the effectiveness of possible $hp$ refinement strategies, the system degrees of freedom are then
increased either by refining the computational grid, or by increasing the finite element degree.
The three plots on each row of Figure \ref{fig:GridSpacing} depict values of $-c_p$ as a function of non dimensional chord wise coordinate $x/c$,
on vertical sections located at non dimensional span wise coordinate values $y/2c=.211$, $y/2c=.611$ and $y/2c=.811$, from left to right.
In the first row plots, the base solution (blue lines) and is compared with the one obtained with bi-linear elements on a grid featuring an
additional uniform refinement step (green lines). The three diagrams suggest that the bi-linear solution on the finer
grid appears significantly less dissipative than the solution on the coarser grid. With the finer grid, the peak $c_p$ value on the wing windward
side corresponding to the flow stagnation point is approaching the theoretical limit $c_p=1$. In addition, also the suction side peak value is significantly
higher on the refined grid simulation. As a result of the increased peak values located close to the leading edge, the entire solution amplitude up
to the trailing edge appears slightly higher, with a consequent increase of the overall area included between the leeward and windward side $c_p$ lines.
A similar behavior is observed also when the base solution degrees of freedom are increased  by
means of an increment of the discretization degree. The three plots in the second row in of Figure \ref{fig:GridSpacing} refer to the comparison between
the base bi-linear solution $c_p$ values, and the pressure coefficient values obtained making use of bi-quadratic finite elements on the same grid.
In the third row plots instead, the base solution $c_p$ field is compared to the one computed making use of the same grid and bi-cubic finite elements.
In both cases, the results show that a refinement strategy based on the increase of the discretization degree is effective in reducing numerical diffusion.
The peak values in the leading edge region, including the one associated with the stagnation point, are in fact sharper, and as a result the overall area
between the leeward and windward curves is increased. The plots in the final row of of Figure \ref{fig:GridSpacing} finally compare the $c_p$ solutions
obtained with the different refinement strategies tested. As a general trend, the diagrams indicate that the solution obtained with the three refinement
strategies appear quite similar. Once again, this seems to suggest that the solver solutions converge as the number of degrees of freedom in the system
is increased. A more detailed look at the plots indicates that the bi-linear solution on the finer grid and the bi-quadratic solution on the coarser grid are
particularly close to each other across all the regions and sections. Instead, the bi-cubic solution seems to have less numerical dissipation, as its peaks
present higher values with respect to the bi-linear and bi-quadratic ones. This could be associated to the fact that, in the case of the present
numerical test, the cubic solution has a slightly higher amount of degrees of freedom with respect to the bi-linear and bi-quadratic solutions, which
share the same spacing among the collocation points.

\subsection{Wake relaxation}

After the convergence and accuracy features of the nonlinear problem resolution algorithm have been described and discussed,
the present section focuses on highlighting and discussing the performance of the wake relaxation algorithm described in
Section \ref{sec:wake_relaxation}. Also in this case, we make use of the geometrical layout described in
Figure \ref{fig:wingSetup}, and compare the results obtained with different discretization degrees. The grid used for the first
bi-linear finite elements test here considered is set up prescribing maximum cell aspect ratio 2.5, 2 uniform refinement cycles,
and 1 adaptive refinement cycle based on CAD curvature. To obtain an even comparison in terms of collocation point distribution
with respect to the bi-linear test case, the grid generated for
the bi-quadratic case features identical refinement parameters, but one less uniform refinement cycle. As for the bi-cubic test case
a similar spacing of the collocation points is obtained imposing maximum cell aspect ratio 1.85, 1 uniform refinement cycle,
and 1 adaptive refinement cycle based on CAD curvature. In all the simulations discussed in the present sections, the wake relaxation
algorithm has been set up with a prescribed wake cell length (see Equation \eqref{eq:wake_relax}) $d=0.15\ $m, and has been arrested after
9 iterations.

\begin{figure}[htb!]
\begin{center}
\includegraphics[width=\textwidth]{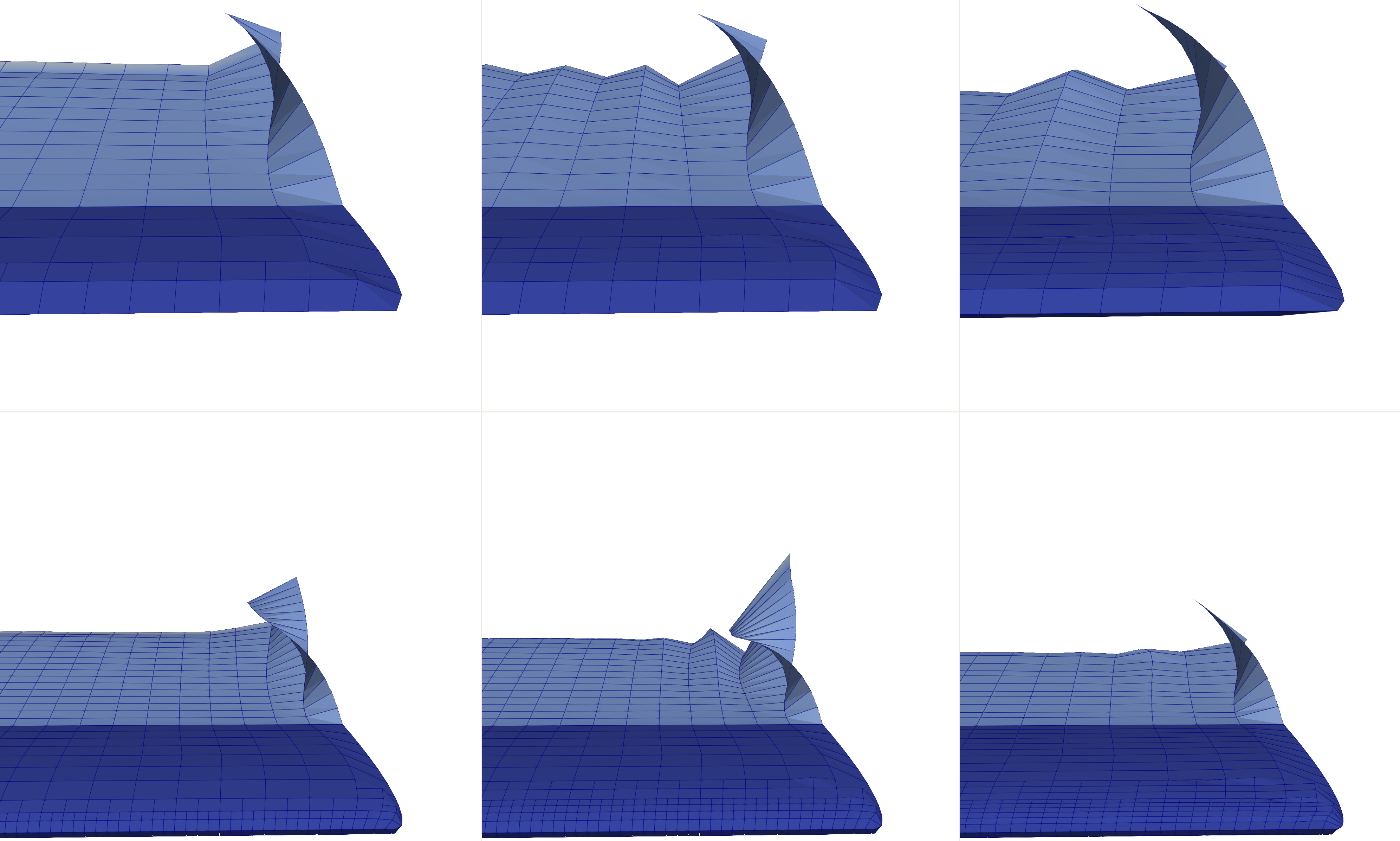}
\end{center}
\caption{A front view of the grid in the --- rounded --- tip region of the wing and of the fully developed wake. The three images
         on the  top refer to the base solution obtained with linear (left), quadratic (middle) and cubic (right) finite elements, respectively.
         The three images on the bottom show the linear (left), quadratic (middle) and cubic (right) solutions obtain with one additional
         uniform refinement level. 
         \label{fig:WingTipFECompare}}
\end{figure}

The final wake geometries obtained with the grids just described are depicted in the top row of Figure \ref{fig:WingTipFECompare}.
The left, right and middle images refer to the solutions obtained with bi-linear, bi-quadratic and bi-cubic finite elements. We point out
that, in order to show all the collocation points of the solutions obtained with higher degree finite elements, in all the following wake geometry
visualizations the higher order cells have been split in multiple bi-linear cells. In this way, each bi-quadratic cell is represented in the visualizations
as a set of 4 bi-linear cells, whereas each bi-cubic cell is visualized by means of 9 bi-linear cells. The wake plots
confirm that the wake velocity computation and wake relaxation algorithm used lead to a qualitatively correct representation of the
wake geometry. The images in fact clearly show the presence of a vortex developing past the wing tip, which curls upward the lateral edge of the
wake vortex sheet. Instead, in the region located downstream with respect to the wing middle portion, the wake appears flat and aligned with
the trailing edge slope. We point out that in the framework of the wake relaxation algorithm developed the wake orientation at the trailing
edge is not imposed, but is instead a result of the computations. Thus, the latter observation is quite satisfactory. At a closer look, the
wake geometries obtained with bi-linear and bi-quadratic finite elements appear very similar, especially in the wake roll up region. The bi-cubic
solution also appears similar to that obtained with lower degree finite elements, despite being computed making use of a different distribution
and number of collocation points. In the central portion of the wake, the bi-quadratic and bi-cubic solutions differ from their bi-linear counterpart,
as they both present a series of oscillations. The plots in the bottom row of \ref{fig:WingTipFECompare}, which have been obtained adding a uniform
refinement step to the grids in the top row, confirm the presence of these oscillations. On the finer grid here considered, the oscillations in the
bi-quadratic solution appear more visible, as they extend to the wake roll up region, which appears wider than the one computed with bi-linear elements.   
These oscillations are the consequence of not fully accurate computations of the velocity at the
wake collocation points. Higher order shape functions change sign multiple times over each cell, which causes cancellation problems when
computing the hypersingular kernel integrals required by the wake velocity Galerkin BEM problem. Clearly, when the oscillations are so severe that
they are able to cause self contact of the wake surface, they induce instabilities that lead to simulations break up. Indeed, the problem
can be significantly reduced both through an increase of the number of quadrature nodes in both the inner and outer loop of the wake
velocity Galerkin BEM problem, and with a choice of shorter cells in the wake. However, we must warn that increasing these parameters leads to
more than linear increment in the computational cost of the wake relaxation algorithm. For such a reason, the execution of the velocity computation and wake relaxation algorithm in the present example were extremely more time consuming when higher order finite elements were employed. In addition, to avoid
instabilities, the bi-cubic simulation on the finer grid required a shorter wake cell length, which in turn resulted in shorter overall wake length.

\begin{figure}[htb!]
\begin{center}
\includegraphics[width=\textwidth]{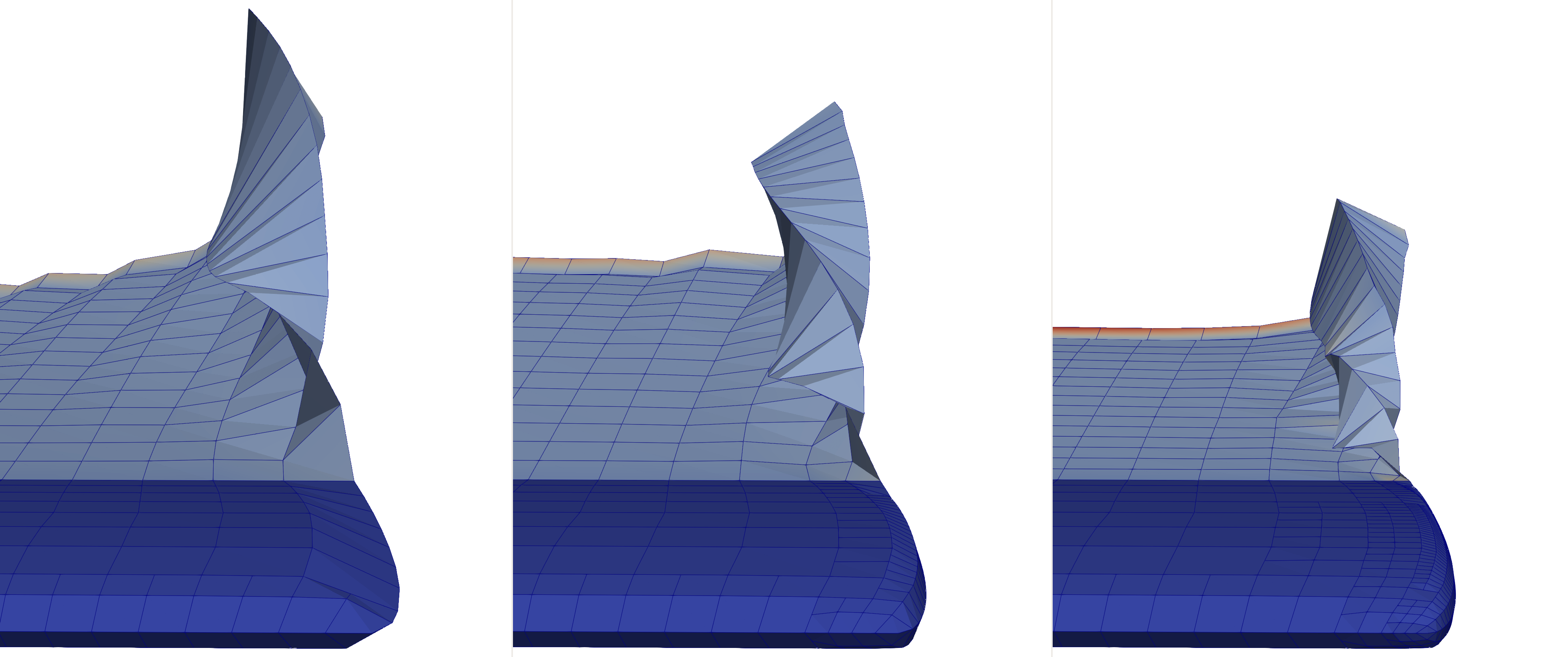}
\end{center}
\caption{A front view of the grid in the --- rounded --- tip region of the wing and of the fully developed wake. The three images
         refer to the solution obtained with the base grid (on the left), with one additional refinement on the tip (in the middle),
         and with two additional refinements on the tip (on the right). 
         \label{fig:WingTipRef}}
\end{figure}

Given the aforementioned considerations, in the remainder of the article we will only present solutions based on bi-linear finite elements, in which
localized grid refinements are used to compensate for the slower convergence rate with respect to higher order elements, highlighted in
Figure \ref{fig:DeltaPhiDistribution}. A combination of $h$-adaptivity and linear elements appear in fact to strike a balance between
accuracy, stability and efficiency of the simulations. Thus, the next test case is designed to showcase the local refinement capabilities of the
solver developed. To this end, three different grids are generated starting from a base grid with maximum cell aspect ratio 1.85, 3 uniform refinement
cycles and 1 adaptive refinement cycle based on CAD curvature. So, the first grid considered is obtained adding one cycle refining cells located in
correspondence with the rounded tip of the wing. Similarly, the second and third computational grids are generated adding to the base grid two and
three local refinement cycles in the tip region, respectively.
Figure \ref{fig:WingTipRef} shows the results obtained using bi-linear elements on the three grids generated. As can be appreciated in the plots,
the rounded tip geometry is reproduced with greater accuracy as the grid is locally refined. This has clearly a consequence on the computed tip
vortex, which appears confined in a smaller region when finer tip grids are considered (right plot), as opposed to the wider curl computed with coarser
grids. The wider wake tip vortex appearing in the coarser mesh simulation (left plot) is also affecting the nearby portion of the wake vortex sheet, which
appears less flat than the one computed with the finer meshes. Clearly, this can have significant effects on the accuracy of the computed downwash caused
by the wing presence. 

\begin{figure}[htb!]
\resizebox{\textwidth}{!}{\input{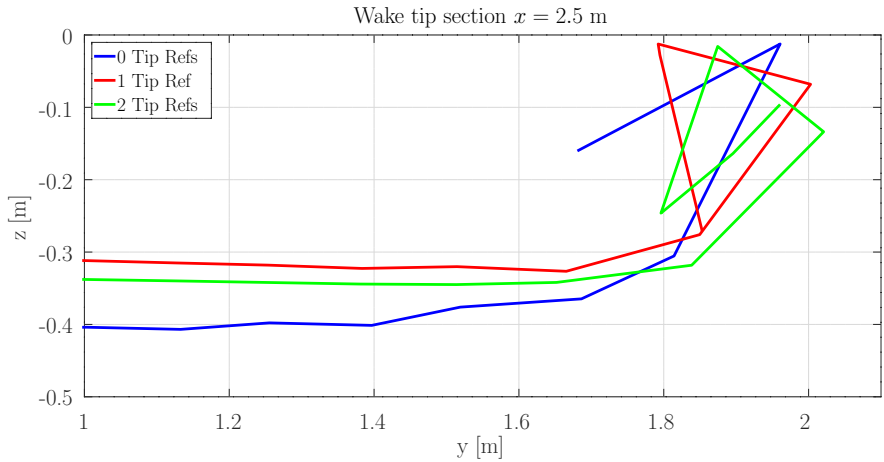}}
\caption{Wing tip refinement effect on a section of the wake located 1.5 chord lengths downstream with respect to the trailing edge.
         The blue, red and green curves in the plot refer to solutions obtained with the base solution, one and two tip refinements, respectively. \label{fig:WingTipRefSection}}
\end{figure}

Further confirmation of the effect that tip refinements can have on the overall wake shape is given by Figure \ref{fig:WingTipRefSection}. The plot
displays a vertical cut of the converged wake geometry on a plane located at constant longitudinal coordinate $x=2.5\ $m. The curves in the plot
show that adding nodes to the final portion of the wake vortex sheet increases the tendency of the wake edge to curl up in a spiraling fashion. As
this results in smaller spiral coils, the influence of the tip vortex presence on the central part of the wake is weakened. Thus, the central
portion of the wake geometry seems to converge to a more flat solution as the tips are refined.

\section{Comparison with experimental data}\label{sec:results}

As a final confirmation of the accuracy of the mathematical model presented in 
this work, the present section is devoted to the comparison of the numerical
results with experimental measurements available in the literature. The test cases considered
are presented in \cite{YipShubertNACA0012} and both feature a semi-wing with constant NACA 0012 airfoil section and chord length $c=1\ $m.
In the first test, the wing has a regular rectangular shape with semi-span $s=2.95c$. The semi-wing is mounted on a wall on one side, and
has a free tip on the other. In the second test, the wing is rotated along the vertical axis $z$ to obtain a swept semi-wing in which
the leading edge is inclined by an angle $\beta=20^\circ$ with respect to the $y$ axis. In such a case, the semi-wing span reported is $s=2.91c$.
In the specific experiments here considered, both wings are set at an angle of attack $\alpha=6.75^\circ$.

\subsubsection{Rectangular wing}

The experimental setup corresponding to the rectangular semi-wing described, is reproduced in the simulations
with a geometric arrangement by all means identical to that of the numerical tests presented in section \ref{sec:convergence}. 
To obtain an equivalent flow on a wing in free air --- as opposed to the experimental one on a semi-wing --- we considered
an aspect ratio of $\Lambda=5.9$, corresponding in our case to a wingspan $b=5.9c=5.9\ $m. 

The geometric CAD model, designed using FreeCad, features in this case a flat cap surface on each wing tip. An initial mesh
made of approximately 10 cells is then refined onto the surface of the CAD model imposing maximum cell aspect ratio 4.5,
2 uniform refinement cycles and 1 adaptive refinement cycle based on CAD curvature. This results in a computational grid
featuring approximately 3800 quadrilateral cells. Bi-linear finite elements are used for this test, which result in  
a resolution system of 4240 degrees of freedom.

\begin{figure}[htb!]
\begin{tabular}{c c}
\resizebox{0.5\linewidth}{!}{
\centerline{
%
%
\begin{tikzpicture}

\begin{axis}[%
width=4.602in,
height=3.506in,
at={(0.772in,0.473in)},
scale only axis,
xmin=0,
xmax=1,
xlabel style={font=\color{white!15!black}},
xlabel={\textbf{$x/c$}},
ymin=-1,
ymax=2.5,
ylabel style={font=\color{white!15!black}},
ylabel={\textbf{$Cp$}},
axis background/.style={fill=white},
xmajorgrids,
ymajorgrids,
legend style={legend cell align=left, align=left, draw=white!15!black}
]
\addplot [color=red, line width=1.3pt, only marks, mark=asterisk, mark options={solid, red}]
  table[row sep=crcr]{%
3.2041e-06	0.26525\\
0.002012150449	1.86434\\
0.005020431025	2.3252\\
0.007996115241	2.24881\\
0.011659896361	2.19847\\
0.016011612369	2.17181\\
0.020884585225	2.1357\\
0.025480459876	1.93394\\
0.043423057924	1.62214\\
0.0676104004	1.39548\\
0.0963978304	1.26595\\
0.126951115204	1.13657\\
0.162640404369	0.98821\\
0.196002769284	0.9135\\
0.234719901441	0.82214\\
0.274486927225	0.73796\\
0.3167100729	0.66327\\
0.361957453641	0.56964\\
0.410960487844	0.5044\\
0.459125497744	0.48188\\
0.506656392804	0.44523\\
0.558261514561	0.39907\\
0.616919135364	0.33624\\
0.674686460449	0.29006\\
0.723164852881	0.22753\\
0.781496664529	0.17906\\
0.8431646976	0.10926\\
0.893962686016	0.05627\\
0.955308805201	-0.0395\\
0.997727292769	-0.11835\\
};
\addlegendentry{EXPERIMENTAL}

\addplot [color=red, line width=1.3pt, only marks, mark=asterisk, mark options={solid, red}, forget plot]
  table[row sep=crcr]{%
3.2041e-06	0.27406\\
0.002012150449	-0.75027\\
0.005020431025	-0.96321\\
0.007996115241	-0.96902\\
0.011659896361	-0.91071\\
0.016011612369	-0.84992\\
0.020884585225	-0.78186\\
0.025480459876	-0.69126\\
0.043423057924	-0.45761\\
0.0676104004	-0.34849\\
0.0963978304	-0.24246\\
0.126951115204	-0.12829\\
0.162640404369	-0.05208\\
0.196002769284	0.00754\\
0.234719901441	0.03099\\
0.274486927225	0.05636\\
0.3167100729	0.09257\\
0.361957453641	0.05976\\
0.410960487844	0.05792\\
0.459125497744	0.07709\\
0.506656392804	0.06886\\
0.558261514561	0.08166\\
0.616919135364	0.07291\\
0.674686460449	0.07188\\
0.723164852881	0.06054\\
0.781496664529	0.03228\\
0.8431646976	0.008\\
0.893962686016	-0.0136\\
0.955308805201	-0.06534\\
0.997727292769	-0.1317\\
};
\addplot [color=blue, line width=1.3pt, mark=square, mark options={solid, blue}]
  table[row sep=crcr]{%
-1.38153269442179e-18	0.37722927\\
0.00139307630005401	1.5409409\\
0.00497626636451105	2.0700188\\
0.00970895471656558	2.2025635\\
0.0150663810641994	1.9665945\\
0.0263588828865331	1.8438565\\
0.0382461606009063	1.5908488\\
0.0620389332725132	1.4355464\\
0.0863171952784596	1.2099799\\
0.136795491697454	1.0349758\\
0.187646818819561	0.90006524\\
0.238626025761697	0.78377837\\
0.289656065840704	0.69065833\\
0.340645114236712	0.60844588\\
0.391616123512058	0.53602999\\
0.442556283598292	0.47060904\\
0.49346149431994	0.41096297\\
0.544274601830179	0.35575721\\
0.595050981677663	0.30405435\\
0.645790889390751	0.25478312\\
0.696495128302566	0.20701948\\
0.747162724736705	0.15876324\\
0.797793651146356	0.11047266\\
0.848385571855541	0.051222619\\
0.898934953143353	0.0019913816\\
0.924216456524101	-0.056182317\\
0.949484926671784	-0.07761658\\
0.974739587626676	-0.23945583\\
};
\addlegendentry{RELAX}

\addplot [color=blue, line width=1.3pt, mark=square, mark options={solid, blue}, forget plot]
  table[row sep=crcr]{%
-1.38153269442179e-18	0.37722927\\
0.00139078451033227	-0.63309056\\
0.00497199352608279	-0.98835063\\
0.00970306388271443	-0.9886604\\
0.0150591283867115	-0.86684859\\
0.0263494907835121	-0.71942151\\
0.0382350648782204	-0.53338349\\
0.0620252922726972	-0.39076266\\
0.0863016539364379	-0.23921859\\
0.136777281480813	-0.14231272\\
0.187627005760579	-0.052414771\\
0.238605384861892	-0.015930912\\
0.289635120095325	0.0092724785\\
0.340624340397061	0.021090131\\
0.391595820677051	0.026191318\\
0.442536743964811	0.026178654\\
0.49344299232772	0.022809621\\
0.544257325499687	0.016908167\\
0.595035095535824	0.0090880217\\
0.645776581449922	-0.00044911631\\
0.696482442542836	-0.011687431\\
0.747151897753328	-0.025511127\\
0.797784720719302	-0.041528553\\
0.848378719422441	-0.066671208\\
0.898930207728643	-0.091461986\\
0.924212849418534	-0.12530527\\
0.949482533932792	-0.12742102\\
0.974738362879276	-0.26022965\\
1.00000000707044	-0.4098599\\
};
\addplot [color=green, line width=1.3pt, mark=square, mark options={solid, green}]
  table[row sep=crcr]{%
-1.38153269442179e-18	0.36693737\\
0.00139307630005401	1.52657\\
0.00497626636451105	2.0553162\\
0.00970895471656558	2.188772\\
0.0150663810641994	1.9549133\\
0.0263588828865331	1.8334908\\
0.0382461606009063	1.5821155\\
0.0620389332725132	1.4278113\\
0.0863171952784596	1.2033885\\
0.136795491697454	1.0291023\\
0.187646818819561	0.89474159\\
0.238626025761697	0.77882242\\
0.289656065840704	0.68597317\\
0.340645114236712	0.60397184\\
0.391616123512058	0.53172648\\
0.442556283598292	0.46644777\\
0.49346149431994	0.40692315\\
0.544274601830179	0.3518241\\
0.595050981677663	0.300217\\
0.645790889390751	0.25103661\\
0.696495128302566	0.20335461\\
0.747162724736705	0.15520416\\
0.797793651146356	0.10693689\\
0.848385571855541	0.048036706\\
0.898934953143353	-0.0020694567\\
0.924216456524101	-0.055619821\\
0.949484926671784	-0.095849514\\
0.974739587626676	-0.20466214\\
};
\addlegendentry{FLAT}

\addplot [color=green, line width=1.3pt, mark=square, mark options={solid, green}, forget plot]
  table[row sep=crcr]{%
-1.38153269442179e-18	0.36693737\\
0.00139078451033227	-0.63767278\\
0.00497199352608279	-0.9892115\\
0.00970306388271443	-0.98797148\\
0.0150591283867115	-0.86558473\\
0.0263494907835121	-0.71802211\\
0.0382350648782204	-0.53222597\\
0.0620252922726972	-0.38988423\\
0.0863016539364379	-0.23879477\\
0.136777281480813	-0.14228751\\
0.187627005760579	-0.052759621\\
0.238605384861892	-0.016533662\\
0.289635120095325	0.0084552411\\
0.340624340397061	0.020094447\\
0.391595820677051	0.025039969\\
0.442536743964811	0.024888182\\
0.49344299232772	0.021391701\\
0.544257325499687	0.01537133\\
0.595035095535824	0.0074378136\\
0.645776581449922	-0.0022070541\\
0.696482442542836	-0.013555963\\
0.747151897753328	-0.027466597\\
0.797784720719302	-0.043650042\\
0.848378719422441	-0.068661459\\
0.898930207728643	-0.094463423\\
0.924212849418534	-0.12387383\\
0.949482533932792	-0.14475088\\
0.974738362879276	-0.2248922\\
1.00000000707044	-0.2961047\\
};
\end{axis}

\node[align=center,font=\large, yshift=2em] (title) 
    at (current bounding box.north)
    {$\alpha$ = 6.75 ${^\circ}$, $\Lambda$ = 6\\ $C_p$ distribution, y/2c = 0.200};
\end{tikzpicture}%
}    
}
&
\resizebox{0.5\linewidth}{!}{
\centerline{
%
%
\begin{tikzpicture}

\begin{axis}[%
width=4.602in,
height=3.506in,
at={(0.772in,0.473in)},
scale only axis,
xmin=0,
xmax=1,
xlabel style={font=\color{white!15!black}},
xlabel={\textbf{$x/c$}},
ymin=-1,
ymax=2.5,
ylabel style={font=\color{white!15!black}},
ylabel={\textbf{$Cp$}},
axis background/.style={fill=white},
xmajorgrids,
ymajorgrids,
legend style={legend cell align=left, align=left, draw=white!15!black}
]
\addplot [color=red, line width=1.3pt, only marks, mark=asterisk, mark options={solid, red}]
  table[row sep=crcr]{%
0.00116	1.57781\\
0.00465	1.80779\\
0.0093	2.09055\\
0.01395	2.18347\\
0.01453	2.18773\\
0.01919	2.1995\\
0.02384	2.15825\\
0.02442	2.13678\\
0.02733	1.99624\\
0.0343	1.89607\\
0.04128	1.78547\\
0.05174	1.62798\\
0.06628	1.45656\\
0.07791	1.35102\\
0.09244	1.25661\\
0.10581	1.19777\\
0.12442	1.13496\\
0.13779	1.0887\\
0.16047	0.997\\
0.17151	0.95605\\
0.19709	0.88865\\
0.2064	0.86758\\
0.23605	0.79621\\
0.24593	0.77321\\
0.27209	0.71911\\
0.2843	0.69651\\
0.31744	0.63604\\
0.32267	0.62571\\
0.36163	0.54288\\
0.37209	0.51952\\
0.40581	0.45881\\
0.41163	0.4546\\
0.45233	0.44729\\
0.46337	0.44267\\
0.50407	0.39355\\
0.50988	0.38657\\
0.55407	0.34967\\
0.56163	0.34223\\
0.60523	0.28491\\
0.61221	0.27728\\
0.66221	0.24138\\
0.66919	0.23588\\
0.72151	0.17461\\
0.72209	0.17387\\
0.77616	0.10935\\
0.77907	0.10595\\
0.83314	0.03932\\
0.83837	0.03214\\
0.88605	-0.03869\\
0.88721	-0.0403\\
0.94128	-0.1145\\
0.94186	-0.11544\\
0.97674	-0.19061\\
0.97733	-0.19213\\
};
\addlegendentry{EXPERIMENTAL}

\addplot [color=red, line width=1.3pt, only marks, mark=asterisk, mark options={solid, red}, forget plot]
  table[row sep=crcr]{%
0.00116	-0.80781\\
0.00465	-1.014\\
0.0093	-0.93154\\
0.01395	-0.86381\\
0.01453	-0.85583\\
0.01919	-0.79313\\
0.02384	-0.7178\\
0.02442	-0.70773\\
0.02733	-0.65833\\
0.0343	-0.55732\\
0.04128	-0.48687\\
0.05174	-0.43193\\
0.06628	-0.38385\\
0.07791	-0.32675\\
0.09244	-0.25138\\
0.10581	-0.20322\\
0.12442	-0.16016\\
0.13779	-0.13753\\
0.16047	-0.09547\\
0.17151	-0.06731\\
0.19709	-0.01015\\
0.2064	-0.0016\\
0.23605	-0.00143\\
0.24593	-0.00193\\
0.27209	0.01024\\
0.2843	0.02124\\
0.31744	0.03957\\
0.32267	0.03823\\
0.36163	0.01293\\
0.37209	0.01474\\
0.40581	0.03636\\
0.41163	0.03805\\
0.45233	0.03033\\
0.46337	0.02768\\
0.50407	0.02723\\
0.50988	0.02864\\
0.55407	0.03887\\
0.56163	0.03772\\
0.60523	0.02104\\
0.61221	0.01921\\
0.66221	0.01204\\
0.66919	0.01048\\
0.72151	-0.00875\\
0.72209	-0.00901\\
0.77616	-0.03248\\
0.77907	-0.03348\\
0.83314	-0.05326\\
0.83837	-0.05619\\
0.88605	-0.08534\\
0.88721	-0.08582\\
0.94128	-0.10955\\
0.94186	-0.11017\\
0.97674	-0.18918\\
0.97733	-0.19101\\
1.00000	-0.19400\\
};
\addplot [color=blue, line width=1.3pt, mark=square, mark options={solid, blue}]
  table[row sep=crcr]{%
-4.57990949428877e-20	0.32724237\\
0.00139307630005401	1.4766173\\
0.00497626636451105	2.0085428\\
0.00970895471656558	2.1478341\\
0.0150663810641994	1.9225931\\
0.0263588828865331	1.8068154\\
0.0382461606009063	1.5619816\\
0.0620389332725132	1.4097658\\
0.0863171952784596	1.1902654\\
0.13678761396123	1.0187736\\
0.187630992234537	0.88646078\\
0.23860228657713	0.77190268\\
0.289624388138373	0.68010503\\
0.340605563180073	0.59895909\\
0.391568681470503	0.52744544\\
0.442500937830766	0.46281978\\
0.493398354450621	0.4038811\\
0.544219378021628	0.34932581\\
0.595003700375795	0.29823348\\
0.645751571987372	0.24955252\\
0.696463674082834	0.2023723\\
0.747139137169198	0.15470074\\
0.79777790426021	0.10703144\\
0.848377716186936	0.048402019\\
0.898934953143353	6.5039632e-05\\
0.924216456524101	-0.054763954\\
0.949484926671784	-0.08011467\\
0.974739587626676	-0.23783754\\
};
\addlegendentry{RELAX}

\addplot [color=blue, line width=1.3pt, mark=square, mark options={solid, blue}, forget plot]
  table[row sep=crcr]{%
-4.57990949428877e-20	0.32724237\\
0.00139078451033227	-0.65749151\\
0.00497199352608279	-0.99329668\\
0.00970306388271443	-0.98416418\\
0.0150591283867115	-0.85757977\\
0.0263494907835121	-0.70807093\\
0.0382350648782204	-0.52174962\\
0.0620252922726972	-0.3784197\\
0.0863016539364379	-0.22847083\\
0.13676940947696	-0.13274331\\
0.187611182259658	-0.043790273\\
0.238581632611066	-0.0081070475\\
0.289603447032044	0.016449865\\
0.340584751107735	0.027705457\\
0.391548362701881	0.032309588\\
0.442481401146082	0.031845741\\
0.493379843583215	0.028055478\\
0.544202100986024	0.02175327\\
0.594987778150088	0.013545571\\
0.645737210794862	0.0036245985\\
0.69645099746172	-0.0080033978\\
0.747128295802251	-0.022228058\\
0.797768970005893	-0.038695075\\
0.848370845564952	-0.064287394\\
0.898930207728643	-0.089864619\\
0.924212849418534	-0.12611163\\
0.949482533932792	-0.12673144\\
0.974738362879276	-0.25809342\\
1.00000000707044	-0.40402722\\
};
\addplot [color=green, line width=1.3pt, mark=square, mark options={solid, green}]
  table[row sep=crcr]{%
-4.57990949428877e-20	0.31768009\\
0.00139307630005401	1.4631367\\
0.00497626636451105	1.9946783\\
0.00970895471656558	2.1347613\\
0.0150663810641994	1.9115076\\
0.0263588828865331	1.7969128\\
0.0382461606009063	1.553604\\
0.0620389332725132	1.4024343\\
0.0863171952784596	1.1839283\\
0.13678761396123	1.0130969\\
0.187630992234537	0.88129264\\
0.23860228657713	0.76707059\\
0.289624388138373	0.67551732\\
0.340605563180073	0.59455889\\
0.391568681470503	0.52319407\\
0.442500937830766	0.45869005\\
0.493398354450621	0.39985296\\
0.544219378021628	0.34538502\\
0.595003700375795	0.29436949\\
0.645751571987372	0.24576052\\
0.696463674082834	0.19864367\\
0.747139137169198	0.151059\\
0.79777790426021	0.10339844\\
0.848377716186936	0.045092434\\
0.898934953143353	-0.004090209\\
0.924216456524101	-0.054590747\\
0.949484926671784	-0.097694851\\
0.974739587626676	-0.20506696\\
};
\addlegendentry{FLAT}

\addplot [color=green, line width=1.3pt, mark=square, mark options={solid, green}, forget plot]
  table[row sep=crcr]{%
-4.57990949428877e-20	0.31768009\\
0.00139078451033227	-0.66166568\\
0.00497199352608279	-0.99400717\\
0.00970306388271443	-0.98345959\\
0.0150591283867115	-0.85637003\\
0.0263494907835121	-0.7067728\\
0.0382350648782204	-0.52071363\\
0.0620252922726972	-0.37762114\\
0.0863016539364379	-0.22815569\\
0.13676940947696	-0.13282581\\
0.187611182259658	-0.044241346\\
0.238581632611066	-0.0088171456\\
0.289603447032044	0.015522179\\
0.340584751107735	0.026594771\\
0.391548362701881	0.031037509\\
0.442481401146082	0.030427868\\
0.493379843583215	0.026502689\\
0.544202100986024	0.020073457\\
0.594987778150088	0.011743786\\
0.645737210794862	0.0017060748\\
0.69645099746172	-0.010041384\\
0.747128295802251	-0.024362689\\
0.797768970005893	-0.041000787\\
0.848370845564952	-0.066479519\\
0.898930207728643	-0.093025036\\
0.924212849418534	-0.12506399\\
0.949482533932792	-0.14348908\\
0.974738362879276	-0.22481401\\
1.00000000707044	-0.29625547\\
};
\end{axis}

\node[align=center,font=\large, yshift=2em] (title) 
    at (current bounding box.north)
    {$\alpha$ = 6.75 ${^\circ}$, $\Lambda$ = 6\\ $C_p$ distribution, y/2c = 0.350};
\end{tikzpicture}%
}    
}
\\
\resizebox{0.5\linewidth}{!}{
\centerline{
%
%
\begin{tikzpicture}

\begin{axis}[%
width=4.602in,
height=3.506in,
at={(0.772in,0.473in)},
scale only axis,
xmin=0,
xmax=1,
xlabel style={font=\color{white!15!black}},
xlabel={\textbf{$x/c$}},
ymin=-1,
ymax=2.5,
ylabel style={font=\color{white!15!black}},
ylabel={\textbf{$Cp$}},
axis background/.style={fill=white},
xmajorgrids,
ymajorgrids,
legend style={legend cell align=left, align=left, draw=white!15!black}
]
\addplot [color=red, line width=1.3pt, only marks, mark=asterisk, mark options={solid, red}]
  table[row sep=crcr]{%
0.00245	1.58515\\
0.00281	1.62926\\
0.00476	1.86212\\
0.00639	2.02747\\
0.00829	2.13608\\
0.00935	2.15228\\
0.01124	2.13792\\
0.01287	2.10472\\
0.01595	2.03301\\
0.01697	2.01555\\
0.02049	1.97691\\
0.02125	1.95824\\
0.02597	1.83616\\
0.02868	1.78075\\
0.04362	1.58421\\
0.04452	1.57268\\
0.06656	1.30465\\
0.06856	1.28731\\
0.09419	1.17485\\
0.09675	1.16836\\
0.12593	1.07908\\
0.1267	1.07615\\
0.15649	0.94795\\
0.15782	0.94232\\
0.19293	0.82434\\
0.19482	0.82037\\
0.23182	0.743\\
0.2323	0.7411\\
0.27226	0.57432\\
0.27293	0.57348\\
0.30692	0.62852\\
0.3123	0.6328\\
0.35157	0.54654\\
0.35694	0.53163\\
0.39858	0.47258\\
0.39923	0.47207\\
0.4403	0.43656\\
0.44623	0.43032\\
0.49611	0.36783\\
0.49616	0.36777\\
0.54374	0.30819\\
0.54547	0.30613\\
0.59541	0.25754\\
0.59544	0.25752\\
0.64417	0.22985\\
0.65301	0.22168\\
0.70294	0.14724\\
0.70351	0.14652\\
0.75932	0.09063\\
0.75934	0.09061\\
0.81043	0.02178\\
0.81338	0.01913\\
0.87036	-0.01835\\
0.8739	-0.02266\\
0.93146	-0.13897\\
0.93263	-0.1417\\
0.96259	-0.20123\\
};
\addlegendentry{EXPERIMENTAL}

\addplot [color=red, line width=1.3pt, only marks, mark=asterisk, mark options={solid, red}, forget plot]
  table[row sep=crcr]{%
0.00245	-0.83433\\
0.00281	-0.88726\\
0.00476	-1.06019\\
0.00639	-1.07688\\
0.00829	-1.04832\\
0.00935	-1.0244\\
0.01124	-0.98889\\
0.01287	-0.97175\\
0.01595	-0.94135\\
0.01697	-0.92858\\
0.02049	-0.86882\\
0.02125	-0.84925\\
0.02597	-0.72742\\
0.02868	-0.69003\\
0.04362	-0.5343\\
0.04452	-0.52625\\
0.06656	-0.36795\\
0.06856	-0.35735\\
0.09419	-0.25511\\
0.09675	-0.2467\\
0.12593	-0.16605\\
0.1267	-0.16446\\
0.15649	-0.10968\\
0.15782	-0.10654\\
0.19293	-0.02658\\
0.19482	-0.02437\\
0.23182	-0.0121\\
0.2323	-0.01183\\
0.27226	0.02668\\
0.27293	0.027\\
0.30692	0.01469\\
0.3123	0.01283\\
0.35157	0.00861\\
0.35694	0.00772\\
0.39858	0.0055\\
0.39923	0.00577\\
0.4403	0.03212\\
0.44623	0.03277\\
0.49611	0.02601\\
0.49616	0.02602\\
0.54374	0.04338\\
0.54547	0.0437\\
0.59541	0.03463\\
0.59544	0.03463\\
0.64417	0.02557\\
0.65301	0.02159\\
0.70294	-0.00412\\
0.70351	-0.00432\\
0.75932	-0.01934\\
0.75934	-0.01935\\
0.81043	-0.06407\\
0.81338	-0.06636\\
0.87036	-0.09694\\
0.8739	-0.09868\\
0.93146	-0.13873\\
0.93263	-0.14009\\
0.96259	-0.20123\\
1.00000 -0.24000\\
};
\addplot [color=blue, line width=1.3pt, mark=square, mark options={solid, blue}]
  table[row sep=crcr]{%
0	0.13768013\\
0.00139307630005401	1.2285708\\
0.00497626636451105	1.7693139\\
0.00970895471656558	1.9329975\\
0.0150663810641994	1.749568\\
0.0263588828865331	1.6582673\\
0.0382461606009063	1.4443314\\
0.0620389332725132	1.3103012\\
0.0863171952784596	1.1103091\\
0.13676673355032	0.95224512\\
0.187589038468952	0.82977641\\
0.238539266902111	0.72185302\\
0.289540357954421	0.63514912\\
0.34051364552262	0.55816585\\
0.391468943681125	0.49023032\\
0.442393396046121	0.42880508\\
0.493283003429768	0.3727943\\
0.544118498322231	0.32098362\\
0.594917132407805	0.27253556\\
0.645679304649632	0.22637682\\
0.696405729297293	0.18192984\\
0.747095745872333	0.13670269\\
0.797749000230727	0.091873407\\
0.848363293022724	0.035945732\\
0.898934953143353	-0.0078868121\\
0.924216456524101	-0.059915166\\
0.949484926671784	-0.08829698\\
0.974739587626676	-0.2327902\\
};
\addlegendentry{RELAX}

\addplot [color=blue, line width=1.3pt, mark=square, mark options={solid, blue}, forget plot]
  table[row sep=crcr]{%
0	0.13768013\\
0.00139078451033227	-0.74634022\\
0.00497199352608279	-1.0070753\\
0.00970306388271443	-0.96180469\\
0.0150591283867115	-0.81820428\\
0.0263494907835121	-0.6598534\\
0.0382350648782204	-0.47291726\\
0.0620252922726972	-0.33079049\\
0.0863016539364379	-0.18510312\\
0.136748519476892	-0.094267577\\
0.187569233232222	-0.0090046404\\
0.23851861664087	0.023520786\\
0.289519408965734	0.045551863\\
0.34049284060909	0.054614462\\
0.391448614321994	0.057291396\\
0.442373867181741	0.055080235\\
0.49326450064693	0.049672686\\
0.544101186303544	0.041843485\\
0.594901242355705	0.032159004\\
0.645664943235695	0.02080497\\
0.696393087921635	0.0076161916\\
0.747084843995798	-0.0080490718\\
0.797740061663123	-0.026317388\\
0.848356349116228	-0.053663727\\
0.898930207728643	-0.082767606\\
0.924212849418534	-0.11953142\\
0.949482533932792	-0.12436555\\
0.974738362879276	-0.24855115\\
1.00000000707044	-0.38190392\\
};
\addplot [color=green, line width=1.3pt, mark=square, mark options={solid, green}]
  table[row sep=crcr]{%
0	0.13053292\\
0.00139307630005401	1.2181437\\
0.00497626636451105	1.7583808\\
0.00970895471656558	1.9225824\\
0.0150663810641994	1.7406144\\
0.0263588828865331	1.6502986\\
0.0382461606009063	1.4374704\\
0.0620389332725132	1.3042215\\
0.0863171952784596	1.1050242\\
0.13676673355032	0.94742876\\
0.187589038468952	0.82532978\\
0.238539266902111	0.71763211\\
0.289540357954421	0.63107777\\
0.34051364552262	0.55419433\\
0.391468943681125	0.48632333\\
0.442393396046121	0.42493635\\
0.493283003429768	0.36894339\\
0.544118498322231	0.31713501\\
0.594917132407805	0.26867649\\
0.645679304649632	0.22250132\\
0.696405729297293	0.17802379\\
0.747095745872333	0.13279212\\
0.797749000230727	0.087888554\\
0.848363293022724	0.03215804\\
0.898934953143353	-0.01245138\\
0.924216456524101	-0.061033193\\
0.949484926671784	-0.1036009\\
0.974739587626676	-0.20765002\\
};
\addlegendentry{FLAT}

\addplot [color=green, line width=1.3pt, mark=square, mark options={solid, green}, forget plot]
  table[row sep=crcr]{%
0	0.13053292\\
0.00139078451033227	-0.74921519\\
0.00497199352608279	-1.0073361\\
0.00970306388271443	-0.96105832\\
0.0150591283867115	-0.81718624\\
0.0263494907835121	-0.65881014\\
0.0382350648782204	-0.47219166\\
0.0620252922726972	-0.33031413\\
0.0863016539364379	-0.18506745\\
0.136748519476892	-0.094621725\\
0.187569233232222	-0.0097197145\\
0.23851861664087	0.022538759\\
0.289519408965734	0.044334605\\
0.34049284060909	0.053188197\\
0.391448614321994	0.05567066\\
0.442373867181741	0.053274341\\
0.49326450064693	0.047686789\\
0.544101186303544	0.039680358\\
0.594901242355705	0.029819058\\
0.645664943235695	0.018289946\\
0.696393087921635	0.0049196584\\
0.747084843995798	-0.010905827\\
0.797740061663123	-0.029393932\\
0.848356349116228	-0.056721248\\
0.898930207728643	-0.08665587\\
0.924212849418534	-0.12006037\\
0.949482533932792	-0.1391838\\
0.974738362879276	-0.22308336\\
1.00000000707044	-0.29653022\\
};
\end{axis}

\node[align=center,font=\large, yshift=2em] (title) 
    at (current bounding box.north)
    {$\alpha$ = 6.75 ${^\circ}$, $\Lambda$ = 6\\ $C_p$ distribution, y/2c = 0.610};
\end{tikzpicture}%
}    
}
&
\resizebox{0.5\linewidth}{!}{
\centerline{
%
%
\begin{tikzpicture}

\begin{axis}[%
width=4.521in,
height=3.566in,
at={(0.758in,0.481in)},
scale only axis,
xmin=0,
xmax=1,
xlabel style={font=\color{white!15!black}},
xlabel={\textbf{$x/c$}},
ymin=-1,
ymax=2,
ylabel style={font=\color{white!15!black}},
ylabel={\textbf{$Cp$}},
axis background/.style={fill=white},
xmajorgrids,
ymajorgrids,
legend style={legend cell align=left, align=left, draw=white!15!black}
]
\addplot [color=red, line width=1.3pt, only marks, mark=asterisk, mark options={solid, red}]
  table[row sep=crcr]{%
6.25e-08	-0.16184\\
0.001969939456	0.92764\\
0.004829277049	1.40113\\
0.008393308225	1.44467\\
0.012195668356	1.60413\\
0.015125556196	1.46407\\
0.020869847296	1.46653\\
0.026514911556	1.33614\\
0.043245282025	1.23963\\
0.0681262201	1.01995\\
0.096943804164	0.92587\\
0.127894925376	0.8197\\
0.160827467089	0.74251\\
0.197489582404	0.63634\\
0.234621078129	0.59296\\
0.278506341169	0.48437\\
0.317756562601	0.46996\\
0.3643691769	0.39518\\
0.408296718361	0.36145\\
0.455485960609	0.31806\\
0.506945424001	0.29157\\
0.557634069001	0.24335\\
0.612583286329	0.2072\\
0.6701387044	0.18072\\
0.727254195264	0.13491\\
0.787720151296	0.08427\\
0.842107216896	0.05053\\
0.899365825801	0.00229\\
0.955044925696	-0.06527\\
0.996073861225	-0.14493\\
};
\addlegendentry{EXPERIMENTAL}

\addplot [color=red, line width=1.3pt, only marks, mark=asterisk, mark options={solid, red}, forget plot]
  table[row sep=crcr]{%
6.25e-08	-0.14744\\
0.001969939456	-0.93864\\
0.004829277049	-0.98485\\
0.008393308225	-0.90121\\
0.012195668356	-0.83394\\
0.015125556196	-0.71536\\
0.020869847296	-0.63194\\
0.026514911556	-0.51907\\
0.043245282025	-0.35378\\
0.0681262201	-0.12972\\
0.096943804164	-0.1088\\
0.127894925376	0.00901\\
0.160827467089	0.00766\\
0.197489582404	0.10849\\
0.234621078129	0.10513\\
0.278506341169	0.12336\\
0.317756562601	0.13417\\
0.3643691769	0.12561\\
0.408296718361	0.10786\\
0.455485960609	0.11038\\
0.506945424001	0.09839\\
0.557634069001	0.11534\\
0.612583286329	0.09115\\
0.6701387044	0.08362\\
0.727254195264	0.06753\\
0.787720151296	0.0421\\
0.842107216896	0.00272\\
0.899365825801	-0.02092\\
0.955044925696	-0.0815\\
0.996073861225	-0.15489\\
};
\addplot [color=blue, line width=1.3pt, mark=square, mark options={solid, blue}]
  table[row sep=crcr]{%
-6.45314452950959e-19	-0.35982135\\
0.00138090587516383	0.53247041\\
0.00493680989673767	1.0730702\\
0.00963919309516276	1.3013992\\
0.01496575456924	1.2246664\\
0.0149657548043148	1.2246664\\
0.0261994378325879	1.2048931\\
0.0380285791218925	1.0721411\\
0.0380285830941663	1.0721411\\
0.0619268224690066	0.97767198\\
0.0863171952784596	0.83138335\\
0.136726115579987	0.70528191\\
0.187508727376033	0.61380017\\
0.238432733086602	0.5270102\\
0.289407820996959	0.45831379\\
0.340382299032716	0.39733273\\
0.39133885665323	0.34436551\\
0.44226558279383	0.29713613\\
0.493157516502823	0.25528187\\
0.544008356853286	0.21615575\\
0.594822388936075	0.18543766\\
0.645600124535529	0.13960861\\
0.696342087865318	0.16051839\\
0.721699600454973	0.1626232\\
0.747048091066029	0.068432547\\
0.772387327394397	0.070001528\\
0.797717247730493	0.044798519\\
0.823037364302751	0.02561919\\
0.848347399524606	0.0033317774\\
0.873646720742548	-0.020810096\\
0.898934953143353	-0.046621352\\
0.92421651693097	-0.082299478\\
0.949484988606639	-0.10520927\\
0.974739649561532	-0.21435177\\
};
\addlegendentry{RELAX}

\addplot [color=blue, line width=1.3pt, mark=square, mark options={solid, blue}, forget plot]
  table[row sep=crcr]{%
-6.45314452950959e-19	-0.35982135\\
0.00137862392891182	-0.93936974\\
0.00493255329796678	-0.99145335\\
0.00963332197481834	-0.84601921\\
0.0149585245579034	-0.6684255\\
0.0261900729166228	-0.48271376\\
0.0380175136195508	-0.30917069\\
0.0380175140897004	-0.30917069\\
0.0619131907703767	-0.13986585\\
0.0863016539364379	-0.014769257\\
0.136707913005504	0.026780572\\
0.187488917813536	0.1022156\\
0.238412069097884	0.11871053\\
0.289386883629768	0.13045447\\
0.340361470860317	0.13095903\\
0.391318562745218	0.12727925\\
0.442246063730287	0.12024167\\
0.493139012244865	0.11110329\\
0.543991057408873	0.10050868\\
0.594806498346702	0.086875446\\
0.645585732832652	0.075350434\\
0.696329373512585	0.041556198\\
0.721687800938375	0.013659894\\
0.747037184457529	0.04922669\\
0.772377455515161	0.031622987\\
0.797708310793818	0.027225792\\
0.8230294151687	0.017380564\\
0.848340492163718	0.0066071874\\
0.873640963591637	-0.0076635755\\
0.898930207728643	-0.02431808\\
0.924212909002642	-0.053363994\\
0.949482593516899	-0.071247697\\
0.974738422463383	-0.19689131\\
1.00000000707044	-0.32422554\\
};

\addplot [color=green, line width=1.3pt, mark=square, mark options={solid, green}]
  table[row sep=crcr]{%
-6.45314452950959e-19	-0.36286998\\
0.00138090587516383	0.52766877\\
0.00493680989673767	1.0679177\\
0.00963919309516276	1.296478\\
0.01496575456924	1.2204391\\
0.0149657548043148	1.220439\\
0.0261994378325879	1.2011492\\
0.0380285791218925	1.0689626\\
0.0380285830941663	1.0689627\\
0.0619268224690066	0.97496039\\
0.0863171952784596	0.82908958\\
0.136726115579987	0.70316035\\
0.187508727376033	0.61186498\\
0.238432733086602	0.52515465\\
0.289407820996959	0.45649359\\
0.340382299032716	0.3955076\\
0.39133885665323	0.3425017\\
0.44226558279383	0.29520056\\
0.493157516502823	0.25323528\\
0.544008356853286	0.21396503\\
0.594822388936075	0.18301585\\
0.645600124535529	0.13702525\\
0.696342087865318	0.1571994\\
0.721699600454973	0.15857646\\
0.747048091066029	0.065029189\\
0.772387327394397	0.065995172\\
0.797717247730493	0.040451739\\
0.823037364302751	0.020802563\\
0.848347399524606	-0.0020594257\\
0.873646720742548	-0.026639607\\
0.898934953143353	-0.053788852\\
0.92421651693097	-0.08784841\\
0.949484988606639	-0.120584\\
0.974739649561532	-0.20557018\\
};
\addlegendentry{FLAT}

\addplot [color=green, line width=1.3pt, mark=square, mark options={solid, green}]
  table[row sep=crcr]{%
-6.45314452950959e-19	-0.36286998\\
0.00137862392891182	-0.94022632\\
0.00493255329796678	-0.99104184\\
0.00963332197481834	-0.84514308\\
0.0149585245579034	-0.6675061\\
0.0261900729166228	-0.48183882\\
0.0380175136195508	-0.30851617\\
0.0380175140897004	-0.30851617\\
0.0619131907703767	-0.13939992\\
0.0863016539364379	-0.014567128\\
0.136707913005504	0.026738096\\
0.187488917813536	0.10193802\\
0.238412069097884	0.11824618\\
0.289386883629768	0.12981057\\
0.340361470860317	0.13013618\\
0.391318562745218	0.12626623\\
0.442246063730287	0.11901854\\
0.493139012244865	0.10963695\\
0.543991057408873	0.098761603\\
0.594806498346702	0.084763125\\
0.645585732832652	0.072885536\\
0.696329373512585	0.038289141\\
0.721687800938375	0.0092966808\\
0.747037184457529	0.045366742\\
0.772377455515161	0.0269947\\
0.797708310793818	0.022072636\\
0.8230294151687	0.011550031\\
0.848340492163718	-3.3464643e-05\\
0.873640963591637	-0.014996725\\
0.898930207728643	-0.033293545\\
0.924212909002642	-0.060900122\\
0.949482593516899	-0.089062326\\
0.974738422463383	-0.18886508\\
1.00000000707044	-0.27772605\\
};
\end{axis}

\node[align=center,font=\large, yshift=2em] (title) 
    at (current bounding box.north)
    {$\alpha$ = 6.75 ${^\circ}$, $\Lambda$ = 6\\ $C_p$ distribution, y/2c = 0.900};
\end{tikzpicture}%
}    
}
\end{tabular}  	
\caption{Comparative analysis between results obtained using $\Pi$-BEM and experimental data 
for a rectangular planform wing. The plots refer to four planar sections normal to the $y$-axis 
located along the semi-wingspan direction at $y/s$ = 0.050, 0.350, 0.610, 0.900.\label{fig:RectangularWing}}
\end{figure}
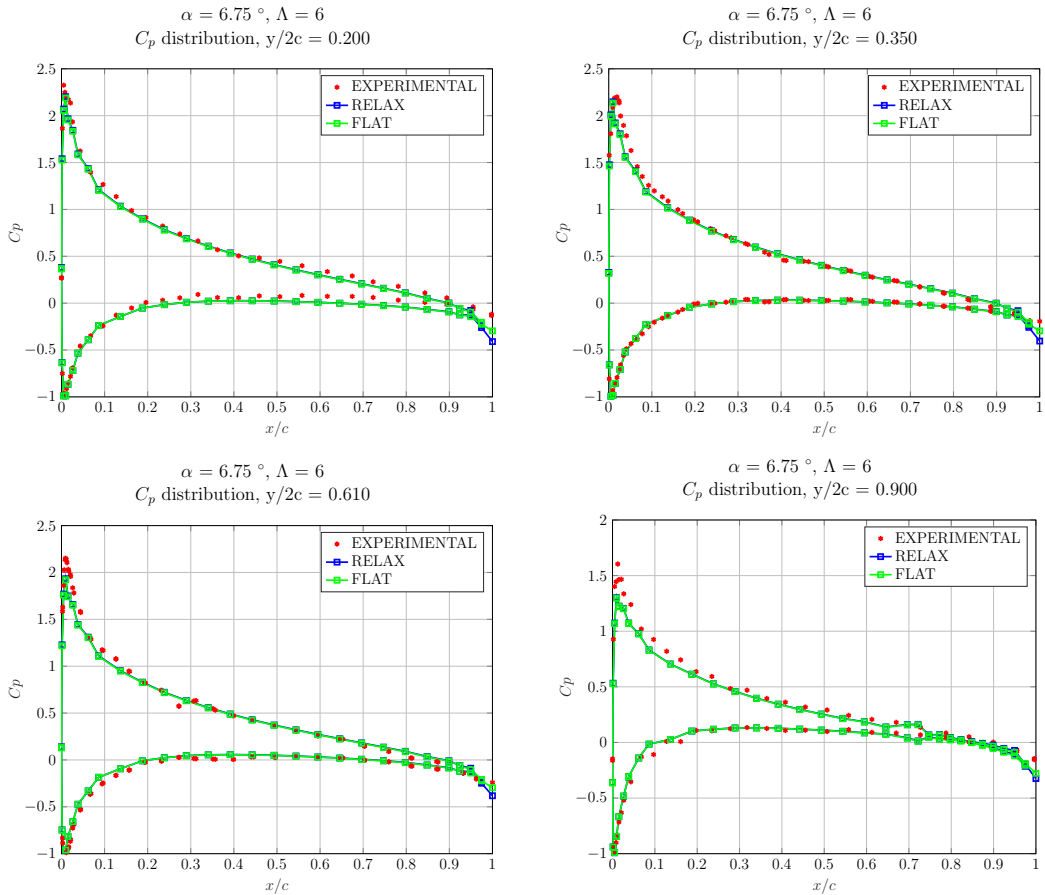

Figure \ref{fig:RectangularWing}, presents a comparison between experimental and numerical values
of $-c_p$ on four cross-sections perpendicular to the wingspan direction. The wing sections are located
at wingspan coordinates $y/c$ = 0.200, 0.350, 0.610, 0.900.
In each section, the numerical values of windward and leeward pressure coefficients obtained imposing a
flat wake (blue lines) and with relaxed wake (green lines) are plotted against their experimental counterparts (red dots).

By a qualitative standpoint, the comparisons appear satisfactory, as the model appears able to reproduce
with accuracy the main features of the experimental curves. By a more quantitative point of view, the
peaks location is well captured and their height, especially on the windward side stagnation point, appears accurately reproduced by the numerical results.
Also the height of the suction side peaks seems reasonably accurate, although slightly underestimated. This is more evident for the peaks on the
sections closer to the tip.
Such a discrepancy is likely due to the potential flow model inability to reproduce viscous effects
associated with the increasingly close tip vortex. In addition, in every plot in Figure \ref{fig:RectangularWing} the numerical
solution appears to be consistently overestimating the experimental pressure at the trailing edge. This should be in part
associated to the sharp trailing edge considered in the potential flow simulations which --- as opposed to the finite thickness one used
in the experimental campaign --- leads to higher pressure recovery. Finally, we point out that the results obtained with flat wake are extremely
close to the ones obtained including the wake relaxation part of the algorithm. The only visible differences can be in fact observed at the
trailing edge, where the relaxed wake solution consistently exhibits higher pressure values than its flat wake counterpart.
On one hand, such a marginal difference might not seem enough to justify
the additional computational cost associated with the wake velocity computation and wake relaxation part of the resolution algorithm.
On the other hand, there surely are different applications in which the wake geometry has more influence on the computed body pressure,
with respect to the isolated wing considered in the present numerical tests. This is definitely the case for propellers, in which the
helical shape of the wake is typically closer to the propeller blades than it is in the example here considered.

\subsection{Swept wing}

\begin{figure}[htb!]
\begin{center}
\includegraphics[width=\textwidth]{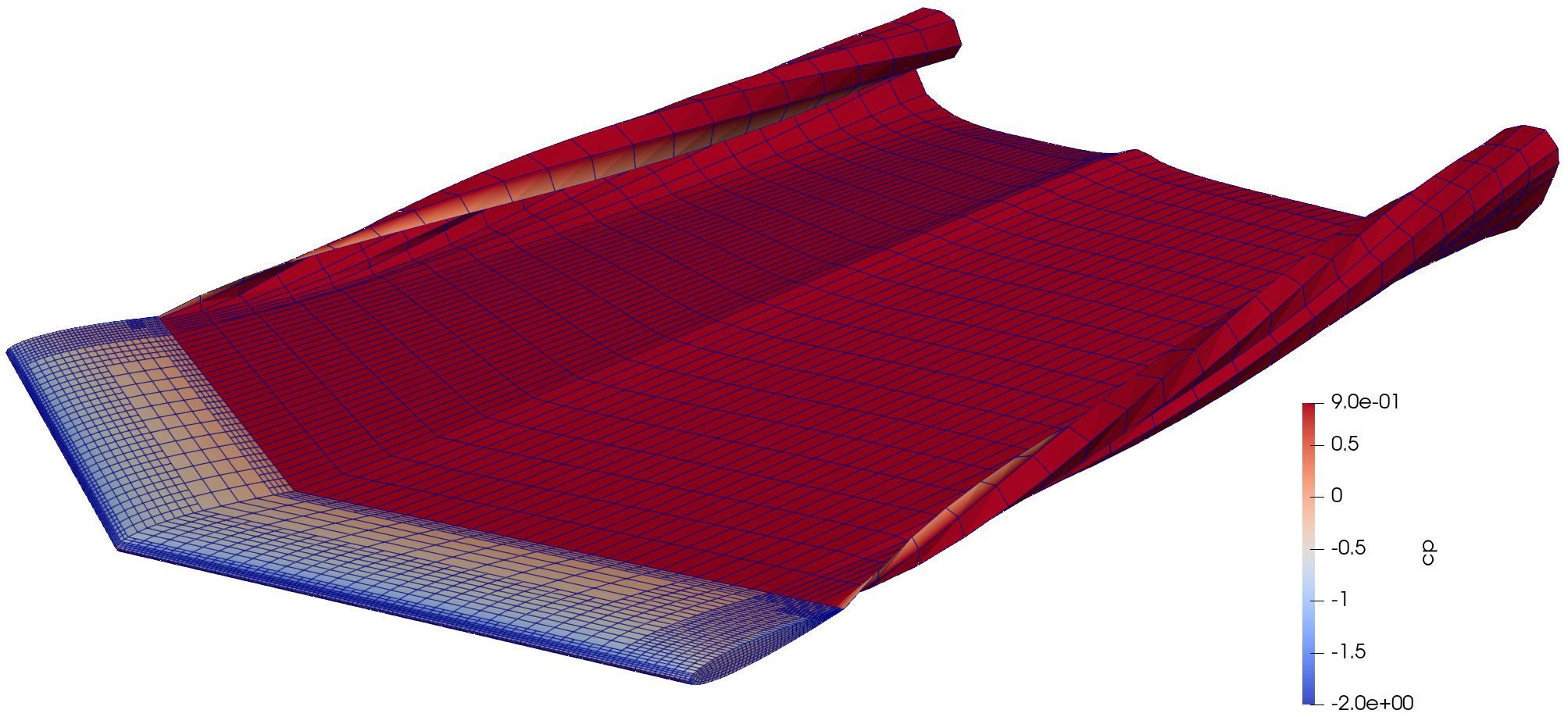}
\end{center}
\caption{A view of the grid obtained on the geometry of the NACA 0012 wing with aspect ratio $\Lambda=5.47$ and sweep angle $\beta=20{^\circ}$, and on
         the fully developed wake following the body. The surfaces are colored according to contours of pressure coefficient $c_p$.  
         \label{fig:VWing}}
\end{figure}

The experimental setup corresponding to the semi-wing with sweep angle $\beta=20^\circ$ previously described, is also in this case reproduced in the simulations
by considering the entire extent of a V-shaped wing, as illustrated in Figure \ref{fig:VWing}. In the experimental setup \cite{YipShubertNACA0012},
the swept wing was obtained first by rotating the NACA 0012 rectangular wing previously considered, and then adding a fairing at the tip. This
resulted in a wing having semi-wing span $s=2.91c$. The geometric CAD model, designed using FreeCad, features again flat cap
surfaces on the wing tips. The initial mesh made of approximately 10 cells is refined onto the surface of the CAD model imposing, also in this case, maximum
cell aspect ratio 4.5, 2 uniform refinement cycles and 1 adaptive refinement cycle based on CAD curvature. This results in a computational grid
featuring approximately 4000 quadrilateral cells. Bi-linear finite elements are used for this test too, which result in  
a resolution system of 4520 degrees of freedom. Figure \ref{fig:VWing} depicts the computational grid on the wing and on the fully relaxed wake. The
body surface is colored according to contours of pressure coefficient $c_p$.

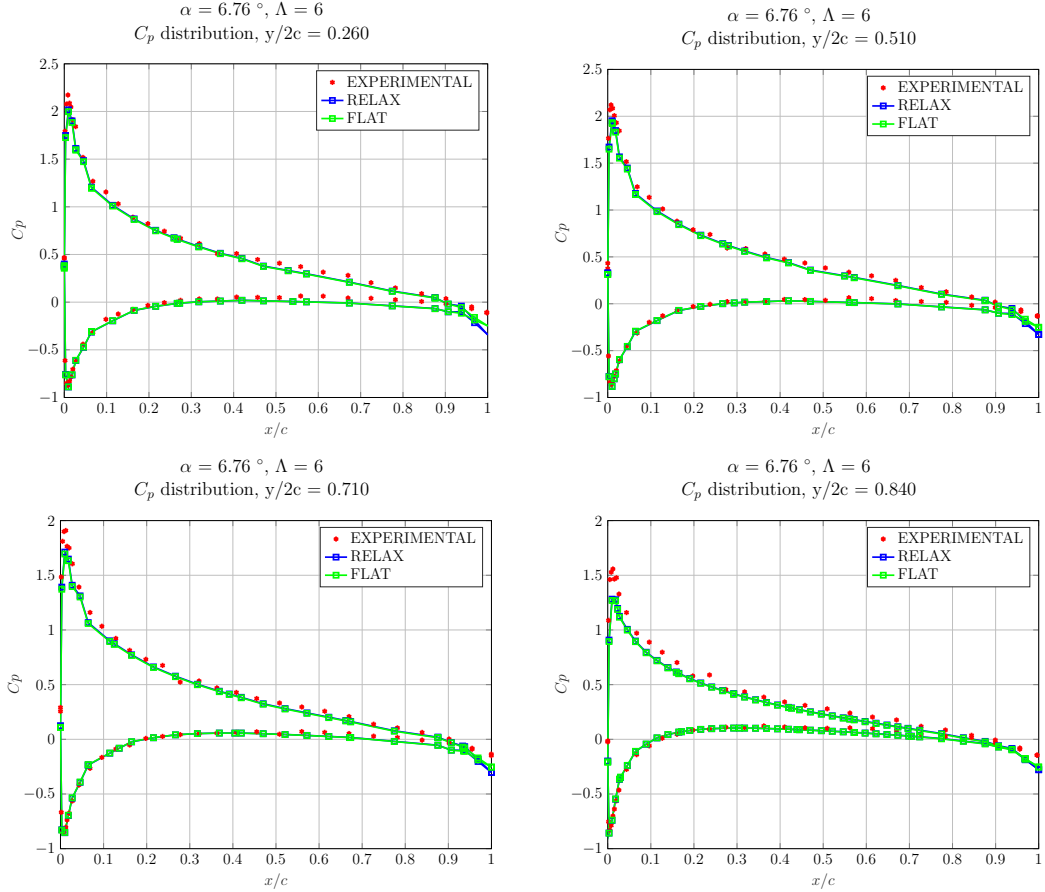
\begin{figure}[htb!]
\begin{tabular}{c c}
\resizebox{0.5\linewidth}{!}{
\centerline{
%
%
\begin{tikzpicture}

\begin{axis}[%
width=4.521in,
height=3.566in,
at={(0.758in,0.481in)},
scale only axis,
xmin=0,
xmax=1,
xlabel style={font=\color{white!15!black}},
xlabel={\textbf{$x/c$}},
ymin=-1,
ymax=2.5,
ylabel style={font=\color{white!15!black}},
ylabel={\textbf{$Cp$}},
axis background/.style={fill=white},
xmajorgrids,
ymajorgrids,
legend style={legend cell align=left, align=left, draw=white!15!black}
]
\addplot [color=red, line width=1.3pt, only marks, mark=asterisk, mark options={solid, red}]
  table[row sep=crcr]{%
1.684804e-06	0.45805\\
0.001845303849	1.79122\\
0.005224253841	2.07791\\
0.008556435001	2.17251\\
0.012355656336	2.08676\\
0.015656265625	2.04618\\
0.020188147225	1.91299\\
0.026774122384	1.84143\\
0.044256378384	1.51787\\
0.067571922916	1.26783\\
0.0982509025	1.15539\\
0.127519695801	1.03125\\
0.162010275025	0.89046\\
0.197453143449	0.82331\\
0.236410860841	0.7443\\
0.274594864324	0.67009\\
0.318883831204	0.61483\\
0.361641065956	0.50741\\
0.406953409041	0.50918\\
0.456589409796	0.44684\\
0.508220706816	0.40825\\
0.558267491929	0.37208\\
0.611605330704	0.31454\\
0.670186185201	0.2807\\
0.724582298176	0.20421\\
0.7820926096	0.14907\\
0.8439361956	0.08917\\
0.901495179841	0.03644\\
0.962155001025	-0.04952\\
0.997887117249	-0.10679\\
};
\addlegendentry{EXPERIMENTAL}

\addplot [color=red, line width=1.3pt, only marks, mark=asterisk, mark options={solid, red}, forget plot]
  table[row sep=crcr]{%
1.684804e-06	0.46789\\
0.001845303849	-0.61298\\
0.005224253841	-0.85044\\
0.008556435001	-0.87185\\
0.012355656336	-0.82908\\
0.015656265625	-0.78461\\
0.020188147225	-0.70365\\
0.026774122384	-0.61994\\
0.044256378384	-0.44128\\
0.067571922916	-0.32002\\
0.0982509025	-0.17996\\
0.127519695801	-0.12456\\
0.162010275025	-0.07655\\
0.197453143449	-0.0342\\
0.236410860841	-0.00542\\
0.274594864324	0.01905\\
0.318883831204	0.03338\\
0.361641065956	0.03259\\
0.406953409041	0.05308\\
0.456589409796	0.05029\\
0.508220706816	0.0474\\
0.558267491929	0.06411\\
0.611605330704	0.06257\\
0.670186185201	0.04292\\
0.724582298176	0.04038\\
0.7820926096	0.02509\\
0.8439361956	0.00662\\
0.901495179841	-0.01764\\
0.962155001025	-0.06842\\
0.997887117249	-0.11208\\
};
\addplot [color=blue, line width=1.3pt, mark=square, mark options={solid, blue}]
  table[row sep=crcr]{%
-1.66707941895439e-08	0.38612875\\
3.06276164482866e-05	0.39978531\\
0.00305491535365421	1.7483197\\
0.00998726205971484	2.0171466\\
0.0182825148861865	1.9013003\\
0.0272149107868748	1.6098701\\
0.0272149110080986	1.6098701\\
0.0455060723811834	1.4848394\\
0.0643378476117348	1.2052388\\
0.0643378757171409	1.2052389\\
0.114390455924057	1.0155914\\
0.165034723029547	0.87327474\\
0.215884562923509	0.75384951\\
0.259229929849465	0.67613125\\
0.266818324376037	0.66252136\\
0.317739423419114	0.58355534\\
0.368650393076084	0.51269507\\
0.419534549479196	0.4592481\\
0.470384928956074	0.37861845\\
0.528789295311777	0.33226529\\
0.571930826644706	0.29806721\\
0.673335237845187	0.21061133\\
0.774599462495295	0.11599799\\
0.875715387449974	0.046672691\\
0.875715480765955	0.046672691\\
0.906810704171077	-0.020146282\\
0.937887603789149	-0.041169465\\
0.968944636181288	-0.18846777\\
1.00000009656696	-0.3415049\\
};
\addlegendentry{RELAX}

\addplot [color=blue, line width=1.3pt, mark=square, mark options={solid, blue}, forget plot]
  table[row sep=crcr]{%
-1.66707941895439e-08	0.38612875\\
3.85725956500438e-05	0.37167037\\
0.00305489170235146	-0.75831097\\
0.00998722453265299	-0.88952428\\
0.0182825320581843	-0.76344591\\
0.0272149029535996	-0.61269939\\
0.0455060500565578	-0.47437924\\
0.0643378741425521	-0.31168479\\
0.0643379021373464	-0.31168479\\
0.114390547656573	-0.19601046\\
0.165034965686852	-0.086726077\\
0.215884927780181	-0.043962728\\
0.266818825758526	-0.012851632\\
0.272541662382939	-0.010997047\\
0.317740036242815	0.0036620595\\
0.368651164956183	0.011481653\\
0.419535375912407	0.019266833\\
0.470385878589371	0.013901423\\
0.539984512059771	0.0066647399\\
0.571931453373099	0.0033784579\\
0.673335704958388	-0.0089545427\\
0.774599465338668	-0.038596615\\
0.875715475965536	-0.06686189\\
0.906810636971452	-0.10163757\\
0.937887588147395	-0.10211483\\
0.968944604933125	-0.2154846\\
1.00000009656696	-0.3415049\\
};
\addplot [color=green, line width=1.3pt, mark=square, mark options={solid, green}]
  table[row sep=crcr]{%
-1.66707941895439e-08	0.37000832\\
3.06276164482866e-05	0.38360226\\
0.00305491535365421	1.7259561\\
0.00998726205971484	1.9975324\\
0.0182825148861865	1.8847501\\
0.0272149107868748	1.5966308\\
0.0272149110080986	1.5966308\\
0.0455060723811834	1.473354\\
0.0643378476117348	1.196228\\
0.0643378757171409	1.196228\\
0.114390455924057	1.008082\\
0.165034723029547	0.86687958\\
0.215884562923509	0.74813926\\
0.259229929849465	0.67084074\\
0.266818324376037	0.65730441\\
0.317739423419114	0.57871628\\
0.368650393076084	0.5081526\\
0.419534549479196	0.45495257\\
0.470384928956074	0.37460777\\
0.528789295311777	0.32841071\\
0.571930826644706	0.29432869\\
0.673335237845187	0.20712097\\
0.774599462495295	0.1129146\\
0.875715387449974	0.043221656\\
0.875715480765955	0.043221653\\
0.906810704171077	-0.020018078\\
0.937887603789149	-0.056238893\\
0.968944636181288	-0.16025279\\
1.00000009656696	-0.24865077\\
};
\addlegendentry{FLAT}

\addplot [color=green, line width=1.3pt, mark=square, mark options={solid, green}, forget plot]
  table[row sep=crcr]{%
-1.66707941895439e-08	0.37000832\\
3.85725956500438e-05	0.35569471\\
0.00305489170235146	-0.76297277\\
0.00998722453265299	-0.88845462\\
0.0182825320581843	-0.76094484\\
0.0272149029535996	-0.60992217\\
0.0455060500565578	-0.471632\\
0.0643378741425521	-0.30945843\\
0.0643379021373464	-0.30945843\\
0.114390547656573	-0.19432476\\
0.165034965686852	-0.085574336\\
0.215884927780181	-0.043196745\\
0.266818825758526	-0.012399404\\
0.272541662382939	-0.010573797\\
0.317740036242815	0.0038569446\\
0.368651164956183	0.011442382\\
0.419535375912407	0.019056266\\
0.470385878589371	0.013443983\\
0.539984512059771	0.0060324967\\
0.571931453373099	0.0026669551\\
0.673335704958388	-0.010008734\\
0.774599465338668	-0.039821234\\
0.875715475965536	-0.068859175\\
0.906810636971452	-0.10042676\\
0.937887588147395	-0.11612267\\
0.968944604933125	-0.18668114\\
1.00000009656696	-0.24865077\\
};
\end{axis}

\node[align=center,font=\large, yshift=2em] (title) 
    at (current bounding box.north)
    {$\alpha$ = 6.76 ${^\circ}$, $\Lambda$ = 6\\ $C_p$ distribution, y/2c = 0.260};
\end{tikzpicture}%
}    
}
&
\resizebox{0.5\linewidth}{!}{
\centerline{
%
%
\begin{tikzpicture}

\begin{axis}[%
width=4.602in,
height=3.506in,
at={(0.772in,0.473in)},
scale only axis,
xmin=0,
xmax=1,
xlabel style={font=\color{white!15!black}},
xlabel={\textbf{$x/c$}},
ymin=-1,
ymax=2.5,
ylabel style={font=\color{white!15!black}},
ylabel={\textbf{$Cp$}},
axis background/.style={fill=white},
xmajorgrids,
ymajorgrids,
legend style={legend cell align=left, align=left, draw=white!15!black}
]
\addplot [color=red, line width=1.3pt, only marks, mark=asterisk, mark options={solid, red}]
  table[row sep=crcr]{%
1.545049e-06	0.4326\\
0.001732640625	1.76507\\
0.004873156864	2.06944\\
0.007709366809	2.12106\\
0.011338416324	2.0852\\
0.015387410116	2.00682\\
0.020052542449	1.93079\\
0.026842890244	1.84521\\
0.043223657409	1.51567\\
0.068537669209	1.24738\\
0.096266231824	1.13524\\
0.127452570025	1.01132\\
0.160671102244	0.88036\\
0.198209272849	0.78958\\
0.23736384	0.73904\\
0.2770864321	0.594\\
0.320357132001	0.58843\\
0.364724197776	0.53089\\
0.409725449604	0.47574\\
0.457320825025	0.43478\\
0.505064641041	0.38204\\
0.559514472049	0.33633\\
0.612171667396	0.29777\\
0.668165647396	0.24975\\
0.7276601809	0.17334\\
0.785494010961	0.12772\\
0.8444874816	0.07265\\
0.898978941025	0.01765\\
0.9586563921	-0.06575\\
0.996804556801	-0.12763\\
};
\addlegendentry{EXPERIMENTAL}

\addplot [color=red, line width=1.3pt, only marks, mark=asterisk, mark options={solid, red}, forget plot]
  table[row sep=crcr]{%
1.545049e-06	0.37473\\
0.001732640625	-0.55734\\
0.004873156864	-0.83535\\
0.007709366809	-0.86198\\
0.011338416324	-0.83956\\
0.015387410116	-0.77308\\
0.020052542449	-0.71284\\
0.026842890244	-0.61323\\
0.043223657409	-0.463\\
0.068537669209	-0.31316\\
0.096266231824	-0.19784\\
0.127452570025	-0.12843\\
0.160671102244	-0.07253\\
0.198209272849	-0.02938\\
0.23736384	-0.00609\\
0.2770864321	0.02329\\
0.320357132001	0.02541\\
0.364724197776	0.02562\\
0.409725449604	0.04552\\
0.457320825025	0.04467\\
0.505064641041	0.03649\\
0.559514472049	0.06481\\
0.612171667396	0.05445\\
0.668165647396	0.03448\\
0.7276601809	0.02436\\
0.785494010961	0.01512\\
0.8444874816	-0.01683\\
0.898978941025	-0.04038\\
0.9586563921	-0.08349\\
0.996804556801	-0.13307\\
};
\addplot [color=blue, line width=1.3pt, mark=square, mark options={solid, blue}]
  table[row sep=crcr]{%
4.40718485039682e-08	0.32800156\\
0.00305491401010027	1.6719398\\
0.00998723360626197	1.952278\\
0.0108129253219215	1.9418803\\
0.0182824495487888	1.8476698\\
0.0272148676794989	1.5669713\\
0.0455059818415336	1.4481984\\
0.0643378249230543	1.1752503\\
0.0643378252548901	1.1752503\\
0.114392465573242	0.99018306\\
0.165040847280624	0.85026622\\
0.215896824209716	0.73222864\\
0.26683867368624	0.64181346\\
0.280160357243089	0.62132794\\
0.317769800918135	0.56336588\\
0.368692731322828	0.49329588\\
0.419590687014326	0.43939087\\
0.470457034601346	0.35985237\\
0.549745000555364	0.2984052\\
0.571996430244576	0.28102967\\
0.673386550433174	0.1952574\\
0.774628846383967	0.10325631\\
0.87571532435139	0.036960937\\
0.875715326563628	0.036960937\\
0.906810625077606	-0.025223905\\
0.93788744274642	-0.049924925\\
0.968944452405112	-0.18716511\\
0.999999897902016	-0.32815316\\
};
\addlegendentry{RELAX}

\addplot [color=blue, line width=1.3pt, mark=square, mark options={solid, blue}, forget plot]
  table[row sep=crcr]{%
4.40718485039682e-08	0.32800156\\
0.00305491401248303	-0.77500534\\
0.00998719963173467	-0.88339645\\
0.0150907237466243	-0.80153036\\
0.0182824712960732	-0.75034559\\
0.0272148896237624	-0.59703642\\
0.0455059860594568	-0.45765257\\
0.0643377785147553	-0.29503062\\
0.114392386981218	-0.18032429\\
0.165040598944605	-0.071624622\\
0.215896467105962	-0.029627906\\
0.266838109200108	0.00078782113\\
0.293534440914652	0.0091374973\\
0.317769099637719	0.016765088\\
0.368692009998319	0.02353102\\
0.419589948586958	0.031887222\\
0.470456028878461	0.025210405\\
0.560629529825513	0.014054058\\
0.57199573567152	0.012693216\\
0.673386241974256	-0.0013172842\\
0.774628684343519	-0.033137634\\
0.875715298927918	-0.063973561\\
0.906810670507207	-0.10133822\\
0.93788741795018	-0.10363954\\
0.968944528051892	-0.2118565\\
0.999999897902016	-0.32815316\\
};
\addplot [color=green, line width=1.3pt, mark=square, mark options={solid, green}]
  table[row sep=crcr]{%
4.40718485039682e-08	0.31193304\\
0.00305491401010027	1.6493003\\
0.00998723360626197	1.93231\\
0.0108129253219215	1.9222174\\
0.0182824495487888	1.83077\\
0.0272148676794989	1.5534438\\
0.0455059818415336	1.436432\\
0.0643378249230543	1.1660188\\
0.0643378252548901	1.1660188\\
0.114392465573242	0.98245537\\
0.165040847280624	0.84365743\\
0.215896824209716	0.7262978\\
0.26683867368624	0.63636011\\
0.280160357243089	0.61596948\\
0.317769800918135	0.55827725\\
0.368692731322828	0.48846591\\
0.419590687014326	0.43482485\\
0.470457034601346	0.35552141\\
0.549745000555364	0.29421887\\
0.571996430244576	0.27688676\\
0.673386550433174	0.1912908\\
0.774628846383967	0.099588819\\
0.87571532435139	0.03293702\\
0.875715326563628	0.032937016\\
0.906810625077606	-0.026340958\\
0.93788744274642	-0.063612863\\
0.968944452405112	-0.1646236\\
0.999999897902016	-0.25131053\\
};
\addlegendentry{FLAT}

\addplot [color=green, line width=1.3pt, mark=square, mark options={solid, green}, forget plot]
  table[row sep=crcr]{%
4.40718485039682e-08	0.31193304\\
0.00305491401248303	-0.77936101\\
0.00998719963173467	-0.88200235\\
0.0150907237466243	-0.79928344\\
0.0182824712960732	-0.74756467\\
0.0272148896237624	-0.5940305\\
0.0455059860594568	-0.45471928\\
0.0643377785147553	-0.29267293\\
0.114392386981218	-0.17855465\\
0.165040598944605	-0.070433684\\
0.215896467105962	-0.028852796\\
0.266838109200108	0.0012231171\\
0.293534440914652	0.0094266105\\
0.317769099637719	0.016921453\\
0.368692009998319	0.023430923\\
0.419589948586958	0.031605437\\
0.470456028878461	0.024662111\\
0.560629529825513	0.013265141\\
0.57199573567152	0.011873906\\
0.673386241974256	-0.0024791709\\
0.774628684343519	-0.034465309\\
0.875715298927918	-0.065921441\\
0.906810670507207	-0.10059602\\
0.93788741795018	-0.11535577\\
0.968944528051892	-0.18812452\\
0.999999897902016	-0.25131053\\
};
\end{axis}

\node[align=center,font=\large, yshift=2em] (title) 
    at (current bounding box.north)
    {$\alpha$ = 6.76 ${^\circ}$, $\Lambda$ = 6\\ $C_p$ distribution, y/2c = 0.510};
\end{tikzpicture}%
}    
}
\\
\resizebox{0.5\linewidth}{!}{
\centerline{
%
%
\begin{tikzpicture}

\begin{axis}[%
width=4.602in,
height=3.506in,
at={(0.772in,0.473in)},
scale only axis,
xmin=-1.29805728388277e-08,
xmax=1,
xlabel style={font=\color{white!15!black}},
xlabel={\textbf{$x/c$}},
ymin=-1,
ymax=2,
ylabel style={font=\color{white!15!black}},
ylabel={\textbf{$Cp$}},
axis background/.style={fill=white},
xmajorgrids,
ymajorgrids,
legend style={legend cell align=left, align=left, draw=white!15!black}
]
\addplot [color=red, line width=1.3pt, only marks, mark=asterisk, mark options={solid, red}]
  table[row sep=crcr]{%
1.20409e-07	0.25624\\
0.001863389889	1.48418\\
0.004857393025	1.81047\\
0.008223950596	1.89986\\
0.012318558121	1.90948\\
0.015263614116	1.76441\\
0.020021401009	1.74986\\
0.027632080441	1.60476\\
0.04289041	1.39192\\
0.068570135881	1.15972\\
0.095978518416	1.0339\\
0.128414005801	0.92258\\
0.160417071441	0.8137\\
0.198066832209	0.7314\\
0.236887597521	0.6757\\
0.277818597225	0.5233\\
0.321600141604	0.53288\\
0.362855640625	0.46994\\
0.408117823281	0.42875\\
0.455230885264	0.37307\\
0.508244945569	0.33188\\
0.559842147076	0.29553\\
0.613020097936	0.25677\\
0.669537972009	0.20834\\
0.726519274321	0.13815\\
0.781807871601	0.10423\\
0.839105464729	0.06306\\
0.901681285761	0.00254\\
0.9550761984	-0.08214\\
0.999634033489	-0.13537\\
};
\addlegendentry{EXPERIMENTAL}

\addplot [color=red, line width=1.3pt, only marks, mark=asterisk, mark options={solid, red}, forget plot]
  table[row sep=crcr]{%
1.20409e-07	0.29123\\
0.001863389889	-0.66725\\
0.004857393025	-0.84922\\
0.008223950596	-0.84579\\
0.012318558121	-0.80406\\
0.015263614116	-0.7379\\
0.020021401009	-0.67568\\
0.027632080441	-0.5639\\
0.04289041	-0.41829\\
0.068570135881	-0.26539\\
0.095978518416	-0.16556\\
0.128414005801	-0.08817\\
0.160417071441	-0.05172\\
0.198066832209	0.00829\\
0.236887597521	0.02648\\
0.277818597225	0.04264\\
0.321600141604	0.04998\\
0.362855640625	0.05666\\
0.408117823281	0.0637\\
0.455230885264	0.06878\\
0.508244945569	0.04662\\
0.559842147076	0.07084\\
0.613020097936	0.06329\\
0.669537972009	0.05289\\
0.726519274321	0.02898\\
0.781807871601	0.0168\\
0.839105464729	-0.00753\\
0.901681285761	-0.03386\\
0.9550761984	-0.09044\\
0.999634033489	-0.14895\\
};
\addplot [color=blue, line width=1.3pt, mark=square, mark options={solid, blue}]
  table[row sep=crcr]{%
-1.29805728388277e-08	0.1225105\\
0.00305489171275275	1.3917456\\
0.00998717843424246	1.711938\\
0.0182825618615516	1.6479404\\
0.0272148112492127	1.4086174\\
0.045505938530895	1.3117592\\
0.0643377136508346	1.0671576\\
0.11437500317899	0.90077323\\
0.124917436841641	0.87432712\\
0.165005803927117	0.77277118\\
0.215844434221016	0.66216052\\
0.266769051008003	0.57760978\\
0.31768458214255	0.50369573\\
0.368592000415678	0.43989015\\
0.391686476628168	0.41440046\\
0.419474095979753	0.38406023\\
0.470324381984058	0.32535771\\
0.521120123243527	0.28114161\\
0.57189622299315	0.24159999\\
0.622611326442285	0.20218477\\
0.662075029509358	0.17190266\\
0.673320485435413	0.16302739\\
0.774595832673263	0.078251936\\
0.875715337278809	0.019089732\\
0.875715339491047	0.019089738\\
0.906810543582924	-0.029053876\\
0.933179436145657	-0.061121896\\
0.937887454567719	-0.067066252\\
0.968944431130516	-0.18339409\\
0.999999947243687	-0.30321276\\
};
\addlegendentry{RELAX}

\addplot [color=blue, line width=1.3pt, mark=square, mark options={solid, blue}, forget plot]
  table[row sep=crcr]{%
-1.29805728388277e-08	0.1225105\\
0.00305485618577828	-0.82717741\\
0.00998726597742339	-0.85397559\\
0.018282481099114	-0.69697934\\
0.0272148423310955	-0.53802586\\
0.0455059378961858	-0.39348644\\
0.0643377312220884	-0.23551342\\
0.0643377323282073	-0.23551342\\
0.114374833474657	-0.12713426\\
0.136225985771922	-0.082437411\\
0.165005475980829	-0.02311923\\
0.215844043492923	0.014300985\\
0.266768549627334	0.040908907\\
0.317683987874277	0.052479304\\
0.368591219054997	0.0571656\\
0.404617738255793	0.056930386\\
0.419473220873706	0.057076775\\
0.470323440971633	0.05164282\\
0.521120085575028	0.044302359\\
0.571895574170043	0.035919189\\
0.622611349774654	0.026305245\\
0.670821085696949	0.017378796\\
0.673320088828077	0.016938722\\
0.774595620116798	-0.018823585\\
0.87571527986556	-0.053757951\\
0.906810533216916	-0.098973297\\
0.935164820718243	-0.10095214\\
0.93788746839797	-0.1010255\\
0.968944323463691	-0.20257136\\
0.999999947243687	-0.30321276\\
};
\addplot [color=green, line width=1.3pt, mark=square, mark options={solid, green}]
  table[row sep=crcr]{%
-1.29805728388277e-08	0.10926241\\
0.00305489171275275	1.372264\\
0.00998717843424246	1.6945288\\
0.0182825618615516	1.633171\\
0.0272148112492127	1.3967514\\
0.045505938530895	1.3014824\\
0.0643377136508346	1.0591052\\
0.11437500317899	0.89400297\\
0.124917436841641	0.86775839\\
0.165005803927117	0.76698184\\
0.215844434221016	0.65694553\\
0.266769051008003	0.57278389\\
0.31768458214255	0.49913919\\
0.368592000415678	0.43552127\\
0.391686476628168	0.41009098\\
0.419474095979753	0.37982544\\
0.470324381984058	0.32123086\\
0.521120123243527	0.27707949\\
0.57189622299315	0.23754878\\
0.622611326442285	0.19812968\\
0.662075029509358	0.16783889\\
0.673320485435413	0.15896392\\
0.774595832673263	0.074219108\\
0.875715337278809	0.014662818\\
0.875715339491047	0.014662824\\
0.906810543582924	-0.031939402\\
0.933179436145657	-0.070615165\\
0.937887454567719	-0.077751078\\
0.968944431130516	-0.17035526\\
0.999999947243687	-0.25389919\\
};
\addlegendentry{FLAT}

\addplot [color=green, line width=1.3pt, mark=square, mark options={solid, green}, forget plot]
  table[row sep=crcr]{%
-1.29805728388277e-08	0.10926241\\
0.00305485618577828	-0.83001024\\
0.00998726597742339	-0.85195416\\
0.018282481099114	-0.69389576\\
0.0272148423310955	-0.53491032\\
0.0455059378961858	-0.39050835\\
0.0643377312220884	-0.23317975\\
0.0643377323282073	-0.23317973\\
0.114374833474657	-0.12540732\\
0.136225985771922	-0.080965362\\
0.165005475980829	-0.021983873\\
0.215844043492923	0.015014698\\
0.266768549627334	0.041275725\\
0.317683987874277	0.052554794\\
0.368591219054997	0.056984723\\
0.404617738255793	0.056588136\\
0.419473220873706	0.056668516\\
0.470323440971633	0.050999917\\
0.521120085575028	0.043482251\\
0.571895574170043	0.034920067\\
0.622611349774654	0.025121523\\
0.670821085696949	0.016017161\\
0.673320088828077	0.015567814\\
0.774595620116798	-0.020441527\\
0.87571527986556	-0.055802148\\
0.906810533216916	-0.099260151\\
0.935164820718243	-0.10840572\\
0.93788746839797	-0.10917865\\
0.968944323463691	-0.18784623\\
0.999999947243687	-0.25389919\\
};
\end{axis}

\node[align=center,font=\large, yshift=2em] (title) 
    at (current bounding box.north)
    {$\alpha$ = 6.76 ${^\circ}$, $\Lambda$ = 6\\ $C_p$ distribution, y/2c = 0.710};
\end{tikzpicture}%
}    
}
&
\resizebox{0.5\linewidth}{!}{
\centerline{
%
%
\begin{tikzpicture}

\begin{axis}[%
width=4.602in,
height=3.506in,
at={(0.772in,0.473in)},
scale only axis,
xmin=0,
xmax=1,
xlabel style={font=\color{white!15!black}},
xlabel={\textbf{$x/c$}},
ymin=-1,
ymax=2,
ylabel style={font=\color{white!15!black}},
ylabel={\textbf{$Cp$}},
axis background/.style={fill=white},
xmajorgrids,
ymajorgrids,
legend style={legend cell align=left, align=left, draw=white!15!black}
]
\addplot [color=red, line width=1.3pt, only marks, mark=asterisk, mark options={solid, red}]
  table[row sep=crcr]{%
3.6e-07	-0.01447\\
0.0018714276	1.08741\\
0.004965175296	1.46107\\
0.008151922944	1.5285\\
0.0119858704	1.55735\\
0.015486806916	1.46564\\
0.020293142116	1.47761\\
0.026051251216	1.328\\
0.043583242756	1.15896\\
0.067083590025	0.97062\\
0.096439439209	0.88837\\
0.126517510249	0.79652\\
0.160668697225	0.70225\\
0.197310974809	0.57905\\
0.236193084009	0.58849\\
0.274096225764	0.45567\\
0.317598746481	0.43378\\
0.362162036401	0.38536\\
0.410400390625	0.34177\\
0.457741552356	0.31024\\
0.509322014224	0.27871\\
0.561898660801	0.23995\\
0.614331466849	0.20362\\
0.6700732164	0.17692\\
0.727172329536	0.11888\\
0.783577269601	0.08496\\
0.844186089616	0.03657\\
0.899341168896	-0.00699\\
0.956167154244	-0.07948\\
0.996129751969	-0.14469\\
};
\addlegendentry{EXPERIMENTAL}

\addplot [color=red, line width=1.3pt, only marks, mark=asterisk, mark options={solid, red}, forget plot]
  table[row sep=crcr]{%
3.6e-07	-0.0231\\
0.0018714276	-0.75368\\
0.004965175296	-0.8308\\
0.008151922944	-0.78775\\
0.0119858704	-0.70106\\
0.015486806916	-0.63796\\
0.020293142116	-0.56238\\
0.026051251216	-0.46442\\
0.043583242756	-0.27503\\
0.067083590025	-0.14043\\
0.096439439209	-0.05982\\
0.126517510249	0.00954\\
0.160668697225	0.04607\\
0.197310974809	0.08408\\
0.236193084009	0.09257\\
0.274096225764	0.11047\\
0.317598746481	0.10588\\
0.362162036401	0.12256\\
0.410400390625	0.11256\\
0.457741552356	0.105\\
0.509322014224	0.10035\\
0.561898660801	0.10495\\
0.614331466849	0.09673\\
0.6700732164	0.07666\\
0.727172329536	0.04789\\
0.783577269601	0.02723\\
0.844186089616	0.01622\\
0.899341168896	-0.03116\\
0.956167154244	-0.09023\\
0.996129751969	-0.14546\\
};
\addplot [color=blue, line width=1.3pt, mark=square, mark options={solid, blue}]
  table[row sep=crcr]{%
1.42904951205836e-08	-0.19988538\\
0.00305015412086629	0.90724903\\
0.00997179307057428	1.281737\\
0.0182496607148598	1.2765747\\
0.0229420186305738	1.1964639\\
0.0271657936133146	1.1246673\\
0.0454834996874235	1.0056448\\
0.0643376810479438	0.8987208\\
0.08923816994892	0.79633844\\
0.114376717422473	0.72191876\\
0.139646714123514	0.6565789\\
0.15913492772325	0.61612481\\
0.165001887373523	0.6035049\\
0.19039155333604	0.55642968\\
0.215815413403402	0.5156396\\
0.241260193873826	0.47927982\\
0.266717076487603	0.4466942\\
0.292169661709477	0.41665971\\
0.292375137288655	0.41644156\\
0.317623268661543	0.38846886\\
0.34307448348037	0.36226329\\
0.368521258555115	0.33749712\\
0.393962480149634	0.31461632\\
0.419396228587831	0.29244858\\
0.42662826577855	0.28669012\\
0.444821912722126	0.27186054\\
0.470238898480765	0.25201994\\
0.495648616009797	0.23274921\\
0.521048915204115	0.21509798\\
0.546440204155708	0.19686814\\
0.561522207786905	0.18732634\\
0.57182225042861	0.18065883\\
0.597194229023843	0.16334729\\
0.622559503461322	0.14750026\\
0.647911317810328	0.13128011\\
0.673259795697713	0.11548235\\
0.696752742979599	0.10074444\\
0.723931798971521	0.083803147\\
0.77456634011886	0.052348588\\
0.825161782361485	0.014772624\\
0.875715200811608	-0.018829651\\
0.906810534633719	-0.058419924\\
0.937887351196414	-0.083595872\\
0.967872030151123	-0.17907007\\
0.968944327759211	-0.18251084\\
0.999999902986348	-0.27786329\\
};
\addlegendentry{RELAX}

\addplot [color=blue, line width=1.3pt, mark=square, mark options={solid, blue}, forget plot]
  table[row sep=crcr]{%
1.42904951205836e-08	-0.19988538\\
0.00305021148678286	-0.85730767\\
0.00997140677696673	-0.74164486\\
0.0182483518801969	-0.54776448\\
0.027163358012248	-0.36690503\\
0.0298039124727381	-0.34837255\\
0.0454820394707993	-0.24085949\\
0.0643377149595442	-0.11189257\\
0.089238666776928	-0.044810317\\
0.11437725521617	0.012443379\\
0.139647243710099	0.043130469\\
0.16500205439394	0.065466508\\
0.171392506693301	0.068853207\\
0.190391607542871	0.080387078\\
0.215815179628247	0.090654247\\
0.241259780879099	0.097168475\\
0.26671650969351	0.10080433\\
0.292168945672547	0.10254912\\
0.305759745861238	0.10254177\\
0.317622354066566	0.10252516\\
0.34307357422354	0.10123197\\
0.368520456269849	0.098966084\\
0.393961551079913	0.095844492\\
0.419395294141091	0.092348136\\
0.439209504345402	0.089072876\\
0.444820896378833	0.088181444\\
0.470237942825613	0.083819196\\
0.495647640822351	0.07911849\\
0.521048056878147	0.073957235\\
0.546439393726922	0.068860896\\
0.571821540318256	0.063073561\\
0.57219105255513	0.062992737\\
0.597194008264183	0.057450499\\
0.622558828362375	0.051193766\\
0.647911036741011	0.044848576\\
0.673259321211346	0.038124565\\
0.6985957922243	0.030689836\\
0.704747056277672	0.028834384\\
0.723931518719358	0.023281407\\
0.774566140986811	0.0068877968\\
0.825161662052223	-0.016032878\\
0.875715176494255	-0.039181113\\
0.906810488959577	-0.069051422\\
0.937887262420622	-0.090459831\\
0.968829420500343	-0.18594289\\
0.968944373628452	-0.18629766\\
0.999999902986348	-0.27786329\\
};
\addplot [color=green, line width=1.3pt, mark=square, mark options={solid, green}]
  table[row sep=crcr]{%
1.42904951205836e-08	-0.20797327\\
0.00305015412086629	0.89404917\\
0.00997179307057428	1.2696679\\
0.0182496607148598	1.2663\\
0.0229420186305738	1.1872844\\
0.0271657936133146	1.1164883\\
0.0454834996874235	0.99872071\\
0.0643376810479438	0.89293224\\
0.08923816994892	0.7912764\\
0.114376717422473	0.71738535\\
0.139646714123514	0.6524092\\
0.15913492772325	0.6121614\\
0.165001887373523	0.59960479\\
0.19039155333604	0.55272591\\
0.215815413403402	0.51208204\\
0.241260193873826	0.47583395\\
0.266717076487603	0.44332984\\
0.292169661709477	0.41335934\\
0.292375137288655	0.41314149\\
0.317623268661543	0.38521692\\
0.34307448348037	0.35904559\\
0.368521258555115	0.33430368\\
0.393962480149634	0.31143072\\
0.419396228587831	0.28926796\\
0.42662826577855	0.28350538\\
0.444821912722126	0.26866633\\
0.470238898480765	0.24880478\\
0.495648616009797	0.22950783\\
0.521048915204115	0.21180908\\
0.546440204155708	0.19354093\\
0.561522207786905	0.1839561\\
0.57182225042861	0.17726299\\
0.597194229023843	0.15990275\\
0.622559503461322	0.14398007\\
0.647911317810328	0.12768672\\
0.673259795697713	0.1118053\\
0.696752742979599	0.096977919\\
0.723931798971521	0.079927869\\
0.77456634011886	0.048235547\\
0.825161782361485	0.010440264\\
0.875715200811608	-0.023534276\\
0.906810534633719	-0.062375899\\
0.937887351196414	-0.091704354\\
0.967872030151123	-0.17488509\\
0.968944327759211	-0.17786971\\
0.999999902986348	-0.25332272\\
};
\addlegendentry{FLAT}

\addplot [color=green, line width=1.3pt, mark=square, mark options={solid, green}, forget plot]
  table[row sep=crcr]{%
1.42904951205836e-08	-0.20797327\\
0.00305021148678286	-0.857804\\
0.00997140677696673	-0.73898333\\
0.0182483518801969	-0.54462224\\
0.027163358012248	-0.36397964\\
0.0298039124727381	-0.34550285\\
0.0454820394707993	-0.23831144\\
0.0643377149595442	-0.10979334\\
0.089238666776928	-0.043102156\\
0.11437725521617	0.01382008\\
0.139647243710099	0.044247821\\
0.16500205439394	0.066369243\\
0.171392506693301	0.069710784\\
0.190391607542871	0.081100725\\
0.215815179628247	0.091200493\\
0.241259780879099	0.097565047\\
0.26671650969351	0.10106458\\
0.292168945672547	0.10268182\\
0.305759745861238	0.10261046\\
0.317622354066566	0.10253784\\
0.34307357422354	0.10112983\\
0.368520456269849	0.098752014\\
0.393961551079913	0.095519863\\
0.419395294141091	0.091912806\\
0.439209504345402	0.088549942\\
0.444820896378833	0.087633967\\
0.470237942825613	0.08315935\\
0.495647640822351	0.07834363\\
0.521048056878147	0.073061533\\
0.546439393726922	0.067843705\\
0.571821540318256	0.061921664\\
0.57219105255513	0.061838921\\
0.597194008264183	0.056169044\\
0.622558828362375	0.049766842\\
0.647911036741011	0.043275144\\
0.673259321211346	0.036396064\\
0.6985957922243	0.028795246\\
0.704747056277672	0.026899513\\
0.723931518719358	0.02122839\\
0.774566140986811	0.0044924454\\
0.825161662052223	-0.018752402\\
0.875715176494255	-0.042305477\\
0.906810488959577	-0.071495578\\
0.937887262420622	-0.096626416\\
0.968829420500343	-0.17989706\\
0.968944373628452	-0.18020512\\
0.999999902986348	-0.25332272\\
};
\end{axis}

\node[align=center,font=\large, yshift=2em] (title) 
    at (current bounding box.north)
    {$\alpha$ = 6.76 ${^\circ}$, $\Lambda$ = 6\\ $C_p$ distribution, y/2c = 0.840};
\end{tikzpicture}%
}    
}
\end{tabular}  	
\caption{Comparative analysis between results obtained using $\Pi$-BEM and experimental data for 
a swept wing with a sweep angle of 20${^\circ}$. The plots refer to four planar sections normal 
to the wing axis, and located at coordinates $y/s$ =  0.260, 0.510, 0.710, 0.840, as measured
at their intersection with the trailing edge.\label{fig:ArrowWing}}
\end{figure}

Figure \ref{fig:ArrowWing} presents a comparison between experimental and numerical values
of $-c_p$ on four cross-sections perpendicular to the wing axis direction. We point out in fact that because in the experimental
setup the swept wing has been obtained merely rotating the NACA 0012 wing model, all the pressure taps in this test are not here
aligned with the incident flow, but are instead inclined of an angle $\beta=20^\circ$ with respect to such a direction.
The wing sections are located at coordinates $y/s$ =  0.260, 0.510, 0.710, 0.840, as measured at their intersection with the trailing edge.
In each section, also in this case the numerical values of windward and leeward pressure coefficients obtained imposing a
flat wake (blue lines) or computing the relaxed wake (green lines) are plotted against their experimental counterparts (red dots).

The results appear quite consistent with the ones obtained on the rectangular wing geometry.
The comparisons with the experimental results appear in fact once again satisfactory, by a qualitative perspective.
The experimental curves are in fact reproduced in all their main traits by the numerical results.
By a more quantitative point of view, also in this case the peaks location appear accurately reproduced by the numerical results.
As for the peaks height, on one hand the windward side pressure peak is captured with good accuracy also in this swept wing case, in which
a stagnation point is missing and the $c_p$ maximum is approximately $0.9$. On the other hand, the peaks on the suction side are slightly
underestimated, as was the case for the rectangular wing. And, also in the present case, the discrepancy seems more evident for the
sections closer to the tip, likely due to viscous effects associated with the tip vortex. Again, we can observe that the sharp
trailing edge considered in the numerical simulations --- as opposed to the experimental finite thickness one --- results in a
noticeable difference in correspondence with the last pressure probe measurement, at every section considered. 
Finally, we remark that also in the present case the effect of wake relaxation appears marginal, except on the trailing
edge region, in which the relaxed wake solution predicts higher pressure recovery.

\section{Conclusions and future perspectives}\label{sec:conclusions}

This work presented a numerical method for the simulation of quasi-potential flow past lifting surfaces. The solver
combines the collocation Boundary Element Method formulation implemented in the open source library $\pi$-BEM \cite{pi-BEM} with
a Galerkin formulation of the nonkinear Kutta condition imposed at the lifting body trailing edge. In addition, Galerkin formulation
is used also for the Hypersingular Boundary Integral Equation solved to obtain the velocities at the wake collocation points needed
for wake relaxation. The latter techniques, borrowed from the Finite Element Method, allow for the evaluation of the solution derivatives
in a way that is independent of the local grid topology. As such, the proposed Boundary Element Method based solver allows for the user selection
of the finite element degree at the start of each simulation. This is also made possible by the direct interface with CAD files, which
are interrogated to obtain new points for high order elements grids on curved surfaces. Numerical experiments on rectangular and V-shaped wings
with NACA 0012 airfoil sections reproduce in a satisfactory way experimental data. The results on rectangular wings also confirm the quality
of the arbitrary order numerical approach implemented. In particular, the higher order solutions show
less numerical dissipation in the leading edge region, which is characterized by the highest gradients. However, in BEM the polynomial shape
functions are multiplied by singular Green functions so as to obtain the linear system matrices and vectors. Thus, because high order shape functions
change sign multiple times within each cell, increased accuracy granted by high order elements can only be obtained at the price of a significant
increase in the Gaussian quadrature nodes used for the BEM integrals. For such a reason, the combination of linear elements and additional local
refinements appears to better balance accuracy and computational cost.
In the near future, a possible way to best exploit the accuracy premium that higher order elements offer is that of implementing suitable
$hp$ refinement algorithms. Indeed, this paper showed that the methodology implemented possesses all the ingredients necessary for the implementation
of such a strategy. Another interesting possible extension of this work, is that of using the CAD interface for the mapping of the
reference cell to three dimensional cells in the calculation of the integrals involved in BEM. This approach, which is also discussed
in \cite{heltaiEtAlACMTOMS2021}, has never been applied to BEM. As isogeometric analysis, such a \emph{heterogeometric} approach would allow for
the calculation of integrals on the exact surface prescribed by the user, but without the need to use NURBS shape functions also for the solution.
Finally, further developments of this work can come from an extension of the range of application of the model.
In particular, the free of a surface in proximity of the wing can be initially introduced as a homogeneous Neumann boundary condition,
and can be then treated with a linearized \cite{Giuliani2015} or fully nonlinear \cite{MOLA2023322} formulation. Including  
the presence of a free surface in close to the wing would make the solver able to estimate the forces generated by
hydrofoils in different relative position with respect to the free surface, and contribute to the design of such artifacts.

\section{Acknowledgments}

We acknowledge the support of Piano Nazionale di Ripresa e Resilienza (PNRR) and ToolsPole for co-funding the PhD scholarship
of Luca Cattarossi. Andrea Mola acknowledges project ROMEU: Reduced Order Models for Environmental and Urban flow, (CUP D53D230188800) funded by
Ministero Università e Ricerca (MUR), Italy.  We finally thank Istituto Nazionale di Alta Matematica "Francesco Severi" (INdAM), Gruppo Nazionale
per il Calcolo Scientifico (GNCS) for the support.

\FloatBarrier

\bibliography{mybibfile}

\begin{thebibliography}{10}
\expandafter\ifx\csname url\endcsname\relax
  \def\url#1{\texttt{#1}}\fi
\expandafter\ifx\csname urlprefix\endcsname\relax\def\urlprefix{URL }\fi
\expandafter\ifx\csname href\endcsname\relax
  \def\href#1#2{#2} \def\path#1{#1}\fi

\bibitem{PhillipsEtAl2001}
W.~F. Phillips, D.~O. Snyder, \href{https://doi.org/10.2514/2.2649}{Modern
  adaptation of prandtl's classic lifting-line theory}, Journal of Aircraft
  37~(4) (2000) 662--670.
\newblock \href {http://arxiv.org/abs/https://doi.org/10.2514/2.2649}
  {\path{arXiv:https://doi.org/10.2514/2.2649}}, \href
  {https://doi.org/10.2514/2.2649} {\path{doi:10.2514/2.2649}}.
\newline\urlprefix\url{https://doi.org/10.2514/2.2649}

\bibitem{LiuEtAl1996}
D.~D. Liu, P.~C. Chen, Z.~X. Yao, D.~Sarhaddi, Recent advances in lifting
  surface methods, The Aeronautical Journal 100~(998) (1996) 327–340.
\newblock \href {https://doi.org/10.1017/S0001924000067038}
  {\path{doi:10.1017/S0001924000067038}}.

\bibitem{HESS19671}
J.~Hess, A.~Smith,
  \href{https://www.sciencedirect.com/science/article/pii/0376042167900036}{Calculation
  of potential flow about arbitrary bodies}, Progress in Aerospace Sciences 8
  (1967) 1--138.
\newblock \href {https://doi.org/https://doi.org/10.1016/0376-0421(67)90003-6}
  {\path{doi:https://doi.org/10.1016/0376-0421(67)90003-6}}.
\newline\urlprefix\url{https://www.sciencedirect.com/science/article/pii/0376042167900036}

\bibitem{pi-BEM}
N.~Giuliani, A.~Mola, L.~Heltai,
  \href{http://www.sciencedirect.com/science/article/pii/S0965997818300371}{$\pi$-bem:
  A flexible parallel implementation for adaptive, geometry aware, and high
  order boundary element methods}, Advances in Engineering Software 121 (2018)
  39 -- 58.
\newblock \href
  {https://doi.org/https://doi.org/10.1016/j.advengsoft.2018.03.008}
  {\path{doi:https://doi.org/10.1016/j.advengsoft.2018.03.008}}.
\newline\urlprefix\url{http://www.sciencedirect.com/science/article/pii/S0965997818300371}

\bibitem{gaggeroBrizzolara2007}
S.~Gaggero, S.~Brizzolara, {Exact Modeling of Trailing Vorticity in Panel
  Method for Marine Propeller}, 2nd International Conference on Marine Research
  and Transportation, 2007, ischia, 28-30 June.

\bibitem{heltaiEtAlACMTOMS2021}
L.~Heltai, W.~Bangerth, M.~Kronbichler, A.~Mola,
  \href{https://doi.org/10.1145/3468428}{Propagating geometry information to
  finite element computations}, ACM Trans. Math. Softw. 47~(4) (sep 2021).
\newblock \href {https://doi.org/10.1145/3468428} {\path{doi:10.1145/3468428}}.
\newline\urlprefix\url{https://doi.org/10.1145/3468428}

\bibitem{morinoGennaretti1992}
\href{https://arc.aiaa.org/doi/abs/10.2514/5.9781600866180.0279.0320}{Boundary
  Integral Equation Methods for Aerodynamics}, pp. 279--320.
\newblock \href
  {http://arxiv.org/abs/https://arc.aiaa.org/doi/pdf/10.2514/5.9781600866180.0279.0320}
  {\path{arXiv:https://arc.aiaa.org/doi/pdf/10.2514/5.9781600866180.0279.0320}},
  \href {https://doi.org/10.2514/5.9781600866180.0279.0320}
  {\path{doi:10.2514/5.9781600866180.0279.0320}}.
\newline\urlprefix\url{https://arc.aiaa.org/doi/abs/10.2514/5.9781600866180.0279.0320}

\bibitem{GENNARETTI1996281}
M.~Gennaretti, F.~Salvatore, L.~Morino,
  \href{https://www.sciencedirect.com/science/article/pii/S0889974696900171}{Forces
  and moments in incompressible quasi-potential flows}, Journal of Fluids and
  Structures 10~(3) (1996) 281--303.
\newblock \href {https://doi.org/https://doi.org/10.1006/jfls.1996.0017}
  {\path{doi:https://doi.org/10.1006/jfls.1996.0017}}.
\newline\urlprefix\url{https://www.sciencedirect.com/science/article/pii/S0889974696900171}

\bibitem{XU1998415}
C.~Xu,
  \href{https://www.sciencedirect.com/science/article/pii/S0093641398000548}{Kutta
  condition for sharp edge flows}, Mechanics Research Communications 25~(4)
  (1998) 415--420.
\newblock \href {https://doi.org/https://doi.org/10.1016/S0093-6413(98)00054-8}
  {\path{doi:https://doi.org/10.1016/S0093-6413(98)00054-8}}.
\newline\urlprefix\url{https://www.sciencedirect.com/science/article/pii/S0093641398000548}

\bibitem{CHOULIARAS2021113556}
S.~Chouliaras, P.~Kaklis, K.~Kostas, A.~Ginnis, C.~Politis,
  \href{https://www.sciencedirect.com/science/article/pii/S0045782520307416}{An
  isogeometric boundary element method for 3d lifting flows using t-splines},
  Computer Methods in Applied Mechanics and Engineering 373 (2021) 113556.
\newblock \href {https://doi.org/https://doi.org/10.1016/j.cma.2020.113556}
  {\path{doi:https://doi.org/10.1016/j.cma.2020.113556}}.
\newline\urlprefix\url{https://www.sciencedirect.com/science/article/pii/S0045782520307416}

\bibitem{waveBem}
A.~Mola, L.~Heltai, A.~DeSimone,
  \href{http://www.sciencedirect.com/science/article/pii/S0955799712001907
  http://linkinghub.elsevier.com/retrieve/pii/S0955799712001907}{{A stable and
  adaptive semi-Lagrangian potential model for unsteady and nonlinear ship-wave
  interactions}}, Eng. Anal. Bound. Elem. 37~(1) (2013) 128--143.
\newblock \href {https://doi.org/10.1016/j.enganabound.2012.09.005}
  {\path{doi:10.1016/j.enganabound.2012.09.005}}.
\newline\urlprefix\url{http://www.sciencedirect.com/science/article/pii/S0955799712001907
  http://linkinghub.elsevier.com/retrieve/pii/S0955799712001907}

\bibitem{molaIsope2016}
A.~Mola, L.~Heltai, A.~De~Simone, {Ship Sinkage and Trim Predictions Based on a
  CAD Interfaced Fully Nonlinear Potential Model}, Vol. All Days of
  International Ocean and Polar Engineering Conference, 2016, iSOPE-I-16-438.
\newblock \href
  {http://arxiv.org/abs/https://onepetro.org/ISOPEIOPEC/proceedings-pdf/ISOPE16/All-ISOPE16/ISOPE-I-16-438/1337253/isope-i-16-438.pdf}
  {\path{arXiv:https://onepetro.org/ISOPEIOPEC/proceedings-pdf/ISOPE16/All-ISOPE16/ISOPE-I-16-438/1337253/isope-i-16-438.pdf}}.

\bibitem{DASSI2014253}
F.~Dassi, A.~Mola, H.~Si,
  \href{https://www.sciencedirect.com/science/article/pii/S1877705814016671}{Curvature-adapted
  remeshing of cad surfaces}, Procedia Engineering 82 (2014) 253--265, 23rd
  International Meshing Roundtable (IMR23).
\newblock \href {https://doi.org/https://doi.org/10.1016/j.proeng.2014.10.388}
  {\path{doi:https://doi.org/10.1016/j.proeng.2014.10.388}}.
\newline\urlprefix\url{https://www.sciencedirect.com/science/article/pii/S1877705814016671}

\bibitem{bernasconiPhD}
D.~J. Bernasconi, \href{https://eprints.soton.ac.uk/466428/}{A higher-order
  potential flow method for thick bodies, thin surfaces and wakes}, Ph.D.
  thesis, University of Southampton (2007).
\newline\urlprefix\url{https://eprints.soton.ac.uk/466428/}

\bibitem{guiggianiEtAl1992}
M.~Guiggiani, G.~Krishnasamy, T.~J. Rudolphi, F.~J. Rizzo,
  \href{https://doi.org/10.1115/1.2893766}{{A General Algorithm for the
  Numerical Solution of Hypersingular Boundary Integral Equations}}, Journal of
  Applied Mechanics 59~(3) (1992) 604--614.
\newblock \href
  {http://arxiv.org/abs/https://asmedigitalcollection.asme.org/appliedmechanics/article-pdf/59/3/604/5462606/604\_1.pdf}
  {\path{arXiv:https://asmedigitalcollection.asme.org/appliedmechanics/article-pdf/59/3/604/5462606/604\_1.pdf}},
  \href {https://doi.org/10.1115/1.2893766} {\path{doi:10.1115/1.2893766}}.
\newline\urlprefix\url{https://doi.org/10.1115/1.2893766}

\bibitem{GRIMBERG20081878}
G.~Grimberg, W.~Pauls, U.~Frisch,
  \href{https://www.sciencedirect.com/science/article/pii/S0167278908000225}{Genesis
  of d’alembert’s paradox and analytical elaboration of the drag problem},
  Physica D: Nonlinear Phenomena 237~(14) (2008) 1878--1886, euler Equations:
  250 Years On.
\newblock \href {https://doi.org/https://doi.org/10.1016/j.physd.2008.01.015}
  {\path{doi:https://doi.org/10.1016/j.physd.2008.01.015}}.
\newline\urlprefix\url{https://www.sciencedirect.com/science/article/pii/S0167278908000225}

\bibitem{MORINO2001805}
L.~Morino, G.~Bernardini,
  \href{https://www.sciencedirect.com/science/article/pii/S0955799701000637}{Singularities
  in bies for the laplace equation; joukowski trailing-edge conjecture
  revisited}, Engineering Analysis with Boundary Elements 25~(9) (2001)
  805--818.
\newblock \href {https://doi.org/https://doi.org/10.1016/S0955-7997(01)00063-7}
  {\path{doi:https://doi.org/10.1016/S0955-7997(01)00063-7}}.
\newline\urlprefix\url{https://www.sciencedirect.com/science/article/pii/S0955799701000637}

\bibitem{Giuliani2015}
N.~Giuliani, A.~Mola, L.~Heltai, L.~Formaggia,
  \href{http://linkinghub.elsevier.com/retrieve/pii/S0955799715001058}{{FEM
  SUPG stabilisation of mixed isoparametric BEMs : Application to linearised
  free surface flows}}, Engineering Analysis with Boundary Elements 59 (2015)
  8--22.
\newblock \href {https://doi.org/10.1016/j.enganabound.2015.04.006}
  {\path{doi:10.1016/j.enganabound.2015.04.006}}.
\newline\urlprefix\url{http://linkinghub.elsevier.com/retrieve/pii/S0955799715001058}

\bibitem{brebbia}
C.~A. Brebbia, {The Boundary Element Method for Engineers}, Pentech Press,
  1978.

\bibitem{pi-BEM-repo}
N.~Giuliani, A.~Mola, L.~Heltai, pi-{BEM}: {P}arallel {BEM} {S}olver.,
  \url{https://github.com/mathLab/pi-BEM} (2021).

\bibitem{dealII92}
D.~Arndt, W.~B.~M. Feder, M.~Fehling, R.~Gassm{\"o}ller, T.~Heister, L.~Heltai,
  M.~Kronbichler, M.~Maier, P.~Munch, J.-P. Pelteret, S.~Sticko, B.~Turcksin,
  D.~Wells, \href{https://dealii.org/deal94-preprint.pdf}{The \texttt{deal.II}
  library, version 9.4}, Journal of Numerical MathematicsAccepted (2022).
\newline\urlprefix\url{https://dealii.org/deal94-preprint.pdf}

\bibitem{BangerthHeisterHeltai-2016-b}
W.~Bangerth, D.~Davydov, T.~Heister, L.~Heltai, G.~Kanschat, M.~Kronbichler,
  M.~Maier, B.~Turcksin, D.~Wells, {The deal.II library, Version 8.4}, Journal
  of Numerical Mathematics 24~(3) (2016) 135--141.

\bibitem{ArndtBangerthDavydov-2017-a}
D.~Arndt, W.~Bangerth, D.~Davydov, T.~Heister, L.~Heltai, M.~Kronbichler,
  M.~Maier, J.-P. Pelteret, B.~Turcksin, D.~Wells, The \texttt{deal.II}
  library, version 8.5, Journal of Numerical Mathematics (2017).

\bibitem{Mola2014}
A.~Mola, L.~Heltai, A.~DeSimone,
  \href{http://digitallibrary.sissa.it/handle/1963/7311}{{A fully nonlinear
  potential model for ship hydrodynamics directly interfaced with CAD data
  structures}}, in: 24th International Ocean and Polar Engineering Conference,
  2014.
\newline\urlprefix\url{http://digitallibrary.sissa.it/handle/1963/7311}

\bibitem{MolaHeltaiDeSimone-2017-a}
A.~Mola, L.~Heltai, A.~DeSimone, Wet and dry transom stern treatment for fully
  nonlinear potential flow simulations of naval hydrodynamics, Journal of Ship
  Research 61~(1) (2017) 1--14.

\bibitem{GiulianiPoliMIThesis}
N.~Giuliani, An hybrid boundary element method for free surface flows, Master's
  thesis, Politecnico di Milano (2013).

\bibitem{grilli2001}
S.~T. Grilli, P.~Guyenne, F.~Dias,
  \href{https://onlinelibrary.wiley.com/doi/abs/10.1002/1097-0363(20010415)35:7%3C829::AID-FLD115%3E3.0.CO;2-2}{{A
  fully non-linear model for three-dimensional overturning waves over an
  arbitrary bottom}}, Journal for Numerical Methods in Fluids 35~(7) (2015)
  29--67.
\newblock \href
  {https://doi.org/https://doi.org/10.1002/1097-0363(20010415)35:7<829::AID-FLD115>3.0.CO;2-2}
  {\path{doi:https://doi.org/10.1002/1097-0363(20010415)35:7<829::AID-FLD115>3.0.CO;2-2}}.
\newline\urlprefix\url{https://onlinelibrary.wiley.com/doi/abs/10.1002/1097-0363(20010415)35:7%3C829::AID-FLD115%3E3.0.CO;2-2}

\bibitem{Mousavi2010}
S.~E. Mousavi, N.~Sukumar, {Generalized Duffy transformation for integrating
  vertex singularities}, Computational Mechanics 45~(2-3) (2010) 127--140.
\newblock \href {https://doi.org/10.1007/s00466-009-0424-1}
  {\path{doi:10.1007/s00466-009-0424-1}}.

\bibitem{YipShubertNACA0012}
L.~Yip, G.~Shubert, Pressure distributions on a 1- by 3-meter semispan wing at
  sweep angles from 0$^\circ$ to 40$^\circ$ in subsonic flow, Tech. Rep. TN
  D-8307, NASA (1976).

\bibitem{MOLA2023322}
A.~Mola, N.~Giuliani, Óscar Crego, G.~Rozza,
  \href{https://www.sciencedirect.com/science/article/pii/S0307904X23002445}{A
  unified steady and unsteady formulation for hydrodynamic potential flow
  simulations with fully nonlinear free surface boundary conditions}, Applied
  Mathematical Modelling 122 (2023) 322--349.
\newblock \href {https://doi.org/https://doi.org/10.1016/j.apm.2023.06.001}
  {\path{doi:https://doi.org/10.1016/j.apm.2023.06.001}}.
\newline\urlprefix\url{https://www.sciencedirect.com/science/article/pii/S0307904X23002445}

\end{thebibliography}

\end{document}